\documentclass[preprint,11pt]{elsarticle2}
\usepackage{amsmath,amssymb,amsthm,mathrsfs,geometry,graphicx,setspace,soul,booktabs,dsfont}
\allowdisplaybreaks

\usepackage{cases} 
\usepackage{accents,comment}
\usepackage{color}

\usepackage{hyperref}
\usepackage{mathtools} 
\usepackage{upgreek,mleftright}
\usepackage[greek,UKenglish]{babel}
\usepackage{teubner} 
\usepackage{comment} 
\usepackage{caption}
\usepackage{subcaption} 
\usepackage{framed} 
\usepackage{multicol} 
\usepackage[linesnumbered,ruled,vlined]{algorithm2e}

\SetCommentSty{mycommfont}
\usepackage{parskip}

\usepackage{tikz}
\usetikzlibrary{fit,calc} %
\usetikzlibrary{calc,intersections,patterns}
\usepackage{graphicx}

\newcommand*{\tikzmk}[1]{\tikz[remember picture,overlay,] \node (#1) {};\ignorespaces}
\newcommand{\boxit}[1]{\tikz[remember picture,overlay]{\node[yshift=3pt,fill=#1,opacity=.15,fit={(A)($(B)+(.125\linewidth,-.8\baselineskip)$)}] {};}\ignorespaces}

\usepackage[scr=boondox,  
            cal=esstix]   
           {mathalpha}
\usepackage{multirow}
\DeclareMathAlphabet\mathbfcal{OMS}{cmsy}{b}{n} 
\usepackage[table,dvipsnames]{xcolor} 

\usepackage{threeparttable} 
\newcommand{\headrow}{\rowcolor{black!20}} 
\biboptions{sort&compress}
\usepackage{framed} 
\usepackage{nomencl} 
\makenomenclature
\renewcommand*\nompreamble{\begin{multicols}{2}}
\renewcommand*\nompostamble{\end{multicols}}

\newcommand{\splitcellc}[2][c]{%
  \begin{tabular}[c]{@{}c@{}}\strut#2\strut\end{tabular}}%

\usepackage{amssymb}

 \usepackage{subcaption}
 
\usepackage[export]{adjustbox} 
\definecolor{clr1}{RGB}{184,187,164}
\usepackage{float} 
\renewcommand{\arraystretch}{1.3}

\journal{J. Mech. Phys. Solids}

\begin{document}
\begin{frontmatter}

 \title{Caveats on formulating finite elasto-plasticity in curvilinear coordinates}

\author[du,uoft]{Giuliano Pretti}
\author[du]{Robert E. Bird}
\author[du]{William M. Coombs}
\author[du]{Charles E. Augarde \texorpdfstring{\corref{cor1}}}

\cortext[cor1]{Corresponding author:
charles.augarde@durham.ac.uk}
\address[du]{Department of Engineering, Durham University\\
 Science Site, South Road, Durham, DH1 3LE, United Kingdom.}
\address[uoft]{Department of Civil \& Mineral Engineering, University of Toronto\\
M5S 1A4 Toronto, ON, Canada.
}

\begin{abstract}
Tensor analysis provides a frame‑invariant foundation for continuum mechanics, yet numerical implementations rely on matrix representations expressed in user‑selected bases. When these bases are non‑Cartesian and non‑orthonormal, additional terms arise that are normally absent or implicit in Cartesian formulations. Using cavity expansion as an initial model problem, this paper clarifies the roles of the deformation gradient, Jacobian, and shifter in finite‑strain kinematics under axisymmetry. These quantities, typically straightforward in Cartesian frames, require more careful treatment in curvilinear coordinates, particularly in applications involving large deformations where axisymmetric reductions provide substantial computational savings.
The formulation is further complicated when anelastic effects are included: the elastic and anelastic components of the deformation gradient and Jacobian must be distinguished, and the Cauchy stress depends on configuration changes beyond the current elastic state. This increases the complexity of the consistent linearisation required for finite element implementation. The paper provides a clear, step‑by‑step procedure for handling these contributions.
The focus is practical rather than geometric. Instead of adopting a differential‑geometric manifold framework, we work within standard Cartesian representations and use explicit changes of basis to obtain the required curvilinear forms. The resulting methodology enables robust finite element analysis of axisymmetric elasto‑plastic problems undergoing finite strains.
\end{abstract}


\end{frontmatter}

\section{Introduction}
\label{sec:introduction}

\noindent Three, two-dimensional spatial approximations of the three-dimensional world are commonly used in solid mechanics problems: plane strain, plane stress and axisymmetry.  All serve to reduce the computational costs of discretised numerical methods for solids such as Finite Elements (FEs), used routinely to solve problems and to make predictions of behaviour. The literature highlights a wide range of applications that address problems distinguished by axisymmetry, emphasising their significance and relevance across various fields of engineering. Key infrastructure in geotechnical engineering, such as piled foundations~\cite{said2009axisymmetric} and tunnel excavation~\cite{zhao2015computational}, and tests such as the Cone Penetration Test (CPT)~\cite{bird2024implicit} or the Free-Fall Cone Penetrometer (FFCP)~\cite{mohapatra2025laboratory} exhibit distinctive axisymmetry that is often leveraged in the design phase for the former and in the data interpretation for the latter.  
Axisymmetric problems also occur in civil engineering, particularly in designing structures like cooling towers~\cite{williamson2008comparison}, silos~\cite{martinez2002simulation}, water tanks, and reservoirs~\cite{yoshida2010axisymmetric}. 
In the realm of metal forming, applications with this symmetry include the upsetting of cylinders~\cite{lee1971analyses}, the extrusion of billets~\cite{oyinbo2015numerical}, and wire drawing~\cite{mcallen2007numerical}. 
Moreover, in biomechanics, axisymmetric modelling is applied to the study of arterial growth~\cite{kumar2025nonlinear} and the optical-mechanical functioning of a crystalline lens~\cite{de2023estimation}.
Expanding the scope to bi-phasic materials, i.e., porous solids containing an interstitial fluid---further increases the range of applications, including fluid injections into boreholes~\cite{seth1968transient}, blood flow through the permeable walls of the vascular network~\cite{kenyon1979mathematical}, and drug injection~\cite{stiles2006two}. Regardless of their specific context, many of these mentioned applications involve some inelastic behaviour, either in the form of plasticity, growth, or thermal expansion/contraction as well as geometrical nonlinearity, i.e., the need for finite measures of strain (large deformations). 
When aiming to create FE programs that can incorporate large deformation elasto-plastic behaviour while maintaining a relatively low computational cost by exploiting the axisymmetry,  additional complexities may arise. 
These challenges are not immediately apparent when using Cartesian basis vectors, and stem primarily from the fact that continuum mechanics is formulated with tensors that are thought to be coordinate-invariant. 
However, FE programs rely on matrices that express the components of these tensors in user-defined bases, and frame invariance must be carefully embedded in the programs.

Starting from the basic example of a cavity expansion problem, this paper exposes the definition and role of quantities such as the deformation gradient, the Jacobian tracking the volume change, and the shifter. 
While these quantities are always present in any formulation, they either are well-known or play a silent role when employing Cartesian basis vectors. 
This work exposes how, if naively extending the use of Cartesian coordinates to a curvilinear system, errors arise in the definition of these quantities, which may lead to fundamental issues with FEM programs.
On top of the difficulty due to the use of curvilinear basis vectors, finite-strain elasto-plasticity, despite having been developed consistently for almost half a century, presents some aspects that are still a \emph{``work in progress''}\cite{davoli2015critical}. 
Among these factors, incorporating isochoric plastic deformations makes the Cauchy stress dependent on quantities other than the current elastic deformations~\cite{yavari2023direct,abe2024reconstruction}. 
This work explores the different definitions of the elastic and plastic Jacobians to track reversible and irreversible volume changes, and shows how to linearise the resulting equilibrium equations under these conditions.

It is important to note that this work does not make use of Riemannian manifolds and the geometric view of continuum mechanics. This choice is motivated not by the importance of these topics but by the manuscript's focus. 
This manuscript emphasises \emph{how} to develop an FE code (or indeed any similar discretised method) for axisymmetric analyses, rather than exploring \emph{why} manifold analysis provides more insights that classical tensor analysis may not capture. This goal also motivates the approach employed in this paper: the Cartesian coordinate system and basis vectors are employed to represent the considered quantities, and change of basis vectors is adopted to compute the corresponding curvilinear representations.

This manuscript adopts the same notation used in~\cite{dvorkin2006nonlinear}, with scalars denoted by italic light-face letters, e.g., $p,P$, Italic bold-face letters indicating vectors, e.g., $\underline{\boldsymbol{x}}$ and tensors, e.g., $\underline{\underline{\boldsymbol{\sigma}}}$, contravariant components of these quantities being indicated by superscripts and covariant by subscripts. The Einstein summation convention between these quantities is always implied.

\section{A cavity expansion problem}

\label{sec:cavity_expansion}
\begin{figure}
\centering
\includegraphics[width=0.65\textwidth,trim= 0cm 0 0 0]{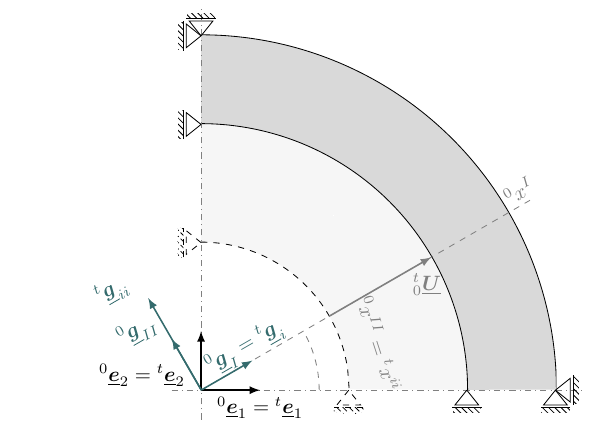}
\caption{Cross section of the cavity expansion problem and relative quantities of interest. The supports shown permit rotational and translation in one direction.}
\label{fig:cavity_expansion}
\end{figure}


\noindent Consider the cavity expansion problem of a hollow cylinder represented in Fig.~\ref{fig:cavity_expansion}, where the inner surface of the cylinder is radially expanded. To describe this example, the matching collinear (between original and current configuration) right-handed Cartesian covariant basis vectors are denoted as ${^0\boldsymbol{\underline{e}}_J}$ and ${^t\boldsymbol{\underline{e}}_j}$, with the components of the position vectors being
\begin{equation}
\left[ {^{0}z^1}, {^{0}z^{2}}, {^{0}z^{3}} \right]^T
 =
\left[ X, Y, Z \right]^T;
\qquad
\left[ {^{t}z^1}, {^{t}z^{2}}, {^{t}z^{3}} \right]^T
 =
\left[ x, y, z \right]^T.
\end{equation}
At the same time, the non-matching collinear right-handed curvilinear (cylindrical) covariant basis vectors are denoted by ${^0\boldsymbol{\underline{g}}_A}$ for the original configuration and by ${^t\boldsymbol{\underline{g}}_a}$ for the current one, with the corresponding components of the position vector labelled as
\begin{equation}
\left[ {^{0}x^I}, {^{0}x^{II}}, {^{0}x^{III}} \right]^T
 =
\left[ R, \Theta, Z \right]^T;
\qquad
\left[{^{t}x^i}, {^{t}x^{ii}}, {^{t}x^{iii}} \right]^T
 =
\left[ r, \theta, z \right]^T.
\end{equation}
The standard geometric relationships between Cartesian and cylindrical coordinates are as follows
\begin{equation}
\begin{dcases}
\label{eq:Cartesian_components_function_cylindrical}
{^{0}z^1} = {^{0}x^I} \cos \left( {^{0}x^{II}} \right);
\\
{^{0}z^2} = {^{0}x^I} \sin \left( {^{0}x^{II}} \right);
\\
{^{0}z^3} = {^{0}x^{III}};
\end{dcases}
\qquad
\begin{cases}
{^{t}z^1} = {^{t}x^i} \cos \left( {^{t}x^{ii}} \right);
\\
{^{t}z^2} = {^{t}x^i} \sin \left( {^{t}x^{ii}} \right);
\\
{^{t}z^3} = {^{t}x^{iii}}.
\end{cases}
\end{equation}
Similarly, their inverse relationships are
\begin{equation}
\begin{dcases}
{^{0}x^I} = \sqrt{ \left({^{0}z^1} \right)^2 + \left({^{0}z^2} \right)^2};
\\
{^{0}x^{II}} = \arctan \left( \dfrac{{^{0}z^2}}{{^{0}z^1}} \right);
\\
{^{0}x^{III}} = {^{0}z^3};
\end{dcases}
\qquad
\begin{dcases}
{^{t}x^i} = \sqrt{ \left({^{t}z^1} \right)^2 + \left({^{t}z^2} \right)^2};
\\
{^{t}x^{ii}} = \arctan \left( \dfrac{{^{t}z^2}}{{^{t}z^1}} \right);
\\
{^{t}x^{iii}} = {^{t}z^3}.
\end{dcases}
\end{equation}

\begin{table}
\centering
\footnotesize{
\begin{threeparttable} 
\setlength\minrowclearance{2.4pt}
\setlength{\tabcolsep}{6pt} 
\begin{tabular}{c|c|c|c|c} 
\headrow
\textbf{\scriptsize{Configuration}} 
& 
{\splitcellc[]{\textbf{Covariant and} \\
\textbf{contravariant} \\ 
\textbf{basis vectors}}} 
&
\textbf{Indices} 
&
{\splitcellc[]{\textbf{Relationship between} \\ \textbf{covariant bases}}}
& 
{\splitcellc[]{\textbf{Metric} \\ \textbf{coefficients}}}
\\
\hline
\multirow{2.6}{*}{Reference}
&
Cartesian: ${^0}\underline{\boldsymbol{e}}_J, {^0}\underline{\boldsymbol{e}}^J$
&
$J, K, \dots = 1,2,3$
&
${^0}\underline{\boldsymbol{e}}_J 
=
\dfrac{\partial  {^0}x^A}{\partial {^0}z^J}
\ {^0}\underline{\boldsymbol{g}}_A$
&
$\delta_{JK} = {^0}\underline{\boldsymbol{e}}_J \cdot {^0}\underline{\boldsymbol{e}}_K  $
\\[1ex]
\cline{2-5}
&
curvilinear: ${^0}\underline{\boldsymbol{g}}_A, {^0}\underline{\boldsymbol{g}}^A$
&
$A, B, \dots = I, II, III$
&
${^0}\underline{\boldsymbol{g}}_A
=
\dfrac{\partial {^0}z^J}{\partial {^0}x^A} \ {^0}\underline{\boldsymbol{e}}_J$
&
${^0}g_{AB} = {^0}\underline{\boldsymbol{g}}_A \cdot {^0}\underline{\boldsymbol{g}}_B  $
\\ [1ex]
\hline
\multirow{2.6}{*}{Current}
&
Cartesian: ${^t}\underline{\boldsymbol{e}}_j {^t}\underline{\boldsymbol{e}}^j$
&
$j, k, \dots = 1,2,3$
&
${^t}\underline{\boldsymbol{e}}_j
=
\dfrac{\partial  {^t}x^a}{\partial {^t}z^j}
\ {^t}\underline{\boldsymbol{g}}_a$
&
$\delta_{jk} = {^t}\underline{\boldsymbol{e}}_j \cdot {^t}\underline{\boldsymbol{e}}_k  $
\\[1ex]
\cline{2-5}
&
curvilinear: ${^t}\underline{\boldsymbol{g}}_a, {^t}\underline{\boldsymbol{g}}^a$
&
$a, b, \dots = i,ii,iii$
&
${^t}\underline{\boldsymbol{g}}_a
=
\dfrac{\partial {^t}z^j}{\partial {^t}x^a} \ {^t}\underline{\boldsymbol{e}}_j$
&
${^t}g_{ab} = {^t}\underline{\boldsymbol{g}}_a \cdot {^t}\underline{\boldsymbol{g}}_b $
\\[1ex]
\hline
\end{tabular}
\caption{Summary of the basis vectors and their indices in the reference and current configurations. 
}
\label{table:basis in configurations}
\end{threeparttable}
}
\end{table}
Table~\ref{table:basis in configurations} summarises the Cartesian and curvilinear basis vectors in different configurations and provides their indices, as well as the relationship between the Cartesian and curvilinear bases. 
The metric coefficients presented in Table~\ref{table:basis in configurations} can be employed to transform basis vectors between covariant and contravariant, as in ${^0}\underline{\boldsymbol{g}}_A = {^0}g_{AB} {^0}\underline{\boldsymbol{g}}^B$, as well as to change the components between vectors and covectors, e.g., ${^0}x_A = {^0}g_{AB} \ {^0}x^B$. 
Owing to their definitions, the covariant and contravariant metric coefficients are the inverse of each other, meaning that, for instance, ${^0}g_{AB} = \left( {^0}g^{AB} \right)^{-1}.$

Returning to the cavity expansion problem, let the components of the displacements on the inner surface be prescribed as follows
\begin{equation}
{^t_0\underline{\boldsymbol{U}}} 
=
\begin{cases}
{^t_0U^{I}} =  \left( \alpha - 1 \right) {^0x^I}; 
\\
{^t_0U^{II}} = 0;
\\
{^t_0U^{III}} = 0,
\end{cases}
\end{equation}
with $\alpha > 0 $ being the parameter controlling the magnitude of the radial displacements.
Similarly, the displacements in  Cartesian components are given by
\begin{equation}
\label{eq:example_cartesian_displacement}
{^t_0\underline{\boldsymbol{U}}} 
=
\begin{cases}
{^t_0U^{1}} = {^t_0U^{I}} \cos \left( {^0x^{II}}\right) = \left( \alpha - 1  \right) {^0}z^1
\\
{^t_0U^{2}} = {^t_0U^{I}} \sin \left( {^0x^{II}}\right) = \left( \alpha - 1  \right) {^0}z^2
\\
{^t_0U^{3}} = 0.
\end{cases}
\end{equation}
The deformation gradient is the two-point tensor described employing the Cartesian or the curvilinear bases as follows
\begin{equation}
\label{eq:deformation_gradient}
{^t_0\underline{\underline{\boldsymbol{X}}}} 
= 
{^t_0X^{j}_{\ J}} 
\ {^t\underline{\boldsymbol{e}}_j} 
\ {^0\underline{\boldsymbol{e}}^J}
= 
{^t_0X^{a}_{\ A}} 
\ {^t\underline{\boldsymbol{g}}_a} 
\ {^0\underline{\boldsymbol{g}}^A}.
\end{equation}
The components of the deformation gradient at the inner surface of the cylinder can be expressed in the more familiar Cartesian environment, and, using the displacement given in Eq.~\eqref{eq:example_cartesian_displacement}, these are
\begin{equation}
\label{eq:example_cartesian_def_gradient}
\left[ {^t_0X^{j}_{\ J}} \right]
= 
\left[ \delta^{j}_{J} \right]
+ 
\left[ \delta^{j}_{K} \
\frac{\partial {^t_0U^K}}{\partial {^0z^{J}}} \right]
= 
\begin{bmatrix}
1 & 0 & 0 
\\
0 & 1 & 0
\\
0 & 0 & 1
\end{bmatrix}
+
\begin{bmatrix}
\dfrac{\partial {^t_0U^1}}{\partial {^0z^{1}}}
& \dfrac{\partial {^t_0U^1}}{\partial {^0z^{2}}}
& 0
\\[8pt]
\dfrac{\partial {^t_0U^2}}{\partial {^0z^{1}}}
& \dfrac{\partial {^t_0U^2}}{\partial {^0z^{2}}}
& 0
\\[8pt]
0 & 0 & 0
\end{bmatrix}
= 
\begin{bmatrix}
\alpha & 0 & 0
\\
0 & \alpha & 0
\\
0 & 0 & 1
\end{bmatrix},
\end{equation}
where $\delta^{j}_{J}$ is the Kronecker delta symbol and is considered here merely as a mean to fix the indices between the current and original configurations (its role and meaning are more detailed in Section~\ref{subsec:def_gradient_decompositions}). 
Employing the change of basis vectors between Cartesian and curvilinear covariant basis vectors for the current configuration (see Table~\ref{table:basis in configurations}), it follows that
\begin{equation}
{^t}\underline{\boldsymbol{e}}_j
=
\dfrac{\partial  {^t}x^a}{\partial {^t}z^j}
\ {^t}\underline{\boldsymbol{g}}_a
= 
\begin{bmatrix}
\dfrac{{^t}z^1}{{^t}x^i} & \dfrac{{^t}z^2}{{^t}x^i} & 0
\\[8pt]
- \dfrac{{^t}z^2}{\left( {^t}x^i \right)^2} & \dfrac{{^t}z^1}{\left( {^t}x^i \right)^2} & 0
\\[8pt]
0 & 0 & 1
\end{bmatrix}
{^t}\underline{\boldsymbol{g}}_a
=
\begin{bmatrix}
\cos \left( {^t}x^{ii} \right) & \sin \left( {^t}x^{ii} \right) & 0
\\[8pt]
- \dfrac{\sin \left( {^t}x^{ii} \right)}{{^t}x^i } & \dfrac{\cos \left( {^t}x^{ii} \right)}{{^t}x^i } & 0
\\[8pt]
0 & 0 & 1
\end{bmatrix}
 {^t}\underline{\boldsymbol{g}}_a.
\end{equation}
Similarly, the Cartesian contravariant basis vector in the original configuration can be expressed as
\begin{equation}
{^0}\underline{\boldsymbol{e}}^J 
=
\delta^{JK} \ {^0}\underline{\boldsymbol{e}}_J 
=
\delta^{JK}
\dfrac{\partial  {^0}x^B}{\partial {^0}z^K}
\ {^0}g_{BA} {^0}\underline{\boldsymbol{g}}^A
=
\begin{bmatrix}
\cos \left( {^0}x^{II} \right) & -{^0}x^I \sin \left( {^0}x^{II} \right) & 0
\\
\sin \left( {^0}x^{II} \right) & {^0}x^I \cos \left( {^0}x^{II} \right) & 0
\\
0 & 0 & 1
\end{bmatrix}
{^0}\underline{\boldsymbol{g}}^A
\end{equation}
Using the above change of bases and Eq.~\eqref{eq:deformation_gradient}, the components of the deformation gradient in curvilinear coordinates can be expressed as
\begin{equation}
\label{eq:deformation_gradient_component_transformation}
{^t_0X^{a}_{\ A}} 
 = \dfrac{\partial  {^t}x^a}{\partial {^t}z^j} \ {^t_0X^{j}_{\ J}} \ \delta^{JK}
\dfrac{\partial  {^0}x^B}{\partial {^0}z^K}
\ {^0}g_{BA} 
= 
\dfrac{\partial  {^t}x^a}{\partial {^t}z^j} \ {^t_0X^{j}_{\ J}} \ 
\dfrac{\partial {^0}z^J}{\partial  {^0}x^A},
\end{equation}
which yields
\begin{equation}
\left[ {^t_0X^{a}_{\ A}} \right]
 = 
{\scriptsize
\begin{bmatrix}
\alpha \left( \cos \left( {^t}x^{ii} \right) \cos \left( {^0}x^{II} \right)  + \sin \left( {^t}x^{ii} \right) \sin \left( {^0}x^{II} \right) \right)
&
\alpha {^0}x^I \left( \sin \left( {^t}x^{ii} \right) \cos \left( {^0}x^{II} \right) - \cos \left( {^t}x^{ii} \right) \sin \left( {^0}x^{II} \right) \right)
& 
0
\\[8pt]
\dfrac{\alpha}{{^t}x^i} \left( \cos \left( {^t}x^{ii} \right) \sin \left( {^0}x^{II} \right) - \sin \left( {^t}x^{ii} \right) \cos \left( {^0}x^{II} \right) \right)
&
\dfrac{\alpha {^0x^I}}{{^tx^i}} \left( \cos \left( {^t}x^{ii} \right) \cos \left( {^0}x^{II} \right)  + \sin \left( {^t}x^{ii} \right) \sin \left( {^0}x^{II} \right) \right)
& 
0
\\
0 & 0 & 1
\end{bmatrix}.
}
\end{equation}
Considering that ${^0}x^{II} = {^t}x^{ii}$ and $\dfrac{\alpha {{^0}x^I}}{{{^t}x^i}} = 1$, it follows that
\begin{equation}
\label{eq:example_curvilinear_def_gradient}
\left[ {^t_0X^{a}_{\ A}} \right] =
\begin{bmatrix}
\alpha & 0 & 0
\\
0 & 1 & 0 
\\
0 & 0 & 1
\end{bmatrix}.
\end{equation}
From the results given in Eqns \eqref{eq:example_cartesian_def_gradient} and \eqref{eq:example_curvilinear_def_gradient}, two main points emerge: $(i)$ $\det [{^t_0}X^i_{\ J}] \neq \det [{^t_0}X^ a_{\ A}]$, and $(ii)$ the out-of-plane ${^t_0}X^{ii}_{\ II}$ component for the deformation gradient in cylindrical coordinates expressed in Eq.~\eqref{eq:example_curvilinear_def_gradient} is equal to one. 

Concerning the first point, it is well-known that the determinant of the deformation gradient governs the transformation between infinitesimal volumes in the Cartesian case (see, for instance,~\cite{ogden1997non}). 
This raises the question about what corrective factors must be applied when only the curvilinear components ${^t_0}X^ a_{\ A}$ are available, with the goal of ensuring that the transformation between volumes remains invariant regardless of the selected coordinate system~\cite{einstein2003meaning}.
As for the second point, this unitary value of the out-of-plane components indicates what the deformation gradient truly represents.
Despite its name, the deformation gradient is \emph{not} a gradient (see, for instance,~\cite{hughes1983mathematical,kanso2007geometric}), which might justify why some of the literature (see, e.g.,~\cite{danielson1997three,celigoj1998assumed,rauchs2016direct,steigmann2017finite}) includes non-unitary out-of-plane components for similar setups\footnote{
Non-unitary components in the out-of-plane components of the deformation gradients appear when the \emph{physical components} of the deformation gradient, denoted by $^{t}_0X^{<aA>}$, are employed. These components are available when the basis vectors are orthogonal but not normal, as is the case for cylindrical basis vectors. 
In this case, the following transformation formula between components holds $^{t}_0X^{<aA>} = ^{t}_0X^{a}_{\ A} \frac{\sqrt{^{t}g_{aa}}}{\sqrt{^{0}g_{AA}}}$, with no sum implied (compare with~\cite{kohler2001domain}). 
However, if the use of these physical components is not specified, it is assumed that the most natural form of the deformation gradient in curvilinear components is $^{t}_0X^{a}_{\ A}$, where the current covariant basis and original contravariant basis are employed.
}. From the two above observations, it becomes clear that using naive definitions of the deformation gradient and volume transformation in curvilinear coordinates can easily result in inaccuracies. This simple result also leads to a number of other consequences that are examined below.

\section{Kinematics and strain measures}
\label{sec:kinematic}
\noindent This section brings together various contributions scattered in~\cite{truesdell1960classical,hughes1983mathematical,dvorkin2006nonlinear,steigmann2017finite,leborgne2023objectivity,steigmann2025principles}, with the scope of clarifying when kinematics in curvilinear components becomes more intricate than in the Cartesian counterpart, and extends, with the same gist, to elasto-plastic deformations.


\subsection{Configurations, coordinate systems and point mapping}
\begin{figure}
\centering
\includegraphics[width=0.9\textwidth, trim=1cm 0cm 0cm 0cm]{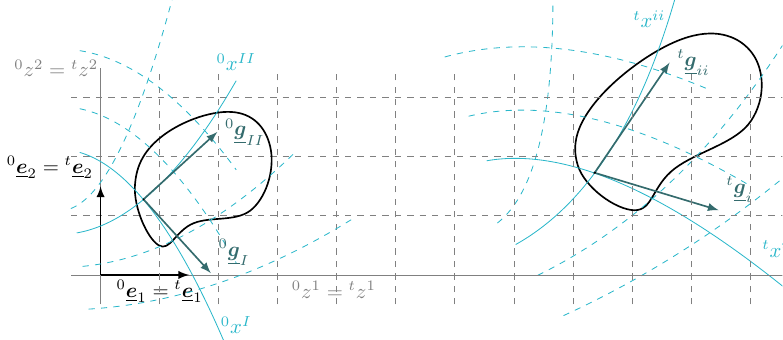}
\caption{Illustration of body $\mathcal{B}$ in the reference and the current configuration. 
Euclidean space is described via Cartesian common basis vectors ${^0\underline{\boldsymbol{e}}_J} ={^t\underline{\boldsymbol{e}}_j}$ tangent to coordinate lines (grey solid lines). 
Curvilinear basis vectors ${^0\underline{\boldsymbol{g}}_A}$, ${^t\underline{\boldsymbol{g}}_a}$ are defined as tangent to the coordinate curves at two generic points in the original and current configuration (light blue solid lines). }
\label{fig:Spaces_and_bases}
\end{figure}
\noindent Let the body $\mathcal{B}$ comprise a fixed number of material points.
Among all the different sub-regions of the Euclidean three-dimensional space $\mathcal{E} \in \mathds{R}^3$ that the body $\mathcal{B}$ can occupy, the reference $^{0}\mathcal{B}$ (at time $t=0$) and current $^{t}\mathcal{B}$ (at generic time $t$) are those of interest for the following developments. 
Material points occupy a position in the reference configuration (a point of the space), generically indicated by $P$.
The mapping used in Lagrangian mechanics, i.e., the \emph{motion}, tracks a material point of the body $\mathcal{B}$ across the different configurations: this is performed by taking the material point of the body $\mathcal{B}$, which is at the point $P$ in the space of the reference configuration $^0\mathcal{B}$ and assigns a point $p$ in the current configuration $^t\mathcal{B}$, i.e., $p= {\phi} \left( P, t \right)$. 
This mapping is observer independent, as it does not require any spatial reference system (and relative bases).
Since this paper focuses on the kinematics at each moment in time, and not on the body's dynamics over time, the dependency from the time is dropped for the rest of this manuscript, i.e, the shorthand notation $p= {\phi} ( P ) $ in adopted in lieu of the above one.

As in  Section~\ref{sec:cavity_expansion}, let the Euclidean space be described by the three right-handed Cartesian basis vectors. 
To distinguish the notation between the configurations, let ${^0}\underline{\boldsymbol{e}}_J$ be the Cartesian basis vectors employed to describe quantities (i.e., vectors and tensors) in the reference configuration and positioned at the point $O$. 
Similarly, let the Cartesian bases vectors ${^t}\underline{\boldsymbol{e}}_j$ be employed in the current configuration and positioned at $o$. 
For simplicity and without loss of generality, this work assumes ${^0}\underline{\boldsymbol{e}}_J = {^t}\underline{\boldsymbol{e}}_j$ with $J, j = 1 \dots 3$ and $O$ and $o$ coincident. 
To each point $P$ in the reference configuration, a position vector can be assigned and its Cartesian bases can be used to express its contravariant components, i.e., $^0\underline{\boldsymbol{z}} = {^0z^{J}} \left( P, O \right) \ {^0}\underline{\boldsymbol{e}}_J \left( O \right) $. 
It can be noticed that, once the origin $O$ is assigned, the initial position vector is purely a function of the point. 
Similarly for the current configuration, the current position vector is $^t\underline{\boldsymbol{z}} = {^tz^{j}} \left(p, o \right) \ {^t}\underline{\boldsymbol{e}}_j \left( o \right) $. 
If, as assumed, the origins of the reference systems are overlapping $O\equiv o$, the dependency from the origin of these functions becomes redundant and can be removed. 
Assuming also that these functions and the mapping $\phi$ are sufficiently regular, the current position can be expressed directly as 
\begin{multline}
^t\underline{\boldsymbol{z}} 
= 
{^tz^{j}} \left( \phi \left( P \right)  \right) \ {^t}\underline{\boldsymbol{e}}_j
=
{^tz^{j}} \left( \phi \left( \left({^0z^{J}}\right)^{-1}  \left( ^0\underline{\boldsymbol{z}}\right) \right)  \right) \ {^t}\underline{\boldsymbol{e}}_j 
= 
{^t\phi^{j}} \left( \left({^0z^{J}}\right)^{-1}  \left( ^0\underline{\boldsymbol{z}}\right) \right)   \ {^t}\underline{\boldsymbol{e}}_j,  
\end{multline}
where the shorthand $\phi^{j} = {^t}z^{j} \left( \phi \right)$ for the composition of functions has been used in the above formula. 
Abusing the previous notation, it is possible to concisely write the above equation as
\begin{equation}
\label{eq:Cartesian_mapping}
{^t}\underline{\boldsymbol{z}} = {^t}\underline{\tilde{\boldsymbol{z}}} \left( {^0}\underline{\boldsymbol{z}} \right),
\end{equation}
where ${^t}\underline{\boldsymbol{z}}$ can be expressed directly as  a function, denoted by ${^t}\underline{\tilde{\boldsymbol{z}}}$, of the initial position vector $
{^0}\underline{\boldsymbol{z}} $.
In the Euclidean context employing position vectors described with persistent Cartesian bases, the above abuse of notation is not harmful, since position vectors and points have trivial relationships, which is guaranteed by the persistency of the flat Euclidean space and its basis vectors ${^0}\underline{\boldsymbol{e}}_J$ and ${^t}\underline{\boldsymbol{e}}_j$.
However, this can lead to misunderstanding when the Euclidean space is not available  or when non-orthonormal basis vectors are employed. 
The latter is the case considered in this manuscript, in which, employing curvilinear basis vectors in the original $^{0}\underline{\boldsymbol{g}}_A \left( \check{O} \right)$ and current configuration $^{t}\underline{\boldsymbol{g}}_a \left( \check{o} \right)$ positioned at the points $\check{O}$ and $\check{o}$, a relationship similar to Eq.~\eqref{eq:Cartesian_mapping}, i.e., relating position vectors in different configurations (which can belong to different curved spaces), is not available.
To avoid confusion, the mapping ${\phi}$ relating points in different configurations must be always privileged, as it does not rely on a Euclidean structure of the space or the choice of a spatial observer (reference system).

To quickly distinguish the position vector in different coordinate systems, let ${^0\underline{\boldsymbol{x}}}$ be the position vector between the point $P$ and the origin of the curvilinear system in the reference configuration. 
Similarly, let ${^t\underline{\boldsymbol{x}}}$ be the position vector between the point $p$ and the origin of the curvilinear system in the current configuration.

\subsection[Deformation gradient: definition]{Deformation gradient: definition\footnotemark}

\footnotetext{
The definition of the deformation gradient considered below is, by no means, the only one proposed in the literature (reference~\cite{truesdell1960classical}, for instance, suggests three different definitions and relative explanations). However, in the authors' opinion, the definition that offers the most clarity and rigour, particularly in the context of of curvilinear coordinates, has been selected below.
}
\begin{figure}
\centering
\includegraphics[width=0.9\textwidth, trim=1cm 0cm 0cm 0cm]{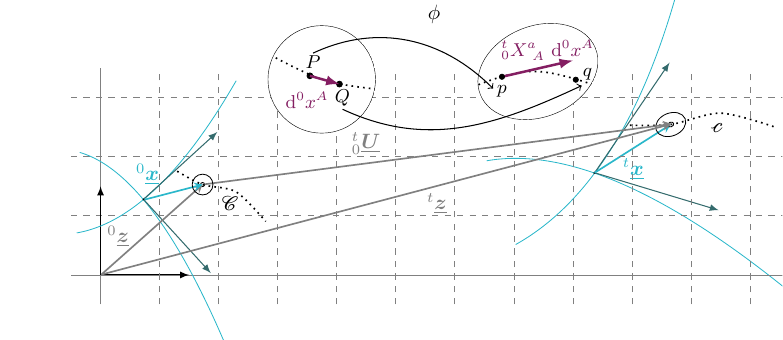}
\caption{Illustration of the initial positions of the point $P$ in the Cartesian ${^0\underline{\boldsymbol{z}}}$ and curvilinear reference ${^0\underline{\boldsymbol{x}}}$ systems. 
The point $p = \phi \left( P \right)$ is defined via the current positions ${^t\underline{\boldsymbol{z}}}$ and ${^t\underline{\boldsymbol{x}}}$ in the Cartesian and coordinate systems, respectively. 
The displacement ${^t_0}\underline{\boldsymbol{U}}$ can  be defined only as the difference between Cartesian position vectors.
Points $P,Q$ belong to the material curve $\mathscr{C}$, which, mapped via $\phi$, result in the points $p,q$ belonging $\mathscr{c}$. 
The magnification shows a tangential portion of material curves defined by $P$ and $Q$ in the original configuration and how the deformation gradient describes its first-order approximation under the mapping $\phi$ in the current configuration.}
\label{fig:deformation_gradient}
\end{figure}
Suppose a curve $\mathscr{C}$ is made by a collection of points in the reference configuration.
If each of these points is mapped via $\phi$ to the current configuration, these points define a curve $\mathscr{c} = \phi \left( \mathscr{C} \right)$. $\mathscr{C}$ and $\mathscr{c}$ are called \emph{material curves} to indicate that the same points constitute them.
If $P, \, Q$ are points belonging to $\mathscr{C}$ and $p = \phi (P)$, $q = \phi (Q)$ their respective mapped points belonging to $\mathscr{c}$, the Taylor expansion up of the first order of the mapping $\phi$ at $Q$ is given by
\begin{equation}
\phi (Q) = \phi (P) + \text{d}\phi(P) \left( Q - P \right) + \mathscr{o} \left( || Q - P || \right),
\end{equation}
where $\mathscr{o} (h)$  identifies the little-o Landau symbol such that for $h \rightarrow0$ then $\mathscr{o} (h)/h \rightarrow 0$. The above equation can be rearranged into 
\begin{equation}
\label{eq:deformatioN-gradient_as_differential}
\dfrac{q - p}{|| Q - P ||} = \text{d}\phi(P) \dfrac{\left( Q - P \right)}{|| Q - P ||} + \dfrac{\mathscr{o} \left( || Q - P || \right)}{|| Q - P ||}.
\end{equation}
As $Q \rightarrow P$, the \emph{differential} $\text{d}\phi(P)$ gives a relationship between tangent vectors to the material curves (see Fig.~\ref{fig:deformation_gradient}). 
Recognising that the basis covectors (and, similarly, the vectors) arise from a differential operator~\cite{hughes1983mathematical,frankel2004geometry}, i.e., $\text{d} {^0}x^{A}  = {^0}\boldsymbol{\underline{g}}^{A}$, it can be observed that the differential of the point mapping $\phi$ yields to a quantity that behaves like a covector
\begin{equation}
\text{d}\phi(P) 
= 
\dfrac{\partial \phi(P)}{\partial {^0}x^{A}} \ \text{d} {^0}x^{A} 
= 
\dfrac{\partial \phi(P)}{\partial {^0}x^{A}} \ {^0}\boldsymbol{\underline{g}}^{A}.
\end{equation}
Based on the above equation, it can be seen how the above differential $\text{d}\phi(P)$ does not depend on the specific choice of reference system~\cite{leborgne2023objectivity}; however, it can be expressed with the assistance of a reference system and its bases, as in the case of ${^0}\boldsymbol{\underline{g}}^{A}$ in the above equation.

Using the composition of functions $\phi^{a} = {^t}x^a \left( \phi \right)$, the \emph{deformation gradient} can be expressed in curvilinear coordinates as follows
\begin{equation}
{^t}x^a \left( \text{d}\phi(P) \right) = \text{d}\phi^{a}(P) = \frac{\partial \phi^{a} (P)}{\partial ^{0}x^{A}} \ \text{d}^{0}x^{A} \coloneqq {^{t}_0X^{a}_{\ A}} \ \text{d}^{0}x^{A},
\end{equation} 
If using Cartesian coordinates and the above-mentioned abuse of notation following Eq.~\eqref{eq:Cartesian_mapping}, it can be seen that
\begin{equation}
\text{d}\phi^{j}(P) 
= 
\frac{\partial \phi^{j} (P)}{\partial ^{0}z^{J}} \ \text{d}^{0}z^{J} 
=  
\frac{\partial \tilde{z}^{j} ({^0}\boldsymbol{\underline{z}})}{\partial ^{0}z^{J}} \ \text{d}^{0}z^{J} \coloneqq 
{^{t}_0X^{j}_{\ J}} \ \text{d}^{0}z^{J}.
\end{equation}

This last abuse of notation, where the vector function $\tilde{z}^{j} ({^0}\boldsymbol{\underline{z}})$ is employed in place of the point mapping $\phi$, might motivate why the deformation gradient is sometimes misunderstood as a derivative of a vector, i.e.,  derivative necessitating a covariant derivative (see Section~\ref{subsec:covariant_derivatives}). 

\subsection{The decompositions of the deformation gradient and the shifter}
\label{subsec:def_gradient_decompositions}
The deformation gradient is a real-valued square matrix. As such, it admits the following polar decompositions
\begin{equation}
\label{eq:RU_decomposition}
{{^t_0}\underline{\underline{\boldsymbol{X}}}} 
= 
{{^t_0}\underline{\underline{\boldsymbol{Q}}}}  \cdot {{^0}\underline{\underline{\boldsymbol{U}}}} 
=  
{{^t}\underline{\underline{\boldsymbol{v}}}} \cdot {{^t_0}\underline{\underline{\boldsymbol{Q}}}},
\end{equation}
with ${{^0}\underline{\underline{\boldsymbol{U}}}}$ and ${{^t}\underline{\underline{\boldsymbol{v}}}}$ being the symmetric right and left stretch tensors, and ${{^t_0}\underline{\underline{\boldsymbol{Q}}}}$ a proper orthogonal matrix indicating a rigid body motion (i.e., rotation and translation). 
At the same time, the following singular value decomposition of the deformation gradient is admitted (see~\cite{brannon2018rotation})
\begin{equation}
\label{eq:SVD_decomposition}
{{^t_0}\underline{\underline{\boldsymbol{X}}}} 
= 
{{^t}\underline{\underline{\boldsymbol{r}}}} \cdot {{^t}\underline{\underline{\boldsymbol{\varsigma}}}} \cdot {{^t_0}\underline{\underline{\boldsymbol{S}}}}.
\end{equation}
Comparing Eqs.~\eqref{eq:RU_decomposition} with~\eqref{eq:SVD_decomposition}, it can be seen that 
\begin{align}
{{^t_0}\underline{\underline{\boldsymbol{Q}}}} 
& = 
{{^t}\underline{\underline{\boldsymbol{r}}}} 
\cdot 
{{^t_0}\underline{\underline{\boldsymbol{S}}}};
\\
{{^0}\underline{\underline{\boldsymbol{U}}}} 
& = 
{{^t_0}\underline{\underline{\boldsymbol{S}}}}^T \cdot  {{^t}\underline{\underline{\boldsymbol{\varsigma}}}} \cdot 
{{^t_0}\underline{\underline{\boldsymbol{S}}}};
\\
{{^t}\underline{\underline{\boldsymbol{v}}}} 
& = 
{{^t}\underline{\underline{\boldsymbol{r}}}} \cdot  {{^t}\underline{\underline{\boldsymbol{\varsigma}}}} \cdot 
{{^t}\underline{\underline{\boldsymbol{r}}}}^T,
\end{align}
where ${{^t}\underline{\underline{\boldsymbol{r}}}}$ is a pure rotation in the current configuration, ${{^t}\underline{\underline{\boldsymbol{\varsigma}}}}$ a diagonal scaling in the same configuration, and 
${{^t_0}\underline{\underline{\boldsymbol{S}}}}$ indicates a translation between the current and the original configuration. 
Owing to its physical meaning, ${{^t_0}\underline{\underline{\boldsymbol{S}}}}$ is referred to as the \emph{shifter} tensor, e.g.,~\cite{truesdell1960classical}. 
As a matter of fact, the components of the shifter in Cartesian coordinates have been already introduced in Eq.~\eqref{eq:example_cartesian_def_gradient} and indicated by $\delta^{j}_{J}$. 
Changing these to curvilinear coordinates (see, again, Table~\ref{table:basis in configurations}) gives the components of the shifter in curvilinear coordinates
\begin{equation}
\label{eq:shifter_components}
{^t_0{S}}^a_{\ A}
= \dfrac{\partial {^t}x^a}{\partial {^t}z^j} \left( \phi \left( P \right) \right) \ \delta^j_{ \ J} \ \dfrac{\partial {^0}z^J}{\partial {^0}x^A} \coloneqq {^t}\underline{\boldsymbol{g}}^a \left( \phi \left( P \right) \right) \cdot {^0}\underline{\boldsymbol{g}}_A.
\end{equation}
Hence, while the shifter components play a silent role in the Cartesian setting (they are employed to make the indices work if their notations differ between the current and original configurations), this is not the case in the most general curvilinear case. 
Following its definition Eq.~\eqref{eq:shifter_components}, it can be seen that the shifter tensor being orthogonal satisfies the following properties ${^t_0\underline{\underline{\boldsymbol{S}}}}^T = {^t_0\underline{\underline{\boldsymbol{S}}}}^{-1} = {^0_t\underline{\underline{\boldsymbol{S}}}} $.

\subsection{The transpose of the deformation gradient and strain measures}
To compute strain and deformation quantities, the transpose of the deformation gradient is also necessary.
The transpose of the deformation gradient relies on the inner product between vectors, which, in turn, relies on the metric. 
If ${{^0}\underline{\boldsymbol{W}}}$ and ${{^t}\underline{\boldsymbol{w}}}$ are two vectors tangent to the material curves in the original configuration and current configurations, the \emph{transpose} or \emph{adjoint} of the deformation gradient is defined by the equation based on the conservation of the following dot products
\begin{equation}
\left( {{^t_0}\underline{\underline{\boldsymbol{X}}}}^T  \cdot {{^t}\underline{\boldsymbol{w}}} \right) \cdot {{^0}\underline{\boldsymbol{W}}} 
\coloneqq 
\left( {{^t_0}\underline{\underline{\boldsymbol{X}}}} \left( \phi^{-1} \left( p \right) \right) \cdot {{^0}\underline{\boldsymbol{W}}} \right) \cdot {{^t}\underline{\boldsymbol{w}}}.
\end{equation}
If employing the curvilinear coordinates, the above definition provides the following components of the transpose of the deformation gradient
\begin{equation}
\nonumber
\left( {{^t_0}X^T}  \right)^A_{\ a} \ {{^t}w^a} \ {{^0}W_A} 
= 
{{^t_0}X}^b_{\ B} \left( \phi^{-1} \left( p \right) \right) \ {{^0}W^B} \ {{^t}w_b},
\end{equation}
which gives
\begin{equation}
\label{eq:def_gradient_transpose_components}
\left( {{^t_0}X^T} \right)^A_{\ a} 
= 
{{^t_0}X}^b_{\ B} \left( \phi^{-1} \left( p \right) \right) \ {{^0}g^{BA}} \left( \phi^{-1} \left( p \right) \right) \ {{^t}g_{ab}}.
\end{equation}
In the case of Cartesian coordinate systems, the transpose of the deformation gradient becomes the well-known matrix where rows and columns are inverted, given that both metric coefficients corresponds to the Kronecker delta (see Table~\ref{table:basis in configurations}).

With a definition of the transpose of the deformation gradient now established, the right and left Cauchy-Green deformation tensors can be defined as follows
\begin{equation}
\label{eq:C_def}
{^0} \underline{\underline{\boldsymbol{C}}} 
\coloneqq
{^t_0}\underline{\underline{\boldsymbol{X}}}^T \cdot {^t_0}\underline{\underline{\boldsymbol{X}}}  
=
\underbrace{ {^t_0}X^b_{\ C} {^0g^{CA}} \ \left( {^tg_{ab}} \left( \phi \left( P \right) \right) \right) \   {^t_0}X^a_{\ B} }_{ = {^0{C}}^A_{\ B}} {^0{\underline{\boldsymbol{g}}}_A} \ {^0{\underline{\boldsymbol{g}}}^B},
\end{equation}
and
\begin{equation}
\label{eq:b_def}
{^t} \underline{\underline{\boldsymbol{b}}} 
\coloneqq
{^t_0}\underline{\underline{\boldsymbol{X}}} \cdot {^t_0}\underline{\underline{\boldsymbol{X}}}^T
=
\underbrace{ {^t_0}X^a_{\ A} \left( \phi^{-1} \left( p \right) \right) \ {^t_0}X ^c_{\ B} \left( \phi^{-1} \left( p \right) \right) \ {^0g^{BA}} \left( \phi^{-1} \left( p \right) \right) \ {^tg_{cb}} }_{= {^t}{b}^a_{\ b}}
\ {^t{\underline{\boldsymbol{g}}}_a} \ {^t{\underline{\boldsymbol{g}}}^b}
.
\end{equation}
The components of the left and right Cauchy-Green deformation measures expressed in the above equations are often termed \emph{mixed} due to the use of covariant and contravariant bases. 
This mixed notation for symmetric measures of strain (and stress) is particularly advantageous, as it eliminates the need to change their components when calculating invariants. 
As an example, let us consider the trace of ${^t} \underline{\underline{\boldsymbol{b}}} $, i.e.,
\begin{equation}
\text{tr} \left( {^t} \underline{\underline{\boldsymbol{b}}} \right) \coloneqq {^t} \underline{\underline{\boldsymbol{b}}}  \boldsymbol{:} {^t}\boldsymbol{\underline{\underline{1}}} 
=
{^t}{b}^a_{\ b} \ \delta^b_{\ a} = {^t}{b}^{ab} \ ^t{g}_{ab},
\end{equation}
where $\boldsymbol{:}$ is the double contraction symbol and the current second order identity matrix ${^t}\boldsymbol{\underline{\underline{1}}}$ becomes the Kronecker delta or the metric coefficients depending on the components of ${^t} \underline{\underline{\boldsymbol{b}}} $.

In addition to the definitions of deformations, this work also adopts the logarithmic strain based on ${^t}\underline{\underline{\boldsymbol{b}}} $ and defined by
\begin{equation}
{^t} \underline{\underline{\boldsymbol{\epsilon}}} 
\coloneqq
\frac{1}{2} \log {^t} \underline{\underline{\boldsymbol{b}}}.
\end{equation}

\subsection{The volume transformation or Jacobian}
\noindent The relationship between the Cartesian and curvilinear components of the deformation gradient Eq.~\eqref{eq:deformation_gradient_component_transformation} can also be employed to establish a relationship between infinitesimal volumes in the original $\text{d}{^0V}$ and current configuration $\text{d}{^tV}$. 
Let the \emph{Jacobian} ${{^t_0}J}$ be defined as the ratio between these two infinitesimal volumes, i.e.,
\begin{equation}
\label{eq:Jacobian_definition}
{{^t_0}J} \ \text{d}{{^0}V} \left( \phi^{-1} \left( p \right) \right) \coloneqq \text{d}{{^t}V}.
\end{equation}
As pointed out in Section~\eqref{sec:cavity_expansion}, it is well-known that ${{^t_0}J} = \det \left[ {{^t_0}X}^j_{\ J} \right]$ gives the ratio of the volumes of the parallelepipeds centred at $P$ and $p = \phi \left( P \right)$. Employing the inverse of Eq.~\eqref{eq:deformation_gradient_component_transformation}, it can be seen that
\begin{equation}
\label{eq:Jacobian_components}
{{^t_0}J} = \det \left[ {{^t_0}X}^j_{\ J} \right] 
=
\det \left[ \dfrac{\partial {^t}z^j}{\partial  {^t}x^a} \ {^t_0X^{a}_{\ A}} \ 
\dfrac{\partial  {^0}x^A}{\partial {^0}z^J}
\right]
= 
\det \left[ {{^t_0}X}^a_{\ A} \right] \dfrac{\sqrt{ \det \left[ {{^t}g_{ab}} \right]}}{\sqrt{ \det \left[ {{^0}g_{AB}} \right]}},
\end{equation}
since $\det \left[ \frac{\partial {^t}z^j}{\partial {^t}x^a}
\right] = \sqrt{\det \left[ {{^t}g_{ab}} \right]}$ and $\det \left[ \dfrac{\partial  {^0}x^A}{\partial {^0}z^J} \right] = \sqrt{\det \left[ {{^0}g^{AB}} \right]} = \sqrt{\det \left[ {{^0}g_{AB}}^{-1} \right]}$.
The presence of the determinant of the metric coefficients is necessary to consider that the transformation of volumes employing curvilinear coordinates does not track volumes of parallelepipeds. 
In the case of cylindrical coordinates, for instance, each small volume is a small portion obtained by the revolution of $2\pi$ of a rectangle. 

Based also on the definition of the strain measures Eq.~\eqref{eq:C_def} and~\eqref{eq:b_def}, it can also be appreciated how 
\begin{equation}
\label{eq:Jacobian_Cauchy-Green}
{{^t_0}J} = \sqrt{\det \left[ {^0{C}}^A_{\ B} \right]}= \sqrt{\det \left[{^t}b^{a}_{\ b} \right] }.
\end{equation}
Either computed via Eq.~\eqref{eq:Jacobian_components} or via Eq.~\eqref{eq:Jacobian_Cauchy-Green}, the metric coefficients (current and original) appear in the calculation of the Jacobian. 
This point is of particular interest for the elasto-plastic case, where it is not trivial to determine the metric coefficient of the intermediate, stress-free configuration. 
This leads to the challenge of computing separately the elastic and plastic parts of the Jacobian (see Section~\ref{subsec:elasto-plasticity}).

\subsection{Covariant derivatives and gradients}
\label{subsec:covariant_derivatives}
\noindent When curvilinear coordinates depending on the local position are employed, the differentiation rule of vectors and tensors must take the dependency of the bases also into account. 
This concept leads to the definition of the Christoffel symbol of the second kind, denoted by ${^0{\Gamma}}^C_{AB}$, as follows
\begin{equation}
\label{eq:Christoffel_symbol_original}
\frac{\partial {^0}\boldsymbol{\underline{g}}_A}{\partial ^{0}x^B}
\coloneqq
{^0{\Gamma}}^C_{AB} \  {^0}\boldsymbol{\underline{g}}_C.
\end{equation}
A similar definition of this quantity can be given for the derivative with respect to the current curvilinear coordinate and is denoted by ${^t{\gamma}^{c}_{ab}}$.
Following its definition, it is also worth noticing that, if the position vector is twice differentiable, these properties for the Christoffel symbol of the second kind holds 
\begin{align}
{^0{\Gamma}}^C_{AB} & = {^0{\Gamma}}^C_{BA};
\\
\label{eq:der_covector_basis}
- {^0{\Gamma}}^C_{AB} \ {^0}\boldsymbol{\underline{g}}^A & = \frac{\partial {^0}\boldsymbol{\underline{g}}^C}{\partial ^{0}x^B}.
\end{align}
If then $^0\boldsymbol{\underline{W}}$ is a vector field, if follows that its derivative is given by
\begin{equation}
\frac{\partial {^0}\boldsymbol{\underline{W}} }{\partial {^0}x^B} 
= 
\frac{\partial }{\partial ^{0}x^B} \left( {^0}W^A \ {^0}\boldsymbol{\underline{g}}_A \right)
= 
\left( \frac{\partial  {^0}W^A}{\partial ^{0}x^B} + {^0{\Gamma}}^A_{BC} \ {^0}W^C  \right) {^0}\boldsymbol{\underline{g}}_A
=
{^0}W^A \bigl|_{B} \ {^0}\boldsymbol{\underline{g}}_A,
\end{equation}
where ${^0}W^A |_{B}$ is also called the \emph{covariant derivatives} of the vector components. 
Following the same rationale, in the case where a second-order tensor is derived, the Christoffel symbol of the second kind appears twice, with signs depending on the use of contravariant or covariant components, according to Eqs.~\eqref{eq:Christoffel_symbol_original} and~\eqref{eq:der_covector_basis}.
In a Cartesian setting, where the dependency of the basis for the local position is not taken into account, the differentiation rule simply leads to a partial derivative, i.e.,
\begin{equation}
{^0{\Gamma}}^J_{KL} = {^t{\gamma}}^j_{kl} = 0.
\end{equation}

Employing the covariant derivative, it is possible to define the \emph{divergence} of a vector field $^0\boldsymbol{\underline{W}}$ (or, more generally, of a contravariant tensor) as the covariant derivative where the last contravariant index is contracted, i.e.,
\begin{equation}
\underline{\textbf{Div}} \left( {^0}\boldsymbol{\underline{W}} \right) = {^0} W^A|_{A}.
\end{equation}
If the vector under consideration can be seen as contraction of a vector with a mixed tensor, i.e., ${^0}W^A = {^0}V^B \ {^0}H_{B}^{\ A}$, the \emph{gradient} of a vector field is naturally defined as follows
\begin{equation}
\underline{\textbf{Div}} \left( {^0}\boldsymbol{\underline{W}} \right) 
= 
\underline{\textbf{Div}} \left( {^0}\boldsymbol{\underline{V}} \cdot {^0}\boldsymbol{\underline{\underline{H}}} \right) 
=
\underline{\textbf{Grad}} \left( {^0}\boldsymbol{\underline{V}} \right) \boldsymbol{:} {^0}\boldsymbol{\underline{\underline{H}}} + 
 {^0}\boldsymbol{\underline{V}} \cdot \underline{\textbf{Div}} \left( {^0}\boldsymbol{\underline{\underline{H}}} \right),
\end{equation}
where, in components, $\underline{\textbf{Grad}} \left( {^0}\boldsymbol{\underline{V}} \right) = \left( {^0}V^B\right)|_{A}  \ {^0}\boldsymbol{\underline{g}}_B \ {^0}\boldsymbol{\underline{g}}^A$. 
The notations $\underline{\textbf{div}} \left( \bullet \right)$ and $\underline{\textbf{grad}} \left( \bullet \right)$ denote similar operators, with the difference that the differentiation is made with respect to the current coordinates. 
Additionally, if the vector ${^t}\boldsymbol{\underline{W}}$ and ${^t}\boldsymbol{\underline{w}}$ both live in the current configuration but satisfy different dependencies, i.e., ${^t}W^a \left( \phi^{-1} \left( p \right) \right) = {^t}w^a \left( p\right)$, then the following relationship hold
\begin{equation}
\label{eq:gradient_transformation}
\underline{\textbf{Grad}} \left( {^t}\boldsymbol{\underline{W}} \right) 
=
{^t}W^a|_B \ {^t}\boldsymbol{\underline{g}}_a \  {^0}\boldsymbol{\underline{g}}^B
=
{^t}w^a|_c \ {^t_0}X^c_{\ B}  \ {^t}\boldsymbol{\underline{g}}_a \  {^0}\boldsymbol{\underline{g}}^B
= \underline{\textbf{grad}} \left( {^t}\boldsymbol{\underline{w}} \right) \cdot {^t_0\boldsymbol{\underline{\underline{X}}}}.
\end{equation}

\subsection{Deformation gradient based on displacements}
\label{subsec:deformation_gradient_on_displacements}
\noindent In the Euclidean persistent environment, it is also possible to define the total displacement as the difference between the current and the original position\footnote{
In a non-Euclidean space, the displacement is not defined~\cite{hughes1983mathematical} because of the lack of a persistent environment between two configurations.
In a geometric view of mechanics, this motivates why the displacements are vectors but not vector fields~\cite{leborgne2023objectivity}, i.e., they do not belong to the tangent bundle of any manifold.
}.
However, even in the Cartesian setting, when it is desired to express the components of the displacement vector, the question of whether to use the basis vectors in the original configuration rather than the current one arises. 
To remove this ambiguity, the components of the shifter ${^t_0\boldsymbol{\underline{\underline{S}}}}$ (see Eq.~\eqref{eq:shifter_components}) can be employed, leading to
\begin{equation}
\label{eq:displacements_Cartesian_components}
^t_0{U}^J = {\delta^J_{\ j}} \left( \phi \left( P \right) \right) \ ^t{z}^j \left( \phi \left( P \right) \right) - {^0{z}}^J;
\qquad
^t_0{u}^j = {^t{z}}^j - {\delta^j_{\ J}} \left( \phi^{-1} \left( p \right) \right) \ {^0{z}}^J \left( \phi^{-1} \left( p \right) \right),
\end{equation}
where the displacement vector ${^t_0\underline{\boldsymbol{U}}}$ is applied in $P$ and ${^t_0\underline{\boldsymbol{u}}}$ in $p$, as clarified by the following equation
\begin{equation}
\label{eq:displacements_and_shifter}
{^t_0\underline{\boldsymbol{u}}} = {^t_0} \boldsymbol{\underline{\underline{S}}} \left( \phi^{-1} \left( p \right) \right) \cdot {^t_0\underline{\boldsymbol{U}}} \left( \phi^{-1} \left( p \right) \right).
\end{equation}
It can be noted that the components of the shifter tensor using the Cartesian basis vectors do not offer any further insight, if not removing the ambiguity regarding the indices. 
This is expected since, in the Cartesian environment, the shifter tensor offers a parallel translation in a flat space, which does not affect the displacement components. 
However, when curvilinear coordinate systems or non-Euclidean spaces are considered, it can be seen that the components of the displacement change since
\begin{equation}
^t_0{U}^A = {^0_t{S}^A_{\ a}} \left( \phi \left( P \right) \right) \ ^t{x}^a \left( \phi \left( P \right) \right) - {^0{x}}^A;
\qquad
^t_0{u}^a = {^t{x}}^a - {^t_0{S}^a_{\ A}} \left( \phi^{-1} \left( p \right) \right) \ {^0{x}}^A \left( \phi^{-1} \left( p \right) \right).
\end{equation}

\subsubsection{Total deformation gradient}
\noindent While the definition of the deformation gradient Eq.~\eqref{eq:deformatioN-gradient_as_differential} is crucial to understand its physical meaning, the vast majority of the FEM programs compute the deformation gradient based on displacements. 
If the total (i.e., between the original and the current configuration) displacements $^t_0\boldsymbol{\underline{U}}$ are available, then the deformation gradient can be computed via the following formula
\begin{equation}
\label{eq:deformation_gradient_components_material_displ}
{^t_0{X}}^{a}_{\ A}  
=
{^t_0{S}}^a_{\ B} \left( \delta^{B}_{\ A} + {^t_0{U}^B}|_{A} \right).
\end{equation}
This calculation can also be derived from the Cartesian case (compare with Eq.~\eqref{eq:example_cartesian_def_gradient}) by extending the differentiation to include curvilinear coordinates (see covariant derivatives in Section~\ref{subsec:covariant_derivatives}) and using the components of the shifter in that context, Eq.~\eqref{eq:shifter_components}. 

The inverse of the deformation gradient is also a quantity of interest. 
Its expression in curvilinear components can be computed by applying the chain rule to the covariant derivative appearing in  Eq.~\eqref{eq:deformation_gradient_components_material_displ}, 
resulting in
\begin{equation}
{^t_0{X}}^{a}_{\ A}  
=
{^t_0{S}}^a_{\ B} \left( \delta^{B}_{\ A} + {^t_0{U}^B}|_{c} \ {^t_0}X^{c}_{\ A} \right),
\end{equation}
which, postmultiplied by $\left( {^t_0}X^{-1} \right)^{A}_{\ \ d}$, gives
\begin{equation}
\delta^a_{ \ d} = {^t_0{S}}^a_{\ A} \ \left( {^t_0}X^{-1} \right)^{A}_{\ \ d} 
+ {^t_0{S}}^a_{\ B} \ {^t_0{U}^B}|_{d}.
\end{equation}
The above equation can be rearranged into
\begin{equation}
\label{eq:deformation_gradient_components_spatial_displ} 
\left( {^t_0{X}^{-1}} \right)^{A}_{\ \ d} 
= {^0_t{S}^{A}_{\ \ a}} \ \left(
\delta^{a}_{\ d} - {^t}{u}^{a}|_d \right),
\end{equation}
by noticing that ${^t_0{S}}^a_{\ B} \left( \phi^{-1} \left( p \right) \right) \ {^t_0{U}^B}|_{d} \left( \phi^{-1} \left( p \right) \right) = {^t_0{u}^a}|_{d}$.

\subsubsection{Incremental deformation gradient}
\noindent Following the incremental procedure via pseudotime discretisation necessary to run non-linear numerical analyses, multiple configurations, in addition to the initial at $t= 0$ and the current one at time $t$, become available. 
Depending on the structure of the FE code or the numerical method adopted\footnote{
This is, for instance, the case of the Material Point Method~\cite{sulsky1994particle}, where the computational Finite Element grid is introduced anew at the beginning of each time step and disposed at its end. Hence, the solution is computed only incrementally, see, e.g.,~\cite{charlton2017igimp}
}, some programs do not have the total displacement ${^t_0\underline{\boldsymbol{U}}}$ available, but only the incremental displacement ${^t_{t_n}}\underline{\boldsymbol{U}}$ between the previously converged configuration at time $t_n$ and the current one at $t$, i.e., ${^t_0\underline{\boldsymbol{U}}} = {^{t_n}_0}\underline{\boldsymbol{U}} + {^t_{t_n}\underline{\boldsymbol{U}}}$. 
For this newly-introduced configuration, let the point ${^{t_n}}P = {^{t_n} \phi} \left( P \right)$, and another set of curvilinear basis vectors be denoted by ${^{t_n}}\boldsymbol{\underline{g}}_{\mathcal{A}}$ with $\mathcal{A} = \widetilde{I},\widetilde{II},\widetilde{III}$. 

In this context, the incremental total deformation gradient cannot be directly computed via Eq.~\ref{eq:deformation_gradient_components_material_displ}, but the chain rule can be exploited as follows
\begin{equation}
\label{eq:deformation_gradient_incrementally_computed}
{^t_0\underline{\underline{\boldsymbol{X}}}} 
= 
{^t_{t_n}\underline{\underline{\boldsymbol{X}}}} \cdot {^{t_n}_0\underline{\underline{\boldsymbol{X}}}}
=
{^t_{t_n}{X}}^{a}_{\ \mathcal{A}} \left( {^{t_n} \phi}^{-1} ( {^{t_n}}P ) \right)  \ {^{t_n}_0{X}}^{\mathcal{A}}_{\ A} \  {^t}\underline{\boldsymbol{g}}_a \left( \phi \left( P \right) \right) 
\ {^0}\underline{\boldsymbol{g}}^A. 
\end{equation}
The incremental deformation gradient ${^t_{t_n}\underline{\underline{\boldsymbol{X}}}}$ follows a similar definition to that of Eq.~\ref{eq:deformation_gradient_components_material_displ}, i.e.,
\begin{equation}
\label{eq:increment_deformation_gradient_components_explicit}
{^t_{t_n}{X}}^{a}_{\ \mathcal{A}}  
 =
{^t_{t_n}{S}}^a_{\ \mathcal{B}} 
\Bigl( 
\delta^{\mathcal{B}}_{\ \mathcal{A}} + \underbrace{
\frac{\partial {^{t}_{t_n}{U}^\mathcal{B}}}{\partial {^{t_n}{x}^{\mathcal{A}}}} + {^{t_n}{\Gamma}}^{\mathcal{B}}_{\mathcal{AC}} \ {^{t}_{t_n}{U}^{\mathcal{C}}}
}_{\coloneqq  ^{t}_{t_n}{U}^{\mathcal{B}} |_{\mathcal{A}} } 
\Bigr),
\end{equation}
with the differences being the incremental shifter tensor of components ${^t_{t_n}{S}}^a_{\ \mathcal{A}} \coloneqq {^t}\underline{\boldsymbol{g}}^a \cdot {^{t_n}}\underline{\boldsymbol{g}}_{\mathcal{A}}$ and the Christoffel symbol of the second kind at the previously-converged configuration, i.e., $\dfrac{\partial {^{t_n}}\boldsymbol{\underline{g}}_{\mathcal{A}}}{\partial ^{t_n}x^{\mathcal{B}}}
\coloneqq
{^{t_n}{\Gamma}}^{\mathcal{C}}_{\mathcal{AB}} \  {^{t_n}}\boldsymbol{\underline{g}}_{\mathcal{C}}$.

Similarly to Eq.~\eqref{eq:deformation_gradient_components_spatial_displ}, the inverse of the incremental deformation gradient is given by
\begin{equation}
\label{eq:inv_incr_deformation_gradient} 
\left( {^t_{t_n}}{X}^{-1} \right)^{\mathcal{A}}_{\ \ b} 
= {{^0_{t_n}}{S}^{\mathcal{A}}_{\ \ a}} \ \left(
\delta^{a}_{\ b} - {^t_{t_n}}{u}^{a}|_b \right).
\end{equation}

\subsection{Elasto-plasticity}
\label{subsec:elasto-plasticity}
\noindent Based on the work of Bilby~\cite{bilby1957continuous}, Kr{\"o}ner~\cite{kroner1959allgemeine} and Lee~\cite{lee1969elastic}, the multiplicative decomposition of the deformation gradient postulates the existence of an intermediate, stress-free configuration, defined by the following equation
\begin{equation}
\label{eq:elasto-plastic_def_gradient}
{^t_0}\boldsymbol{\underline{\underline{X}}} = {^t_0}\boldsymbol{\underline{\underline{X}}}^e \cdot 
{^0_0}\boldsymbol{\underline{\underline{X}}}^p
\end{equation}
where ${^t_0}\boldsymbol{\underline{\underline{X}}}^e$ is the elastic part of the deformation gradient and ${^0_0}\boldsymbol{\underline{\underline{X}}}^p$ its plastic counterpart.
However, generally speaking, this intermediate configuration presents two main issues: $(i)$ it might not be explicitly known~\cite{goodbrake2021mathematical}; and $(ii)$ it is not uniquely defined~\cite{casey1980remark}. While the second point is less critical because the stress remains independent of the specific decomposition, the first point poses a significant challenge. 
If the intermediate configuration is unknown, it becomes impossible to determine its relative metric, which, in turn, makes it difficult to define the elastic and plastic components of the Jacobian.
Nonetheless, following the approach proposed in~\cite{goodbrake2021mathematical}, the existence of the metric $\mathring{g}_{AB}$, defined by 
\begin{equation}
\label{eq:goodbrake_metric}
\mathring{g}_{AB} 
\coloneqq
{^0}g_{CD}  \ \left( {^0_0}X^p \right)^{C}_{\ A} \ \left( {^0_0}X^p \right)^{D}_{\ B},
\end{equation}
is equivalent to the composition of the mapping $\phi =  \tilde{\phi} \ ( \text{Id}_{^0\mathcal{B}} )$. 
The differential of the mappings $\text{Id}_{^0\mathcal{B}}$ and $ ( \tilde{\phi} )$ defines the elastic and plastic deformation gradients, i.e.,
\begin{equation}
\left( {^t_0}X^{e} \right)^{a}_{\ B} \coloneqq \dfrac{\partial \tilde{\phi}^{a}}{\partial {^0}x^{B}},
\qquad 
\left( {^0_0}X^p \right)^{A}_{\ B} \coloneqq \dfrac{\partial ( \text{Id}_{^0\mathcal{B}} )^{A}}{ \partial {^0}x^{B}}.
\end{equation}
The metric coefficients $\mathring{g}_{AB}$ are crucial for defining quantities that require a transpose operator of the deformation gradient ${^t_0}\boldsymbol{\underline{\underline{X}}}$, as it can be seen in the case of the elastic left Cauchy-Green deformation tensor, i.e., 
\begin{align}
\nonumber
( {^t}b^e )^{a}_{ \ b} 
& \coloneqq \left( {^t_0}X^e \right)^{a}_{\ A} \ \left( ( {^t_0}X^e )^T \right)^{A}_{\ b} 
=
\left( {^t_0}X^e \right)^{a}_{\ A} \ \left( {^t_0}X^e \right)^{c}_{\ B} \ {^t}g_{cb} \ {^0}g^{AB}
\\
\nonumber
& =
\left( {^t_0}X^e \right)^{a}_{\ A} \ \left( {^t_0}X^e \right)^{c}_{\ B} \ {^t}g_{cb} \ \mathring{g}^{CD}  \ \left( {^0_0}X^p \right)^{A}_{\ C} \ \left( {^0_0}X^p  \right)^{B}_{\ D}
\\
& = 
{^t_0}X^{a}_{\ C} \ {^t_0}X^{c}_{\ D} \ {^t}g_{cb} \ \mathring{g}^{CD}.
\end{align}
Taking the square root of the determinant of the mixed components of $ \boldsymbol{\underline{\underline{b}}}^e$ gives
\begin{equation}
\label{eq:elastic_Jacobian}
\sqrt{
\det [ ( {^t}b^e )^{a}_{ \ b} ] 
} 
=
\det \left[ ( {^t_0}X^e )^{a}_{\ A} \right] \frac{ \sqrt{ \det [{^t}g_{cb}] } }{\sqrt{ \det [{^0}g_{AB}] }} 
=
\det [ {^t_0}X^{a}_{\ C} ]  \frac{ \sqrt{\det [{^t}g_{cb}] } }{ \sqrt{ \det [\mathring{g}_{CD}] } },
\end{equation}
which, by comparison with the total Jacobian~\eqref{eq:Jacobian_components}, gives
\begin{equation}
{{^t_0}J} = 
\det \left[ {{^t_0}X}^a_{\ A} \right] \dfrac{\sqrt{ \det \left[ {{^t}g_{ab}} \right]}}{\sqrt{ \det \left[ {{^0}g_{AB}} \right]}}
=
\underbrace{ \det [ {^t_0}X^{a}_{\ C} ]  \frac{ \sqrt{\det [{^t}g_{cb}] } }{ \sqrt{ \det [\mathring{g}_{CD}] } }}_{ = \sqrt{
\det [ ( {^t}b^e )^{a}_{ \ b} ] 
} } \underbrace{ \frac{  \sqrt{ \det [\mathring{g}_{CD}] } }{ \sqrt{ \det \left[ {{^0}g_{AB}} \right]}} }_{ = \det [ ( {^0_0}X^p )^{A}_{\ B} ]},
\end{equation}
or, equivalently, 
\begin{equation}
\label{eq:elastic_and_plastic_Jacobians}
J^e = \sqrt{
\det [ ( {^t}b^e )^{a}_{ \ b} ] 
},
\qquad
J^p = \det [ ( {^0_0}X^p )^{A}_{\ B} ].
\end{equation}
The above definition of the elastic and plastic parts of the Jacobian provided by Eq.~\eqref{eq:elastic_and_plastic_Jacobians} are fundamental in cases where isochoric plastic deformations are expected in the context of curvilinear coordinates, i.e., when the Cauchy stress does not entirely depend only on current elastic deformations~\cite{abe2024reconstruction} (this is elaborated on in Section~\ref{subsec:constitutive_equations} and in~\ref{sec:theormodynamics}).

However, while Eq.~\eqref{eq:elastic_Jacobian} offers two methods for computing the elastic part of the Jacobian, neither is practical for linearisation purposes. 
If adopting the first definition, which involves the $\det \left[ ( {^t_0}X^e )^{a}_{\ A} \right]$, it is necessary to compute the derivative $\dfrac{\partial ({^t_0}X^e )^{a}_{\ A}}{\partial {^t_0}X ^{b}_{\ B}}$. While in the case of isotropic elastic behaviour ${^t_0}\underline{\underline{\boldsymbol{X}}}^e$ can be computed by ${^t_0}\underline{\underline{\boldsymbol{X}}}^e = {^t_0}\underline{\underline{\boldsymbol{v}}}^e \cdot {^t_0}\underline{\underline{\boldsymbol{Q}}}^{e, tr}$ (see~\cite{desouza2011computational}), the non-trivial derivative $\dfrac{\partial {^t_0}\underline{\underline{\boldsymbol{v}}}^e}{\partial{^t_0}\underline{\underline{\boldsymbol{X}}}}$ would necessarily require the constitutive relationship and additional kinematic quantities for linearisation. 
On the other hand, if employing the definition based on the total deformation gradient, deriving $\det \left[ \mathring{g}_{AB} \right]$ requires obtaining the derivative of $\dfrac{\partial {^0_0}\underline{\underline{\boldsymbol{X}}}^p}{\partial {^t_0}\underline{\underline{\boldsymbol{X}}}}$, which presents similar challenges to those encountered with the first case.

By exploiting the relationship between $ {^t_0} \underline{\underline{\boldsymbol{b}}}^e$ and the elastic part of the logarithmic strain  ${^t} \underline{\underline{\boldsymbol{\epsilon}}}^e 
\coloneqq
\frac{1}{2} \log {^t} \underline{\underline{\boldsymbol{b}}}^e$, this work employs the following formula to linearise the elastic portion of the Jacobian
\begin{equation}
J^e = \sqrt{
\det [ ( {^t}b^e )^{a}_{ \ b} ] 
}
=
\det [ \sqrt{( {^t}b^e )^{a}_{ \ b} } ] 
=
\det [  \exp ( {^t}\epsilon^e  )^{a}_{ \ b}  ]
=
\exp ( {^t}\epsilon_v^{e} ). 
\end{equation}
with ${^t}\epsilon_v^{e} \coloneqq \text{tr} \left( {^t}\boldsymbol{\underline{\underline{\epsilon}}}^{e} \right).$
As detailed in~\ref{subsec:linearisation_internal_force_vector}, this choice is the one which minimises the extra components with respect to the standard Hencky elasto-plastic material illustrated in the pivotal works~\cite{simo1988frameworkI,simo1988frameworkII}.  
Moreover, these additional components can be directly obtained from the constitutive subroutine, as explained in~\ref{subsec:linearisation_constitutive_equations}.

\section{Stress measures, constitutive equations and forms of equilibrium}
\label{sec:stresses_and_weak_form}

\noindent While the previous section outlined the kinematics of elasto-plastic deformations in non-orthonormal curvilinear coordinates, this section focuses on stress measures and equilibria. 
Specifically, some measures of stress of interest for this manuscript are introduced, and their relationships are stated in general coordinates. 
A qualitative discussion of the adopted constitutive equations for elastically and plastically isotropic materials follows. 
Based on these premises, the strong and weak forms\footnote{
In terms of presentation, this work favours introducing the equations in the strong forms first, followed by the weak forms. 
While it is common in physics and continuum mechanics to present them in the opposite order (see, e.g.,~\cite{sadik2025generalized}), this approach is generally adopted in the context of the FEM (see, for instance,~\cite{oden2006finite}).
} for three-dimensional problems in curvilinear coordinates are introduced in different configurations.

\subsection{Stress measures}
\noindent The Cauchy theorem states that, if the current traction per unit area only depends on the current position $p$, time $t$ and outer normal vector $^t\underline{\boldsymbol{n}}$, then there is a unique measure of stress, the \emph{Cauchy stress}, that depends on the position $p$ and time $t$ such that
\begin{equation}
\label{eq:Cauhcy_theorem}
\int_{\partial ^{t}\mathcal{B}} {^t}t^a \ {^t}\boldsymbol{\underline{g}}_a \ \text{d} {{^t}A} = \int_{\partial ^{t}\mathcal{B}} {^t}\sigma^{a}_{\ b} \ {^t}n^b \ {^t}\boldsymbol{\underline{g}}_a \ \text{d} {{^t}A},
\end{equation}
where ${^t}\sigma^{a}_{\ b}$ are the mixed components of the second-order tensor ${{^t}\underline{\underline{\boldsymbol{\sigma}}}}$.
Additionally, let us consider the Nanson's formula\footnote{
In a geometric view of mechanics, Nanson's formula expresses the pull-back of a contravariant second-order tensor which describes the infinitesimal oriented area via the map $\phi$ (see~\cite{dhas2020geometric}). However, this work is concerned only with the consequences of Nanson's formula for the components of stress tensors. 
} for transforming the oriented areas $^0\underline{\boldsymbol{n}} \ \text{d} {^0A}$ and $^t\underline{\boldsymbol{n}} \ \text{d} {^t A}$ of outer normal vectors $^0\underline{\boldsymbol{n}}$ and $^t\underline{\boldsymbol{n}}$\footnote{
In some ways, unit normals are more naturally represented as covectors rather than vectors (see~\cite{hughes1983mathematical}).
However, treating them as vectors in Eq.~\eqref{eq:Nanson_formula} is correct because the transpose of the deformation gradient provides the necessary metric coefficient to convert them into covectors.
}, that is
\begin{equation}
\label{eq:Nanson_formula}
{{^t_0}\underline{\underline{\boldsymbol{X}}}}^T \cdot {^t\underline{\boldsymbol{n}}} \ \text{d} {^t A}
= 
{^t_0}{J} \ {^0\underline{\boldsymbol{n}}} \ \text{d} {^0 A} \left( \phi^{-1} \left( p \right) \right),
\end{equation}
with the caveat that, in curvilinear components, Eq.~\eqref{eq:def_gradient_transpose_components} for the transpose of the deformation gradient and Eq.~\eqref{eq:Jacobian_components} for the Jacobian must be considered. 

Employing Cauchy's theorem Eq.~\eqref{eq:Cauhcy_theorem} and Nanson's formula Eq.~\eqref{eq:Nanson_formula}, the definition of the first Piola-Kirchoff stress tensor follows, i.e.,
\begin{equation}
\label{eq:first_PK_stress_def}
{^t_0}P^{aA} \left( \phi^{-1} \left( p \right) \right)
\coloneqq 
{^t_0}{J} \left({^t_0}X^{-1} \right)^{A}_{ \ \ b} \ {^t}\sigma^{ab}
\end{equation}
Another possible measure of stress of interest for this manuscript is defined as follows
\begin{equation}
\label{eq:semi_kirchhoff_stress}
\boldsymbol{\underline{\underline{\zeta}}} \coloneqq J^{e} \ {^t}\boldsymbol{\underline{\underline{\sigma}}}.
\end{equation}
As explored in~\cite{bennett2016finite} and based on the initial idea in~\cite{epstein1990energy}, the use of this latter measure of stress in the context of finite-strain elasto-plasticity with associated flow rule is motivated by two features: $(i)$ a constitutive relationship based on $\boldsymbol{\underline{\underline{\zeta}}}$ does not violate the second law of thermodynamics (i.e., the Clausius-Duhem inequality is satisfied, see~\ref{sec:theormodynamics}); $(ii)$ a material described by an elastic relationship between the stress tensor $\boldsymbol{\underline{\underline{\zeta}}}$ and the logarithmic elastic strain ${^t}\boldsymbol{\underline{\underline{\epsilon}}}^{e}$ preserves a stress-free intermediate configuration.  
Reference~\cite{bennett2016finite} also highlights how, if adopting a constitutive relationship between the Kirchhoff stress $\boldsymbol{\underline{\underline{\tau}}} \coloneqq {^t_0}J \ {^t}\boldsymbol{\underline{\underline{\sigma}}}$ and the same measure of elastic strain\footnote{
Constitutive elastic relationships between Kirchhoff stress and logarithmic elastic deformation---often referred to as \emph{Hencky material}~\cite{hencky1933elastic} when linear---were initially and correctly proposed in the context of J-2 plasticity in~\cite{simo1988frameworkI,simo1988frameworkII}. 
This approach has been widely discussed in the literature since it mimics  the small-strain additive elasto-plastic framework. 
However, their direct extension to materials obeying to pressure-dependent yield functions presents the above-discussed issue for the intermediate configuration. 
}, this material preserves an intermediate stress-free configuration only in the case of isochoric plasticity. 

\subsection{Constitutive equations}
\label{subsec:constitutive_equations}
\noindent  This manuscript considers an elastically and plastically isotropic material, with elastic behaviour defined by a stored energy function (hyperelastic material) and plasticity governed by the principle of maximum dissipation. 
In the attempt of not overloading the main part of the manuscript, details on the constitutive equations considered in this work and their relation to the Clausius-Duhem inequality via the Coleman-Noll procedure are reported in~\ref{sec:theormodynamics}. 
This section will focus solely on how the different assumptions regarding the constitutive behaviour impose some constraints and permit other simplifications.

The assumption of isotropic elasticity\footnote{
For the reader interested in elastically anisotropic anelastic solids, reference~\cite{yavari2023direct} discusses the necessity of including the plastic deformation gradient in the structural tensors that describe anisotropic behaviour.
}  implies that the stored energy function governing the constitutive equation can be expressed as a function of the left Cauchy-Green deformation tensor ${^t}\underline{\underline{\boldsymbol{b}}}^e$ (see~\cite{simo1998numerical}), and, more specifically, as a function of its principal invariants, denoted by $I_1^{b^e}, I_2^{b^e}$, and $ I_3^{b^e}$.
By considering the assumption of isotropic plasticity, the need for a constitutive equation for the plastic spin tensor is excluded, as this quantity is set as zero\footnote{For a discussion on anisotropic plasticity and the role of the plastic spin tensor, see~\cite{dafalias1984plastic,dafalias1985plastic} for polycrystalline materials and~\cite{bennett2019anisotropic} for geomaterials.}.

On top of the above hypothesis about isotropy, by assuming that the part of the stored energy function depending on elastic deformations is independent from the other part of the stored energy function relative to kinematic hardening, 
leads to considering a so-called \emph{uncoupled} (or decoupled) material (see, for instance,~\cite{lubliner1972thermodynamic}). 
The assumption of uncoupled material, in turn, is in compliance with considering the principle of maximum dissipation~\cite{hill1948variational,mandel1964contribution} for the evolution equation, which results in associative plasticity\footnote{
As demonstrated in~\cite{collins1997application,collins2002associated}, when considering a coupled material, non-associated flow rules naturally arise and the principle of maximum dissipation cannot be straightforwardly applied. For an extension of uncoupled materials to finite strain elasto-plasticity, see, for instance,~\cite{oliynyk2020finite} and~\cite{pretti2024preserving}.
}.

\subsection{Strong form of balance of rate of linear momentum}
\noindent Balance of rate of linear momentum can be stated at every configuration.
If the original configuration $^{0}\mathcal{B}$ is chosen, then it is often referred to as \emph{Total Lagrangian} (TL), while \emph{Updated Lagrangian} (UL) is the label employed when the current configuration $^{t}\mathcal{B}$ is selected.
In the case where inertia forces are neglected, the UL strong form of the balance of rate of linear momentum is as follows
\begin{equation}
\label{eq:UL_strong_form}
\underline{\textbf{div}} \left( {^t}\boldsymbol{\underline{\underline{\sigma}}} \right) + {^t}\rho \ ^t\boldsymbol{\underline{b}} \approx \boldsymbol{0} \quad \text{on} \ ^{t}\mathcal{B},
\end{equation}
where the components of the divergence with respect to the current position are given by $\underline{\textbf{div}} ( {^t}\boldsymbol{\underline{\underline{\sigma}}} )  = {^t}\sigma^{ab}|_b \ {^t}\boldsymbol{\underline{g}}_a = \left( \dfrac{\partial {^t}\sigma^{ab}}{\partial ^tx^b} + {^t}\sigma^{ac} \ {^t}\gamma^b_{bc}  
+ {^t}\sigma^{bc} \ {^t}\gamma^a_{bc} \right) {^t}\boldsymbol{\underline{g}}_a $ (see Section~\ref{subsec:covariant_derivatives}). 
In Eq.~\eqref{eq:UL_strong_form}, ${^t}\rho$ indicates the current density and $^t\boldsymbol{\underline{b}}$ the body forces per unit mass.
This form of balance of rate of linear momentum is usually complemented by the following Neumann and Dirichlet boundary conditions (BCs)
\begin{alignat}{2}
\label{eq:Neumann_BC_current}
{^t}\sigma^{a}_{\ b} \ {^t}n^b & = {^t}\bar{t}^a \quad & \text{on} \ \partial {^t}\mathcal{B}^N;
\\
\label{eq:Dirichlet_BC_current}
{^t}x^a \left( p \right)  & = \bar{\phi}^{a} \quad & \text{on} \ \partial {^t}\mathcal{B}^D,
\end{alignat}
with the current boundary partitioned as follows $\partial ^{t}\mathcal{B} = \partial {^t}\mathcal{B}^D \cup \partial {^t}\mathcal{B}^N$ and $\partial {^t}\mathcal{B}^D \cap \partial {^t}\mathcal{B}^N = \emptyset$.

To compute the TL form of balance of rate of linear momentum, Eq.~\eqref{eq:gradient_transformation} is recalled, as well as the definition of the first Piola-Kirchhoff stress in Eq.~\eqref{eq:first_PK_stress_def}.
Adopting these formulae (which constitute the base for the Piola identity, see~\cite{hughes1983mathematical}), the TL formulation of the balance of rate of linear momentum can be stated as
\begin{equation}
\label{eq:TL_strong_form}
\underline{\textbf{Div}} \left( {^t_0}\boldsymbol{\underline{\underline{P}}} \right) + {^0}\rho \ ^t\boldsymbol{\underline{B}} \approx \boldsymbol{0} \quad \text{on} \ ^{0}\mathcal{B},
\end{equation}
where the components of the divergence with respect to the initial position are given by $ \underline{\textbf{Div}} ( {^t_0}\boldsymbol{\underline{\underline{P}}} ) = {^t_0}P^{aA}|_A \ {^t}\boldsymbol{\underline{g}}_a = \left( \dfrac{\partial {^t_0}P^{aA}}{\partial ^0x^A} 
+ {^t_0}P^{aA} \ {^0}\Gamma^B_{AB} 
+ {^t_0}P^{bA} \ {^t}\gamma^a_{bc} \ {^t_0}X^{c}_{A} \right) {^t}\boldsymbol{\underline{g}}_a$. 
The initial density ${^0}\rho$ obeys to the principle of mass conservation and the volume transformation via the Jacobian Eq.~\eqref{eq:Jacobian_definition}, resulting in ${^0}\rho \left( \phi^{-1} \left( p \right) \right)  = {^t_0}J \ {^t}\rho  $. The body forces per unit mass are the same as in the UL case, but their dependency is different, i.e., ${^t}B^a \left( \phi^{-1} \left( p \right) \right) = {^t}b^a \left( p \right)$. 
The BCs Eqs.~\eqref{eq:Neumann_BC_current}-\eqref{eq:Dirichlet_BC_current} become as follow
\begin{alignat}{2}
\label{eq:Neumann_BC_original}
{^t_0}P^{aA} \ {^0}n_A & = {^t}\bar{T}^a \qquad & \text{on} \ \partial {^0}\mathcal{B}^N;
\\
\phi^a \left( P \right) & = \bar{\phi}^{a} \qquad & \text{on} \ \partial {^0}\mathcal{B}^D,
\end{alignat}
under the assumption that Neumann boundary transforms as $\phi \left( \partial {^0}\mathcal{B}^N \right) = \partial {^t}\mathcal{B}^N$ and a similar transformation holds for the Dirichlet boundary.
The relationship between the tractions vectors in the current and original configurations follows from Nanson's formula Eq.~\eqref{eq:Nanson_formula}, resulting in 
\begin{equation}
\label{eq:Piola_transformation_applied_tractions}
{^t}\bar{T}^a \ \text{d} {^0 A} \left( \phi^{-1} \left( p \right) \right) = {^t}\bar{t}^a \ \text{d} {^t A}.
\end{equation}

Although the balance of angular momentum is seldom included in the set of equations solved in a FEM program, it provides, in the case of Cauchy continuum, the symmetry of the Cauchy stress ${^t}\sigma^{ab} = {^t}\sigma^{ba}$, and the following symmetry for the first Piola-Kirchhoff stress, that is ${^t_0}P^{aA} \ {^t_0}X^{b}_{A} = {^t_0}P^{bA} \ {^t_0}X^{a}_{A}$.

Eqs.~\eqref{eq:UL_strong_form} and~\eqref{eq:TL_strong_form} are not the only two possible strong forms of equilibrium. 
Other forms of equilibrium can be also considered, either by considering different configurations (see, for instance,~\cite{simo2005stress,menzel2007configurational,coombs2020lagrangian}) or different bases (see, e.g.,~\cite{kohler2001domain}). 
However, once the transformation rules are established, transitioning from one form of equilibrium to another becomes straightforward. 

\subsection{Weak forms of balance of rate of linear momenta}
To achieve a weak statement of the previous strong forms of balance of rate of linear momentum based on an irreducible (i.e., single variable) form, the following infinite-dimensional space is introduced, that is, the \emph{space of configurations} 
${\mathscr{V}_{\bar{\phi}}} \left( {^0}\mathcal{B} \right)$ (see~\cite{simo1998computational}), defined by
\begin{equation}
{\mathscr{V}_{\bar{\phi}}} \left( {^0}\mathcal{B} \right) 
\coloneqq 
\Bigl\{
\phi^{a} : {^0}\mathcal{B} \rightarrow \mathds{R}^{3}
\bigm| 
\phi^{a}|_{\partial {^0}\mathcal{B}^D}
= 
\bar{\phi}^a
\Bigr\},
\end{equation}
where $\bigl( H^1 \left( {^0}\mathcal{B} \right) \bigr)^3 \supset
{\mathscr{V}_{\bar{\phi}}} \left( {^0}\mathcal{B} \right)$. $H^1 \left( {^0}\mathcal{B} \right)$ is the Sobolev space of functions in $H^0 \left( {^0}\mathcal{B} \right)$ with first weak partial derivatives also in $H^0 \left( {^0}\mathcal{B} \right)$, where $H^0 \left( {^0}\mathcal{B} \right)$ is the space of square integrable functions over ${^0}\mathcal{B}$. 
Following a similar definition to the space of configuration, the spaces of the test functions of material \emph{admissible variations} $\mathscr{V}_{0} \left( {^0}\mathcal{B} \right)$, and its current counterpart $\mathscr{V}_{0} \left( {^t}\mathcal{B} \right)$ can be introduced.
Using the divergence theorem, the zero value of the test functions on the Dirichlet boundary and the Neumann BCs Eqs.~\eqref{eq:Neumann_BC_current} and~\eqref{eq:Neumann_BC_original}, the following weak forms hold: find $\phi^a \in {\mathscr{V}_{\bar{\phi}}} \left( {^0}\mathcal{B} \right)$, such that
\begin{multline}
\label{eq:updated_lagrangian_3D}
\mathscr{G} \left( \phi; {^t}\underline{\boldsymbol{w}} \right) \coloneqq
\int_{{^t}\mathcal{B}} 
\left( 
\underline{\textbf{grad}} \left( {^t}\underline{\boldsymbol{w}} \right) 
\boldsymbol{:}
{^t}\underline{\underline{\boldsymbol{\sigma}}} - 
{^t\rho} \ {^t}\underline{\boldsymbol{w}} \cdot {^t} \underline{\boldsymbol{b}} \right) 
\text{d} {^t V} 
- 
\int_{\partial ^{t}\mathcal{B}^N}
{^t}\underline{\boldsymbol{w}} 
\cdot {^t}\bar{\underline{\boldsymbol{t}}} 
\ \text{d} {^t A} = 0, 
\quad \forall \ {^t}\underline{\boldsymbol{w}} \in \mathscr{V}_{0} \left( {^t}\mathcal{B} \right), 
\end{multline}
and, similarly,
\begin{multline}
\mathscr{G} \left( \phi; {^t}\underline{\boldsymbol{W}} \right) \coloneqq
\label{eq:total_lagrangian_3D}
\int_{^{0}\mathcal{B}} 
\left( 
\underline{\textbf{Grad}} \left( {^t}\underline{\boldsymbol{W}} \right)
\boldsymbol{:}
{^t_0}\underline{\underline{\boldsymbol{P}}} 
- 
{^0\rho} \ {^t}\underline{\boldsymbol{W}} \cdot {^t} \underline{\boldsymbol{B}} \right) 
\text{d} {^0 V} 
\\
- \int_{\partial ^{0}\mathcal{B}^N} {^t}\underline{\boldsymbol{W}} \cdot {^t}\bar{\underline{\boldsymbol{T}}} \ \text{d} {^0 A} = 0, 
\quad \forall \ {^t}\underline{\boldsymbol{W}} \in \mathscr{V}_{0} \left( {^0}\mathcal{B} \right),
\end{multline}
where ${^t}W^a \left( \phi^{-1} \left( p \right) \right) = {^t}w^a \left( p \right)$.

\section{Discrete weak form for axisymmetric problems}
The three-dimensional weak forms discussed in the previous section are below particularised to the case of an axisymmetric problem. 
Leveraging their symmetry, these formulations are discretised as a two-dimensional problem, with the out-of-plane components carefully evaluated as needed. 
To enhance clarity, the linearisation of the considered non-linear problems and the algorithmic workflows for the FEM are also outlined. 

\subsection{Axisymmetric assumption}
\begin{figure}
\centering
\includegraphics[width=0.5\textwidth, trim=0cm 0cm 0cm 0cm]{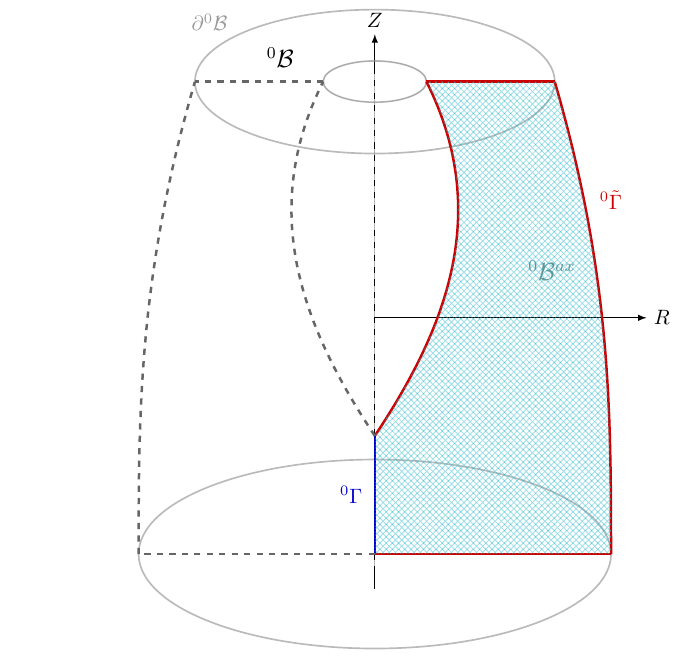}
\caption{Three-dimensional solid ${^0}\mathcal{B}$, obtained by revolving the in-plane cross section ${^0}\mathcal{B}^{ax}$ around the $Z$-axis. Boundaries ${^0}\tilde{\Gamma}$ and ${^0}\Gamma$ of the in-plane section are also illustrated.}
\label{fig:solid_of_revolution}
\end{figure}

Let the three-dimensional body in the reference configuration $^{0}\mathcal{B}$ obtained by the revolution of an arbitrary in-plane cross-section shape (whose domain, denoted by ${^0}\mathcal{B}^{ax}$, lies in the RZ-plane) and therefore symmetric with respect to the $Z$-axis (see Fig.~\ref{fig:solid_of_revolution}).

Let ${^0\Gamma}$ denote the eventual intersection of $ ^{0}\mathcal{B}$ with the $Z$-axis, so that $\Gamma_0 \subset \partial {^0}\mathcal{B}^{ax}$, while ${^0}{\tilde{\Gamma}} \coloneqq \partial{^0}\mathcal{B}^{ax} \setminus {^0\Gamma}$ represents the intersection of the in-plane domain boundary with that of the three-dimensional domain.
This notation is particularly convenient because the initial boundary ${^0}\tilde{\Gamma}$ inherits the Neumann and Dirichlet partitions from $\partial {^0}\mathcal{B}$, i.e., ${^0\tilde{\Gamma}} = {^0\tilde{\Gamma}}_N \cup {^0\tilde{\Gamma}}_D$ (see, for a comparison,~\cite{deparis2004numerical}).

On top of these assumptions on the initial geometry, let us assume that the prescribed motion, the applied tractions and the body forces do not depend on out-of-plane components, i.e.,
\begin{equation}
\begin{dcases}
\begin{alignedat}{3}
\bar{\phi} \left( {^{0}z^1}, {^{0}z^{2}}, {^{0}z^{3}}\right) 
& = 
\bar{\phi}^{ax} \left( {^{0}x^I}, {^{0}x^{III}}\right) 
\qquad 
& \text{on} & \ {^0}\tilde{\Gamma}^D;
\\
{^t}\bar{\boldsymbol{\underline{T}}} \left( {^{0}z^1}, {^{0}z^{2}}, {^{0}z^{3}}\right) 
&
= 
{^t}\bar{\boldsymbol{\underline{T}}}^{ax} \left( {^{0}x^I}, {^{0}x^{III}}\right)
& \text{on} & \ {^0}\tilde{\Gamma}^N;
\\
{^t}\boldsymbol{\underline{B}} \left( {^{0}z^1}, {^{0}z^{2}}, {^{0}z^{3}}\right) 
& = 
{^t}\boldsymbol{\underline{B}}^{ax} \left( {^{0}x^I}, {^{0}x^{III}}\right) 
& \text{on} & \ {^0}\mathcal{B}^{ax},
\end{alignedat}
\end{dcases}
\end{equation}
and do not induce any out-of-plane motion, i.e., if $P \in {^0}\mathcal{B}^{ax}$, it follows that $p = \phi^{ax} \left( P \right) \in {^t}\mathcal{B}^{ax}$. 
Under these conditions, the volume integrals of any quantity $\left( \bullet \right)$ can be expressed as follows
\begin{multline}
\int_{{^0}\mathcal{B}}  \left( \bullet (P) \right) \ \text{d}{^0}V 
=
\int_{{^t}\mathcal{B}}  \left( \bullet \left( \phi(P) \right) \right) \ \text{d}{^t}V 
\\
=
2 \pi
\int_{{^0}\mathcal{B}^{ax}}  \left( \bullet^{ax} (P) \right) \  \text{d}^{0}{V}^{ax}
=
2 \pi
\int_{{^t}\mathcal{B}^{ax}}  \left( \bullet^{ax} \left( \phi^{ax}(P) \right)  \right) \ \text{d}^{t}{V}^{ax}
\end{multline}
with
\begin{equation}
\text{d}^{0}{V}^{ax} 
= 
{^0}{x^I} \ \text{d}^{0}{x^I} \ \text{d}^{0}{x^{III}};
\qquad
\text{d}^{t}{V}^{ax}
=
{^t}{x^i} \ \text{d}^{t}{x^i} \ \text{d}^{t}{x^{iii}}
.
\end{equation}
As for the integral over the surfaces, the outer normal vector ${^0}\underline{\boldsymbol{n}}^{ax}$ to boundary of the in-plane section $\partial {^0}\mathcal{B}^{ax}$ can be decomposed as follows
${^0}\underline{\boldsymbol{n}}^{ax} = {^0}n^{I} \ {^{0}\boldsymbol{\underline{g}}^{I}} + {^0}n^{III} \ {^{0}\boldsymbol{\underline{g}}^{III}}$, which implies that the oriented surfaces can be decomposed into
\begin{equation}
\begin{dcases}
{^0}n^{I} \ \text{d}{^0}A 
= 
{^0}\underline{\boldsymbol{n}}^{ax} \cdot {^{0}\boldsymbol{\underline{g}}^{I}} \ \text{d}{^0}A 
= 
2 \pi \ {^0}{x^{I}} \ \text{d}{^0}{x^{III}} 
\coloneqq 
2 \pi \ {^0}n^{I} \ \text{d}{^0}A^{ax, I};
\\
{^0}n^{III} \ \text{d}{^0}A 
= 
{^0}\underline{\boldsymbol{n}}^{ax} \cdot {^{0}\boldsymbol{\underline{g}}^{III}} \ \text{d}{^0}A 
= 
2 \pi \ {^0}{x^{I}} \ \text{d}{^0}{x^{I}} 
\coloneqq 
2 \pi \ {^0}n^{III} \ \text{d}{^0}A^{ax, III}.
\end{dcases}
\end{equation}
Owing to this equation, the surface integral of quantities $\left( \bullet \right)$ obeying to the Piola transformation, such as the applied tractions (see Eq.~\eqref{eq:Piola_transformation_applied_tractions}), becomes for the axisymmetric case 
\begin{multline}
\int_{\partial{^0}\mathcal{B}}  \left( \bullet (P) \right) \ {^0}\underline{\boldsymbol{n}} \ \text{d}{^0}A 
=
\int_{\partial{^t}\mathcal{B}}  \left( \bullet \left( \phi(P) \right) \right) \ {^t}\underline{\boldsymbol{n}} \ \text{d}{^t}A 
\\
=
2 \pi
\int_{{^0}\tilde{\Gamma}}  \left( \bullet^{ax} (P) \right) \ {^0}\underline{\boldsymbol{n}}^{ax} \  \text{d}^{0}{A}^{ax}
=
2 \pi
\int_{{^t}\tilde{\Gamma}}  \left( \bullet^{ax} \left( \phi^{ax}(P) \right) \ {^t}\underline{\boldsymbol{n}}^{ax} \right) \ \text{d} {^{t}A}^{ax},
\end{multline}
with the current outer normal vector in the axisymmetric case transforming according to Nanson's formula~\eqref{eq:Nanson_formula}.

The assumption of a motion of the axisymmetric kind $\phi^{ax} \left( {^{0}x^I}, {^{0}x^{III}}\right)$ implies also that the displacements do not exhibit any out-of-plane component, i.e.,
\begin{equation}
\label{eq:axisymmetric_displacements}
{^t_0}\boldsymbol{\underline{U}}^{ax} = {^t_0}U^{I} \ {^0} \boldsymbol{\underline{g}}_{I} + {^t_0}U^{III} \ {^0} \boldsymbol{\underline{g}}_{III} = ( {^t_0}U^{ax} )^{\tilde{A}} \ \boldsymbol{\underline{g}}_{\tilde{A}},
\end{equation}
where the indices $\tilde{A} = R, Z = I, III$, and $\tilde{a} = r,z = i, iii$ are introduced for the curvilinear coordinates when the out-of-plane components $\Theta$ and $\theta$ are irrelevant.

Owing to Eq.~\eqref{eq:axisymmetric_displacements}, it can be seen how the in-plane components of the displacement are not scaled by the shifter, that is
\begin{equation}
\label{eq:same_curvilinear_components}
\begin{dcases}
{^t_0}U^{I} = {^t_0}u^{i};
\\
{^t_0}U^{III} = {^t_0}u^{iii},
\end{dcases}
\qquad \text{or} \quad
{^t_0}U^{\tilde{A}} = {^t_0}u^{\tilde{a}}
\end{equation}
since ${^t_0}S^{i}_{\ I} = {^t_0}S^{iii}_{\ III} = 1$ or, more concisely, ${^t_0}S^{\tilde{a}}_{\ \tilde{A}} = \delta^{\tilde{a}}_{\ \tilde{A}}$.
Hence, while the in-plane shifter ${^t_0}S^{\tilde{a}}_{\ \tilde{A}}$ does not change the non-zero in-plane components of the displacements, the three-dimensional $\left( {^t_0}S^{ax} \right)^a_{\ A}$ shifter (given in Eq.~\eqref{eq:axisymmetric_shifter} in the case of axisymmetry) actively contributes to the definitions of the deformation gradient Eq.~\eqref{eq:deformation_gradient_components_material_displ} or its incremental form Eq.~\eqref{eq:increment_deformation_gradient_components_explicit}.

\subsection{Axisymmetric weak forms}
Before showing the equivalence between the axisymmetric weak forms and their three-dimensional counterparts provided in Eqs.~\eqref{eq:updated_lagrangian_3D} and~\eqref{eq:total_lagrangian_3D}, the following space of functions is introduced ${\mathscr{W}_{\bar{\phi}}} \left( {^0}\mathcal{B}^{ax} \right) \coloneqq
\bigl( H^{1} \left( {^0}\mathcal{B}^{ax} \right) \bigr)^2 \cap
{\mathscr{V}_{\bar{\phi}}} \left( {^0}\mathcal{B}^{ax} \right)$, where
\begin{equation}
\label{eq:axisymmetric_space}
{\mathscr{V}_{\bar{\phi}}} \left( {^0}\mathcal{B}^{ax} \right)
\coloneqq 
\biggl\{
\left( \phi^{ax}\right)^{\tilde{a}} : {^0}\mathcal{B}^{ax} \rightarrow \mathds{R}^{2}
\bigm| 
\phi^{i}|_{{^0}{\Gamma}} =  0,
\
\phi^{\tilde{a}}|_{{^0}\tilde{\Gamma}_D}
= 
\bar{\phi}^{\tilde{a}},
\ 
\left( \int_{{^0}\mathcal{B}^{ax}} \dfrac{(\phi^{i})^2}{{^0}x^{I}} \ \text{d}{^0}V \right)^{1/2} < \infty
\biggr\}.
\end{equation}
This space requires the square-integrability of the out-of-plane components arising from the covariant derivative (compare with~\cite{bernardi1999spectral}) and that the parts of the domain which intersect the axis of symmetry must not disrupt this property\footnote{
It is important to note that $^0\Gamma$ represents a line, as it is formed by the intersection of the in-plane domain ${^0}\mathcal{B}^{ax}$ and the $Z$-axis.
Consequently, this line cannot evolve into a higher-order object (e.g., a surface) and the requirement $\phi^{i}|_{{^0}{\Gamma}} =  0$ appearing in Eq.~\eqref{eq:axisymmetric_space} enforces this restriction. 
However, if external factors force an initially small surface in proximity of $^0\Gamma$ (which, in the limit, can be represented by $^0\Gamma$ itself) to move away from the axis of symmetry then the condition $\phi^{i}|_{{^0}{\Gamma}} =  0$ must be disregarded, and $^0\Gamma$ becomes a Neumann boundary with homogeneous conditions. 
This is the case, for instance, of the simulating the Cone Penetration Test (see, for example,~\cite{bird2024implicit}), in which the axis of symmetry can become a contact surface.
}.

Under the assumptions made throughout this section, the equivalence between the three-dimensional forms expressed by Eqs.~\eqref{eq:updated_lagrangian_3D} and~\eqref{eq:total_lagrangian_3D} and their axisymmetric counterparts then hold, i.e.,
\begin{equation}
\mathscr{G} \left( \phi; {^t}\underline{\boldsymbol{w}} \right) = 2 \pi \ \mathscr{G}^{ax} \left( \phi^{ax}; {^t}\underline{\boldsymbol{w}}^{ax} \right), \qquad \mathscr{G} \left( \phi; {^t}\underline{\boldsymbol{W}} \right) = 2 \pi \ \mathscr{G}^{ax} \left( \phi^{ax}; {^t}\underline{\boldsymbol{W}}^{ax} \right),
\end{equation}
where the axisymmetric weak forms are given by: find 
$\left( \phi^{ax} \right)^{\tilde{a}} \in 
{\mathscr{W}_{\bar{\phi}}} \left( {^0}\mathcal{B}^{ax} \right)$ such that
\begin{multline}
\label{eq:updated_lagrangian_2D}
\mathscr{G}^{ax} \left( \phi^{ax}; {^t}\underline{\boldsymbol{w}}^{ax} \right) \coloneqq
\int_{^{t}\mathcal{B}^{ax}} 
\left( 
\underline{\textbf{grad}} \left( {^t}\underline{\boldsymbol{w}}^{ax} \right) 
\boldsymbol{:}
{^t}\underline{\underline{\boldsymbol{\sigma}}} 
- 
{^t\rho} \ {^t}\underline{\boldsymbol{w}}^{ax} \cdot {^t} \underline{\boldsymbol{b}}^{ax} \right) 
\ \text{d} {^t}V^{ax} 
\\
- \int_{^{t}\tilde{\Gamma}^N} 
{^t}\underline{\boldsymbol{w}}^{ax} \cdot {^t}\bar{\underline{\boldsymbol{t}}}^{ax} 
\ \text{d} {^t} A^{ax} = 0, 
\qquad 
\forall \ {^t}\underline{\boldsymbol{w}}^{ax} \in \mathscr{W}_{0} \left( {^t}\mathcal{B}^{ax} \right), 
\end{multline}
or, similarly,
\begin{multline}
\label{eq:total_lagrangian_2D}
\mathscr{G}^{ax} \left( \phi^{ax}; {^t}\underline{\boldsymbol{W}}^{ax} \right) \coloneqq
\int_{^{0}\mathcal{B}^{ax}} 
\left( \underline{\textbf{Grad}} \left( {^t}\underline{\boldsymbol{W}}^{ax} \right) 
\boldsymbol{:}
{^t_0}\underline{\underline{\boldsymbol{P}}}
- 
{^0\rho} \ {^t}\underline{\boldsymbol{W}}^{ax} \cdot {^t} \underline{\boldsymbol{B}}^{ax} \right) 
\ \text{d}{^0}V^{ax}
\\
- 
\int_{ {^0}\tilde{\Gamma}^N} {^t}\underline{\boldsymbol{W}}^{ax} \cdot {^t}\bar{\underline{\boldsymbol{T}}}^{ax} \ \text{d} {^0}A^{ax} = 0, \qquad \forall \ {^t}\underline{\boldsymbol{W}}^{ax} \in \mathscr{W}_{0} \left( {^0}\mathcal{B}^{ax} \right). 
\end{multline}

\subsection{Discretisation}
\label{subsec:discretisation}
Let the in-plane domain $^0\mathcal{B}^{ax}$ be decomposed into $n^{els}$ elements belonging to the set $^{h}\mathcal{T} \left( {^0}\mathcal{B}^{ax} \right)$, with the arbitrary element denoted by $T$, such that the discrete domain is given by the union of these elements, i.e., $^h({^0}\mathcal{B}^{ax}) \coloneqq \cup_{T \in ^{h}\mathcal{T} \left( {^0}\mathcal{B}^{ax} \right)} T$.
By introducing the space of polynomials of order $k$ over the element $T$ (denoted by $P^{k}\left( T \right)$), the finite element dimensional subspace ${^h}{\mathscr{W}_{\bar{\boldsymbol{\underline{U}}}}} \left( {^0}\mathcal{B}^{ax} \right)$ of ${\mathscr{W}_{\bar{\phi}}} \left( {^0}\mathcal{B}^{ax} \right)$ can be defined as follows
\begin{multline}
^h{\mathscr{W}_{\bar{\boldsymbol{\underline{U}}}}} \left( {^0}\mathcal{B}^{ax} \right)
\coloneqq 
\bigl\{ 
{^h}\mleft( \phi^{ax} \mright)^{\tilde{a}} 
= 
{^h}\mleft({^0}x^{\tilde{a}} \mright) 
+ 
\delta^{\tilde{a}}_{\ \tilde{A}} \ {^h}\mleft( {^t_0}U^{ax} \mright)^{\tilde{A}} 
\in 
{\mathscr{W}_{\bar{\phi}}} \left( {^0}\mathcal{B}^{ax} \right) : 
\\
{^h}\mleft( {^t_0}{\boldsymbol{\underline{U}}}^{ax} \mright)|_T \in P^{k}\left( T \right), \forall T \in {^h}\mathcal{T} \left( {^0}\mathcal{B}^{ax} \right)
\bigr\}.
\end{multline}
Following a similar definition for the current configuration (denoted by $^h{\mathscr{W}_{\bar{\boldsymbol{\underline{u}}}}} \left( {^t}\mathcal{B}^{ax} \right)$)
it is possible to express the discrete test and trial functions via the the shape functions $\Uppsi_{\mathsf{N}} (P)$ and $\uppsi_{\mathsf{N}} \left( \phi^{ax} \left( P \right) \right)$  centred at the arbitrary node $\mathsf{N}$
\begin{alignat}{2}
{^h}\mleft( {^t}w^{ax} \left( \phi^{ax} \left( P \right) \right) \mright)^{\tilde{a}} & = {^t}\mathrm{w}^{{\tilde{a}}\mathsf{N}} \ \uppsi_\mathsf{N}\left( \phi^{ax} \left( P \right) \right);
\qquad 
{^h}\mleft( {^t{W}^{ax}} \left( P \right) \mright)^{{\tilde{a}}} & = {^t}\mathrm{W}^{{\tilde{a}}\mathsf{N}}\ \Uppsi_\mathsf{N}\left( P \right);
\\
{^h}\mleft( {^t_0}{u}^{ax} \left( \phi^{ax} \left( P \right) \right) \mright)^{\tilde{a}} & = {^t_0}\mathrm{u}^{{\tilde{a}}\mathsf{N}} \ \uppsi_\mathsf{N}\left( \phi^{ax} \left( P \right) \right);
\qquad
{^h}\mleft( {^t_0}{U}^{ax} \left( P \right) \mright)^{\tilde{A}} & = {^t_0}\mathrm{U}^{{\tilde{A}}\mathsf{N}} \ \Uppsi_\mathsf{N}\left(  P \right),
\end{alignat}
where the summation on shape functions' nodes $\mathsf{N},\mathsf{M},\mathsf{O} \dots$ has no covariant-contravariant meaning, but a similar notation for these summations is kept for consistency.
It can also be noticed that the discrete function ${^h}\mleft( {^t_0}{\underline{\boldsymbol{u}}}^{ax} \mright)$ and ${^h}\mleft( {^t_0}{\underline{\boldsymbol{U}}}^{ax} \mright)$ are defined is the RZ-plane, which means that the displacements can be computed as a difference between position vectors, as also indicated in Eq.~\eqref{eq:displacements_Cartesian_components} and discussed in Section~\ref{subsec:deformation_gradient_on_displacements}.
Hence, computing a solution for the motion ${^h}\mleft(\phi^{ax} \left( P\right)\mright)^{\tilde{a}}$ or any of the displacements ${^h}\mleft({^t_0}{\underline{\boldsymbol{u}}}^{ax} \left( \phi^{ax} \left( P \right) \right) \mright)$ and ${^h}\mleft( {^t_0}{\underline{\boldsymbol{U}}}^{ax} \left( P \right) \mright)$ has no appreciable consequences since their domains of definitions are flat. 
Moreover, in light of Eq.~\eqref{eq:same_curvilinear_components}, there are no differences between the numerical values of the coefficients interpolating ${^h}\mleft( {^t_0}{u}^{ax} \left( \phi^{ax} \left( P \right) \right) \mright)^{\tilde{a}}$ and ${^h}\mleft( {^t_0}{U}^{ax} \left( P \right) \mright)^{\tilde{A}}$, but keeping a difference in notation between these facilitates their distinction, especially when covariant derivatives are computed (see Eq.~\eqref{eq:gradient_transformation} and discussion in Section~\ref{subsec:algorithmic_workflow}).
In light of the introduced discretisations, the weak forms Eq.~\eqref{eq:updated_lagrangian_2D} and~\eqref{eq:total_lagrangian_2D} now read
\begin{multline}
\label{eq:updated_lagrangian_discrete_equilibrium}
\mathscr{G}^{ax} 
\left( {^h}\mleft({^t_0}\boldsymbol{\underline{u}}^{ax}\mright); {^h}\mleft({^t}\boldsymbol{\underline{w}}^{ax}\mright) \right) =
{^t}\mathbf{\underline{w}}^{\mathsf{N}} \cdot \Bigl(
\int_{^{t}\mathcal{B}^{ax}} 
\left( 
\underline{\textbf{grad}} \left( \uppsi_\mathsf{N} \right) 
\cdot
{^t}\underline{\underline{\boldsymbol{\sigma}}} 
- 
{^t\rho} \ \uppsi_\mathsf{N} \ {^t}\underline{\boldsymbol{b}}^{ax} \right) 
\ \text{d} {^t}V^{ax} 
\\
- \int_{^{t}\tilde{\Gamma}^N} 
\uppsi_\mathsf{N} \ {^t}\bar{\underline{\boldsymbol{t}}}^{ax}
\ \text{d} {^t} A^{ax} \Bigr) = 0, 
\qquad 
\forall \ {^h}\mleft({^t}\boldsymbol{\underline{w}}^{ax} \mright) \in {^h}\mathscr{W}_{0} \left( {^t}\mathcal{B}^{ax} \right), 
\end{multline}
and
\begin{multline}
\label{eq:total_lagrangian_discrete_equilibrium}
\mathscr{G}^{ax} 
\left( {^h}\mleft( {^t_0}\boldsymbol{\underline{u}}^{ax} \mright); {^h}\mleft({^t}\boldsymbol{\underline{W}}^{ax}\mright) \right)
=
{^t}\mathbf{\underline{W}^{\mathsf{N}}} \cdot \Bigl(
\int_{^{0}\mathcal{B}^{ax}} 
\left(
{^t_0}\underline{\underline{\boldsymbol{P}}} \cdot
 \underline{\textbf{Grad}} \left( \Uppsi_\mathsf{N} \right) 
- 
{^0\rho} \ \Uppsi_\mathsf{N} \ {^t}\underline{\boldsymbol{B}}^{ax} \right) 
\ \text{d}{^0}V^{ax}
\\
- 
\int_{ {^0}\tilde{\Gamma}^N} \Uppsi_\mathsf{N} \ {^t}\bar{\underline{\boldsymbol{T}}}^{ax} \ \text{d} {^0}A^{ax}
\Bigr)
= 0, \qquad \forall \ {^h}\mleft({^t}\boldsymbol{\underline{W}}^{ax}\mright) \in {^h}\mathscr{W}_{0} \left( {^0}\mathcal{B}^{ax} \right).
\end{multline}
To reduce the notation in the following sections, let us defined the residual $\mathbf{r}$ relative to the above functionals as
\begin{multline}
\mathscr{G}^{ax} 
\left( {^h}\mleft( {^t_0}\boldsymbol{\underline{u}}^{ax} \mright); {^h}\mleft({^t}\boldsymbol{\underline{W}}^{ax}\mright) \right)
\coloneqq
{^t}\mathbf{\underline{W}^{\mathsf{N}}} \cdot \mathbf{r}_{\mathsf{N}} \left( {^h}\mleft( {^t_0}\boldsymbol{\underline{U}}^{ax} \mright) \right)
\\
= 
{^t}\mathbf{\underline{w}^{\mathsf{N}}} \cdot \mathbf{r}_{\mathsf{N}} \left( {^h}\mleft( {^t_0}\boldsymbol{\underline{u}}^{ax} \mright) \right)
\coloneqq
\mathscr{G}^{ax} 
\left( {^h}\mleft( {^t_0}\boldsymbol{\underline{u}}^{ax} \mright); {^h}\mleft({^t}\boldsymbol{\underline{w}}^{ax}\mright) \right).
\end{multline}
The non-linearity and path-dependency of the equations motivate their formulation in an incremental form, leading to the introduction of a (pseudo) time-discretisation. 
This manuscript adopts a Backward-Euler time segmentation of the variables between the current time $t$ and the previously converged time $t_n$.

\subsection{Linearisation}
To solve the non-linear problems~\eqref{eq:updated_lagrangian_discrete_equilibrium} and~\eqref{eq:total_lagrangian_discrete_equilibrium}, an iterative procedure is required. This manuscript adopts a Newton-Raphson (NR) iterative process, which involves linearising the above equations at the solution $\left( {^h}\mleft( {^t_0}\check{\boldsymbol{\underline{u}}}^{ax} \mright); {^h}\mleft({^t}\boldsymbol{\underline{w}}^{ax}\mright) \right)$ and setting these to be zero, i.e.,
\begin{multline}
\label{eq:NR_updated_lagrangian}
0 = \mathscr{G}^{ax} 
\left( {^h}\mleft( {^t_0}\check{\boldsymbol{\underline{u}}}^{ax} \mright); {^h}\mleft({^t}\boldsymbol{\underline{w}}^{ax}\mright) \right)  
\approx
\mathscr{G}^{ax} 
\left( {^h}\mleft( {^t_0}\boldsymbol{\underline{u}}^{ax} \mright); {^h}\mleft({^t}\boldsymbol{\underline{w}}^{ax}\mright) \right) 
\\
+
\delta \left( \mathscr{G}^{ax} 
\left( 
{^h}\mleft( {^t_0}\boldsymbol{\underline{u}}^{ax} \mright); {^h}\mleft({^t}\boldsymbol{\underline{w}}^{ax}\mright) \right)  
\right)
\cdot
\delta {^h}\mleft( {^t_0}\boldsymbol{\underline{u}}^{ax} \mright).
\end{multline}
This procedure is repeated for $\left( k+1 \right)$ iterations, until an adequate convergence criterion is met.

The variation corresponding to the second term on the right-hand side (RHS) of the above equation involves a covariant derivative again (see~\cite{hughes1983mathematical}), and it is fully detailed in~\ref{sec:linearisation}.

\subsection{Algorithmic workflow}
\label{subsec:algorithmic_workflow}
To bridge the gap between the discussion of the axisymmetric problem carried out so far and the development of the corresponding FEM program, Algorithms~\ref{algorithm:workflow_UL} and~\ref{algorithm:workflow_TL} are introduced, where only the steps that depart more significantly from a Cartesian coordinate code are detailed. 
Building on the same idea to make the material more accessible,~\ref{sec:cylindrical_coordinate_values} reports an unambiguous list of quantities necessary for an axisymmetric program.

\SetNlSty{texttt}{[}{]UL:}
\SetAlgoNlRelativeSize{-1}
\begin{algorithm}[]
\small{
\setstretch{1.0}
$t \leftarrow$ 0 \tcp*[r]{0 the time} \
\While{$t \leq T$}{
$\left( k + 1 \right)  \leftarrow$ 1 \tcp*[r]{initiate the NR counter} \

\While{Outer NR}{ 
\If{$t > 0$}{
${^t_0}\mathrm{u}^{\tilde{a}\mathsf{N}}|^{(k+1)} 
= 
{^t_0}\mathrm{u}^{\tilde{a}\mathsf{N}}|^{(k)}
-
\left( ( \mathrm{K}|^{(k)} )^{\tilde{a}\mathsf{N}}_{ \quad \ \tilde{b} \mathsf{M}} \right)^{-1}  ( \mathrm{r}|^{(k)} )^{\tilde{b} \mathsf{M}}  
$
\tcp*[r]{solve system} 
}
el\_count $\leftarrow$ 1\;
 \ForEach(
 \tcp*[f]{loop over elements}
 ){
 element $T \in {^h}\mathcal{T} \left( {^t}\mathcal{B}^{ax} \right)$
 }{
gp\_count $\leftarrow$ 1\;
    \ForEach(
    \tcp*[f]{loop over gps}
    ){
    gp
    }{
    \setstretch{1.3}
     $\dfrac{\partial {^t}x^{\tilde{a}}}{\partial \xi^{\tilde{j}}} = \dfrac{\partial \uppsi_{\mathsf{N}}}{\partial \xi^{\tilde{j}}} \left( {^0}x^{\tilde{a}\mathsf{N}} + {^t_0}\mathrm{u}^{\tilde{a}\mathsf{N}} \right)$
     \tcp*[r]{der. of current pos. wrt parent coord. ${\tilde{j}} = R,Z$} 
     \
    $\dfrac{\partial \uppsi_{\mathsf{N}}}{\partial {^t}x^{\tilde{a}}} = \dfrac{\partial \uppsi_{\mathsf{N}}}{\partial \xi^{\tilde{j}}} \left( \dfrac{\partial {^t}x^{\tilde{a}}}{\partial \xi^{\tilde{j}}} \right)^{-1}$
    \tcp*[r]{update derivative of shape functions}
    \
    ${^t}x^i = \uppsi_{\mathsf{N}} \ {^t_0}\mathrm{u}^{i\mathsf{N}}|^{(k+1)} + {^0}x^I$\tcp*[r]{update radial position}
    \
    $(^t\gamma^{ax} )^{ii}_{i \ ii} = \dfrac{1}{{^t}x^i}$
    \tcp*[r]{update non-zero Christoffel symbol val.}
    \
    $(\delta^a_{\tilde{b}} \ \uppsi_\mathsf{N} )|_{c} = \delta^a_{\tilde{b}} \dfrac{\partial \uppsi_{\mathsf{N}}}{\partial {^t}x^c} + \delta^d_{\tilde{b}} \  {^t}\gamma^{a}_{dc} \ \uppsi_{\mathsf{N}}$ 
    \tcp*[r]{update gradient of shape functions}
    \
    $({^{t_n}_t}S^{ax})^{\mathcal{A}}_{\ \ b} = \text{diag}(1,({^tx^i})/({^{t_n}x^\mathcal{I}}),1)$\tcp*[r]{inverse of incr. shifter}
    \
    $\left( ({^t_{t_n}X^{ax}} )^{-1} \right)^{\mathcal{A}}_{\ \ a} = ({^{t_n}_t}S^{ax})^{\mathcal{A}}_{\ \ b} \left(\delta^{b}_{a} - (\delta^b_{\tilde{c}} \ \uppsi_\mathsf{N} )|_{a} \ {^t_{t_n}} \mathrm{u}^{\tilde{c}\mathsf{N}}|^{(k+1)} \right)$
    \tcp*[r]{inverse of incr. def. grad.}
    \
    $({^t_{0}X^{ax}})^{a}_{\ A}  = ({^t_{t_n}X^{ax}} )^{a}_{\ \mathcal{A}} \ ({^{t_n}_0X^{ax}} )^{\mathcal{A}}_{\ A}$\tcp*[r]{def. grad.}
    \
    $^t_{0}J = \det [({^t_{0}X^{ax}})^{a}_{\ A} ] {^t}x^i / {^0}x^I$
    \tcp*[r]{Jacobian}
    \
    ${^t}g^{ax}_{ab} = \text{diag}(1,({^tx^i})^2,1)$\tcp*[r]{current metric coefficients}
    \
    $\left( ({^t_{t_n}X^{ax}} )^{T} \right)^{\mathcal{A}}_{\ \ a} = {^t}g^{ax}_{ab} \  ({^t_{t_n}X^{ax}} )^{b}_{\ \mathcal{B}} \ ({^{t_n}}g^{ax})^{\mathcal{AB}}$\tcp*[r]{transpose of incr. def. grad.}
    \
    $({^t}b^{e, tr} )^{a}_{ \ b} = ({^t_{t_n}X^{ax}})^{a}_{\ \mathcal{A}} \ ^{t_n}b^{\mathcal{A}}_{ \ \mathcal{B}} \ \left( ({^t_{t_n}X^{ax}} )^{T} \right)^{\mathcal{B}}_{\ \ b}$
    \tcp*[r]{trial elastic strain}
    \
    $^{t}z^{tr} = {^{t_n}}z$\tcp*[r]{trial strain-like internal variable}
    \tikzmk{A}\SetKwFunction{FMain}{Constitutive relationship}
    \SetKwProg{Fn}{Function}{}{}
    \Fn{\FMain{$( b^{e, tr} )^{a}_{\ b}, z^{tr}$}
    }{
    \textbf{return:} $\mathtt{(} \zeta_{a}^{\ b}, (b^{e} )^{a}_{ \ b}, z,  J^{e}, \left( {^t}\mathcal{D}^{alg} \right)^{\ b \ f}_{a \ e}, \dfrac{\partial \left( {^t}\epsilon^e \right)^{a}_{\ b}}{\partial \left( {^t}\epsilon^{e, tr} \right)^{c}_{ \ d} }\mathtt{)} $ \tikzmk{B}\boxit{gray}
        }
    $J^p = J (J^e)^{-1}$ \tcp*[r]{plastic Jacobian}
    \setstretch{1}
    gp\_count ++
 }
  el\_count ++
 }
 $\left( k + 1 \right) $ ++
 }
 $t$ ++ 
 }
 }
\caption{Algorithmic workflow for the UL formulation}
\label{algorithm:workflow_UL}
\end{algorithm}

In Algorithms~\ref{algorithm:workflow_UL} and~\ref{algorithm:workflow_TL}, it was decided to select the contravariant components of the displacements and the mixed components of the symmetric left Cauchy-Green deformation tensor $^t\underline{\underline{\boldsymbol{b}}}$ and the stress tensor $\underline{\underline{\boldsymbol{\zeta}}}$. 
Programs based on different choices are also possible, and the relationships presented throughout Sections~\ref{sec:kinematic}-\ref{sec:stresses_and_weak_form} should help the reader clarify the relationships among the considered variables and extend Algorithms~\ref{algorithm:workflow_UL} and~\ref{algorithm:workflow_TL} to other options.

Gaussian quadrature rule is considered in Algorithms~\ref{algorithm:workflow_UL} and~\ref{algorithm:workflow_TL} for numerical integration of the quantities appearing in Eqs.~\eqref{eq:updated_lagrangian_discrete_equilibrium} and~\eqref{eq:total_lagrangian_discrete_equilibrium}. The generical Gauss Point is indicated by $gp$.

As standard in TL formulations, the partial spatial derivatives of the shape functions, evaluated with respect to the initial position, can be computed outside the NR loop (see Algorithm~\ref{algorithm:workflow_TL}). 
Conversely, in the UL formulation, the spatial partial derivative with respect to the current position must be updated within this loop (see Algorithm~\ref{algorithm:workflow_UL}). 
However, in both cases, the covariant derivative of the test functions in the current configuration (see Eq.~\ref{eq:gradient_transformation}) requires the current Christoffel symbols $^t\gamma^a_{bc}$ and, as such, must be updated within the NR loop. 

Additionally, to explore more options, the UL formulation computes the deformation gradient incrementally using Eq.~\eqref{eq:inv_incr_deformation_gradient}, whereas the TL form computes the total deformation gradient directly using Eq.~\eqref{eq:deformation_gradient_components_material_displ}.
These formulations are not constrained by these choices, and Algorithms~\ref{algorithm:workflow_UL} and~\ref{algorithm:workflow_TL} are also designed with the purpose of comparing these different approaches.

\SetNlSty{texttt}{[}{]TL:}
\SetAlgoNlRelativeSize{-1}
\begin{algorithm}[]
\small{
\setstretch{1.4}
el\_count $\leftarrow$ 1\;
 \ForEach(
 \tcp*[f]{loop over elements}
 ){
 element $T \in {^h}\mathcal{T} \left( {^t}\mathcal{B}^{ax} \right)$
 }{
gp\_count $\leftarrow$ 1\;
    \ForEach(
    \tcp*[f]{loop over gps}
    ){
    gp
    }{ 
${^0}g^{ax}_{AB} = \text{diag}(1,({^0x^I})^2,1)$\tcp*[r]{initial metric coefficients}
\
$\dfrac{\partial {^0}x^{\tilde{A}}}{\partial \xi^{\tilde{j}}} = \dfrac{\partial \Uppsi_{\mathsf{N}}}{\partial \xi^{\tilde{j}}} \ {^0}x^{\tilde{A}\mathsf{N}}$;
$\qquad \qquad \qquad \qquad$
\tcp*[f]{der. of original pos. wrt parent coord. ${\tilde{j}} = R,Z$} 
\\
$\dfrac{\partial \Uppsi_{\mathsf{N}}}{\partial {^0}x^{\tilde{A}}} = \dfrac{\partial \Uppsi_{\mathsf{N}}}{\partial \xi^{\tilde{j}}} \left( \dfrac{\partial {^0}x^{\tilde{A}}}{\partial \xi^{\tilde{j}}} \right)^{-1}$
\tcp*[r]{compute derivative of shape functions}
    gp\_count ++
 }
  el\_count ++
 }
\setstretch{1.0}
$t \leftarrow$ 0 \tcp*[r]{0 the load-step} \
\While{$t \leq T$}{
$\left( k + 1 \right)  \leftarrow$ 1 \tcp*[r]{initiate the NR counter} \
\While{Outer NR}{
Steps \texttt{\scriptsize{[5]UL-[8]UL}}\;
 \ForEach(
 \tcp*[f]{loop over elements}
 ){
 element $T \in {^h}\mathcal{T} \left( {^0}\mathcal{B}^{ax} \right)$
 }{
gp\_count $\leftarrow$ 1\;
    \ForEach(
    \tcp*[f]{loop over gps}
    ){
    gp
    }{
    Steps \texttt{\scriptsize{[14]UL-[15]UL}}\;
    \setstretch{1.4}
    $(\delta^a_{\tilde{b}} \ \Uppsi_\mathsf{N} )|_{A} = \delta^a_{\tilde{b}} \dfrac{\partial \Uppsi_{\mathsf{N}}}{\partial {^0}x^A} + \delta^d_{\tilde{b}} \  {^t}\gamma^{a}_{dA} \ \Uppsi_{\mathsf{N}}$
    \tcp*[r]{compute gradient of shape functions}
    \
    $({^{t}_0}S^{ax})^{a}_{\ \ B} = \text{diag}(1,({^{0}x^{I}})/({^tx^i}),1)$\tcp*[r]{total shifter}
    \
    $\left( {^t_{0}X^{ax}}  \right)^{a}_{\ A} = ({^{t}_0}S^{ax})^{a}_{\  B} \left(\delta^{B}_{C} + ({^{0}_t}S^{ax})^{B}_{\ \ b} \ (\delta^b_{\tilde{c}} \ \Uppsi_\mathsf{N} )|_{A} \ {^t_{0}} \mathrm{u}^{\tilde{c}\mathsf{N}}|^{(k+1)} \right)$
    \tcp*[r]{def. grad.}
    \ 
    Steps \texttt{\scriptsize{[20]UL-[21]UL}}\;
    $\left( ({^t_{0}X^{ax}} )^{T} \right)^{A}_{\ \ a} = {^t}g^{ax}_{ab} \  ({^t_{0}X^{ax}} )^{b}_{\ \mathcal{B}} \ ({^{0}}g^{ax})^{AB}$\tcp*[r]{transpose of  def. grad.}
    \
    $( {^t}b^{e, tr} )^{a}_{ \ b} = ({^t_{0}X^{ax}})^{a}_{\ A} \ \ \left( ({^t_{0}X^{ax}} )^{T} \right)^{A}_{\ \ b}$
    \tcp*[r]{trial elastic strain}
    \
    Steps \texttt{\scriptsize{[23]UL-[26]UL}}\;
    gp\_count ++
 }
  el\_count ++
 }
 $\left( k + 1 \right) $ ++
 }
 $t$ ++
 }
 }
\caption{Algorithmic workflow for the TL formulation}
\label{algorithm:workflow_TL}
\end{algorithm}

\section{Numerical example: three-dimensional vs axisymmetric formulation for a thick-walled cylinder}

\begin{figure}
\centering
\includegraphics[width=0.6\textwidth, trim=0cm 0cm 0cm 0cm]{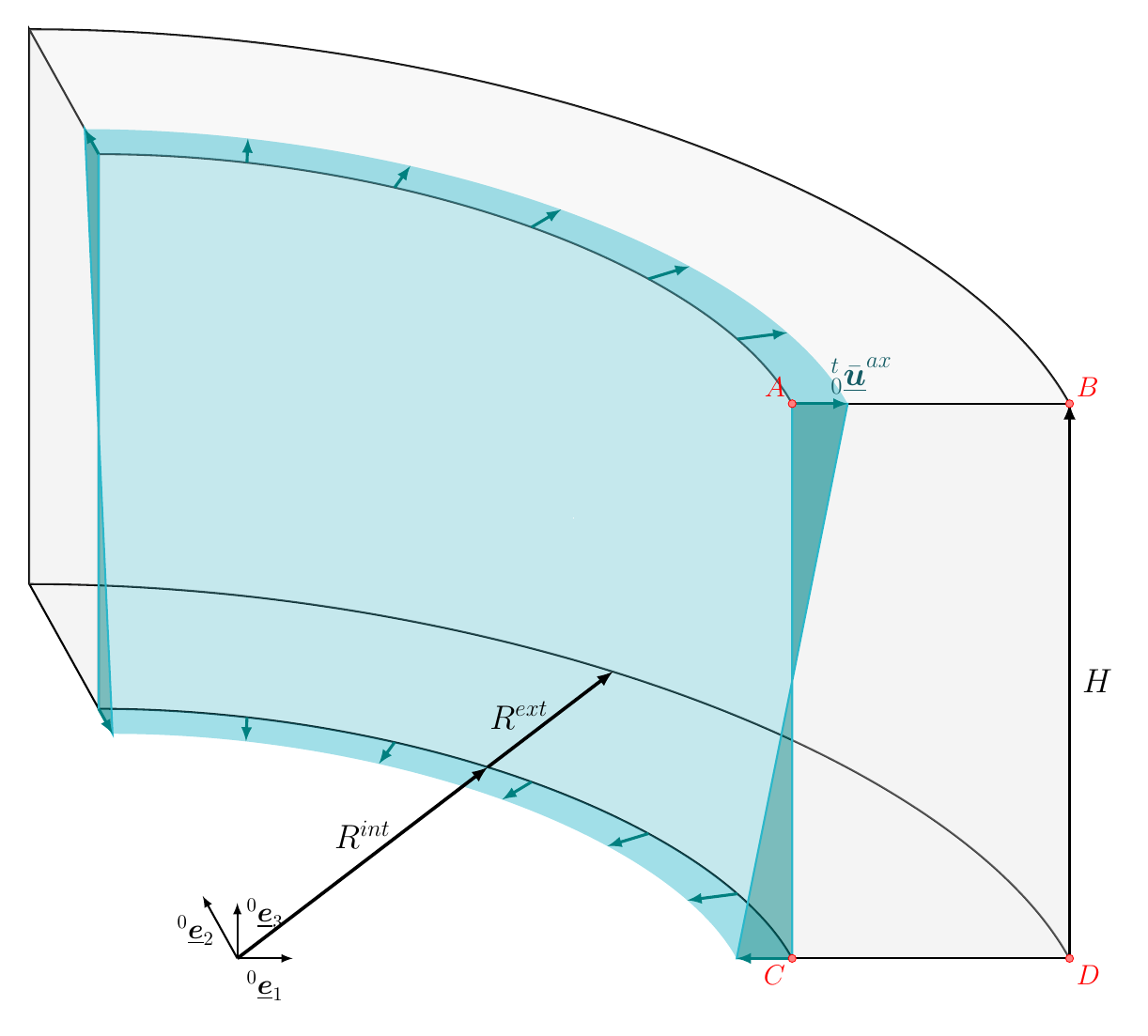}
\caption{Three-dimensional setup of the considered hollow cylinder. Displacements varying along $Z$ on the internal wall are applied.}
\label{fig:setup}
\end{figure}

\paragraph*{Example scope}
The purpose of the example illustrated in Fig.~\ref{fig:setup} is twofold: $(i)$ to compare the results coming from the proposed axisymmetric formulation with those obtained from a three-dimensional analysis; and $(ii)$ to check that the linearisation proposed in~\ref{sec:linearisation} is correct.
The verification of the axisymmetric analysis against another implementation is necessary because an analytical solution is not available for the material under consideration. 
Analytical solutions for axisymmetric problems in the context of finite strain rely on the incompressibility constraint for the volume, as in the case of elastic rubber-like material (see~\cite{haughton1979bifurcation}) or in the case of J-2 plasticity, where the elastic volume contribution is neglected (see~\cite{steigmann2023course}).

\begin{table}
\centering
\small{
\begin{threeparttable} 
\begin{tabular}{ccc} 
\headrow
\multicolumn{3}{c}{Parameter Settings}
\\ \hiderowcolors 
\hline
\multirow{2}{*}{Discretisation}
& $n_{els}^R$, $n_{els}^\Theta${}\textsuperscript{\textdagger}, $n_{els}^Z$ & $5$, $20$\textsuperscript{\textdagger}, $10$\\
& $T$ & $30$\\
\hline 
\multirow{3}{*}{\splitcellc[]{Geometry \\ \& BCs}}
& $R^{int}$, $R^{ext}$ & $10$ m, $15$ m \\
& $H$ & $10$ m \\
& $({^t_0}\bar{u}^i)^{ax}$ & $1$ m\\
\hline
\multirow{3}{*}{Material Parameters}
& $E, \, \nu$ & $1.375 \times 10^9$ Pa, $ 0.375$ \\
& $\kappa$, $\alpha$, $m$ & $0$, $-$, $1$\\
& $H$, $^{0}p_c$ & $765\times 10^6$ Pa, $24 \times 10^7$ Pa\\
\hline
\end{tabular}
\begin{tablenotes}
\item {}\textsuperscript{\textdagger} Out-of-plane discretisation is considered only for the three-dimensional analysis.
\end{tablenotes}
\caption{Summary of the parameters considered in the numerical example.}
\label{table:parameters_setups}
\end{threeparttable}
}
\end{table}

\paragraph*{Setups}
The initial setup for the three-dimensional analysis consists of a quarter of a thick-walled cylinder (see Fig.~\ref{fig:setup} for its illustration) with a radial prescribed displacement on the internal surface. 
This applied Dirichlet BCs varies linearly along the cylinder's height according to the law
\begin{equation}
{^t_0\bar{\underline{\boldsymbol{u}}}}^{ax} 
=
\begin{bmatrix}
( {^t_0}\bar{u}^1)^{ax} = ({^t_0}\bar{u}^i)^{ax} \cos (^0x^{II}) (^0x^{III}-\frac{H}{2}) \frac{2}{H}
\\
( {^t_0}\bar{u}^2)^{ax} = ({^t_0}\bar{u}^i)^{ax} \sin (^0x^{II}) (^0x^{III}-\frac{H}{2}) \frac{2}{H}
\\
( {^t_0}\bar{u}^3)^{ax} = 0
\end{bmatrix}.
\end{equation}
The magnitude of the displacements is linearly ramped up to its final value $({^t_0}\bar{u}^i)^{ax}$ in $T$ equally-divided time-steps (the input parameters used in these numerical simulations are reported in Table~\ref{table:parameters_setups}).
Rollers are placed on the lateral, bottom, and top surfaces, while homogeneous Neumann boundary conditions are applied to the external wall. 
The discretisation of the cylinder is achieved by twenty-noded serendipity hexagonal elements\footnote{It is well acknowledged in the literature that the use of serendipity elements on non-affine grids does not lead to optimal convergence rates with mesh refinement (see~\cite{boffi2013mixed}). 
However, since the current work does not explore this aspect of convergence, the use of serendipity elements to compare results remains valid. 
At most, this shortcoming of serendipity elements in the three-dimensional case is a further point in favour of the proposed axisymmetric formulation, which does not involve a non-affine grid.} with $3 \times 3 \times 3$ Gaussian integration. This discretisation is obtained by uniformly segmenting the RZ-plane using  $n_{els}^R \times n_{els}^Z$ elements and the right angle using $n_{els}^\Theta$ elements.
For the axisymmetric initial setup, only a two-dimensional planar section along the radius of the three-dimensional analysis is considered. 
In this case, eight-noded serendipity quadrilateral elements with full $3 \times 3$ Gauss points are adopted.

The NR error for the examples under consideration (and plotted for some time-steps in Fig.~\ref{fig:NR_convergence}) is as follows 
\begin{equation}
err^{(k+1)} \coloneqq \dfrac{|| \mathbf{r}^{(k+1)} ||}{|| \mathbf{r}^{(1)}||} \leq tol
\end{equation}
with $tol = 10^{-8}$ and $|| \left( \bullet \right) ||$ being the Euclidean norm of $\left( \bullet \right)$.

The material considered for these simulations is an extension of a decoupled, compressible Neo-Hookean material, with a hardening law that, depending on the considered parameters, can be either linear or exponential, i.e., 
\begin{equation}
\label{eq:strain_energy}
 \bar{\rho} \ \psi \left( {^t}\boldsymbol{\underline{\underline{b}}}^e, {^t}z \right) 
= 
\underbrace{
\frac{K}{2} \left( {^t}\epsilon_v^e \right)^2 + \frac{3}{2} G \left( {^t}\epsilon_q^e \right)^2
}_{ 
= \Psi \left( {^t}\underline{\underline{\boldsymbol{b}}}^e ({^t}\boldsymbol{\underline{\underline{\epsilon}}}^e) \right)}
+ 
\underbrace{\dfrac{H}{2} \left( {^t}z^2 + \frac{\kappa}{\alpha} \left( \exp \left( - \alpha \, {^t}z \right) -1 \right)^2 \right)}_{
= \tilde{\Psi} \left( {^t}z \right)}
,
\end{equation}
where $K = \frac{E}{3 \left( 1- 2 \nu \right) } >0$ and $G = \frac{E}{2 \left( 1 + \nu \right) } >0$ are the constant bulk and the shear moduli, $H$ is the modulus governing the linear hardening when $\kappa = 0$, while, in the case of exponential hardening when $\kappa \neq 0$, this is law is governed by the triad of parameters $H, \kappa,$ and $\alpha$. 
The invariants of the elastic logarithmic strain 
${^t}\epsilon_v^{e}$ and ${^t}\epsilon_q^{e}$ are defined as follows
\begin{equation}
{^t}\epsilon_v^{e} \coloneqq \text{tr} \left( {^t}\boldsymbol{\underline{\underline{\epsilon}}}^{e} \right); 
\qquad
{^t}\boldsymbol{\underline{\underline{e}}}^{e}  \coloneqq {^t}\boldsymbol{\underline{\underline{\epsilon}}}^{e}  - \frac{{^t}\epsilon_v^{e}}{3} \ {^t}\boldsymbol{1};
\qquad
{^t}\epsilon_q^{e} \coloneqq \sqrt{\frac{2}{3}  \, {^t}\boldsymbol{\underline{\underline{e}}}^{e} \boldsymbol{:} {^t}\boldsymbol{\underline{\underline{e}}}^{e}
},
\end{equation}
while ${^t}z$ appearing in Eq.~\eqref{eq:strain_energy} is the strain-like internal variable.

The yield function governing both possible incipient plasticity and the flow rule is that of a Modified Cam-Clay (MCC) (see, e.g.,~\cite{ortiz2004variational}, for its formulation in the context of finite strain mechanics) i.e.,
\begin{equation}
\label{eq:MCC}
\Phi \left( \boldsymbol{\underline{\underline{\xi}}}^{\zeta}, \beta^{\zeta}\right) = 
\left( \dfrac{q^{\xi^{\zeta}}}{m} \right)^2 
+ 
p^{\xi^{\zeta}} \left( p^{\xi^{\zeta}} - p_c \right)
\end{equation}
where the invariants of the stress measure $\left( \bullet \right)$ are defined as follows
\begin{equation}
\label{eq:stress_invariants}
p^{\left( \bullet \right)} 
\coloneqq 
\frac{1}{3} \text{tr} 
\left( \bullet \right); 
\qquad
\boldsymbol{\underline{\underline{s}}}^{\left( \bullet \right)} 
\coloneqq
\left( \bullet \right) - p^{\left( \bullet \right)} \ \boldsymbol{\underline{\underline{1}}};
\qquad
\rho^{\bullet}
\coloneqq 
\sqrt{\boldsymbol{\underline{\underline{s}}}^{\left( \bullet \right)} \boldsymbol{:} \boldsymbol{\underline{\underline{s}}}^{\left( \bullet \right)}}
\qquad
q^{\left( \bullet \right)} 
\coloneqq
\sqrt{\dfrac{3}{2}} \, \rho^{\bullet},
\end{equation}
and the consolidation pressure evolves according to
\begin{equation}
p_c = {^0}p_c+ \beta^{\zeta}.
\end{equation}
The definition of the Eshelby-zeta stress tensor $\boldsymbol{\underline{\underline{\xi}}}^{\zeta}$ is provided in Eq.~\eqref{eq:CD_2} in~\ref{sec:theormodynamics}. 
This section also explains why the yield function depends on this stress measure. Similarly, the stress-like internal variable $\beta^{\zeta}$ is defined in Eq.~\eqref{eq:Colemann_Noll_2} in~\ref{sec:theormodynamics}.

Despite the convention on signs in geotechnics, which considers positive compressive pressures (see, for instance,~\cite{abe2024reconstruction}), this work aligns with the convention on signs more common in solid mechanics, which states positive tensile pressures, as it is clear from the above definition of pressures, Eq.~\eqref{eq:stress_invariants}).
The parameters used for the constitutive equations Eqs.~\eqref{eq:strain_energy} and~\eqref{eq:MCC} and reported in Table~\ref{table:parameters_setups} are those employed by~\cite{bennett2016finite} and calibrated based on the test run for a boulder clay in~\cite{mun2015compression}. 
The only parameter that differs from those in~\cite{bennett2016finite} is the initial consolidation pressure ${^0}p_c$, which was artificially increased to run large-deformation analyses with fewer time steps.

\paragraph*{Results discussion}

\begin{figure}
\centering
\begin{subfigure}[t]{0.475\textwidth}
\includegraphics[width=\textwidth, trim=0cm 0cm 0cm 0cm]{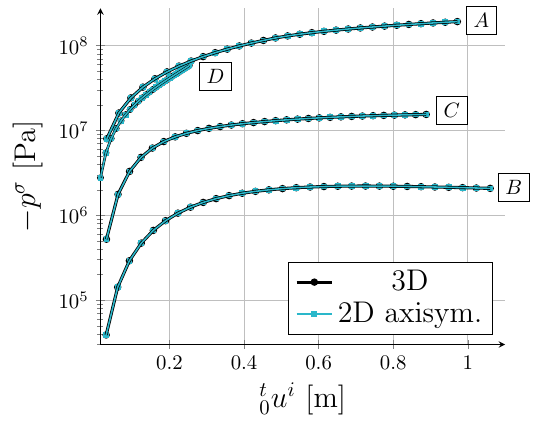}
\caption{Negative Cauchy pressures (logarithmic scale) along the radial displacements for points $A,B,C,$ and $D$.}
\label{subfig:pressure_points}
\end{subfigure}
\hfill
\begin{subfigure}[t]{0.475\textwidth}
\includegraphics[width=\textwidth,  trim=0cm 0cm 0cm 0cm]{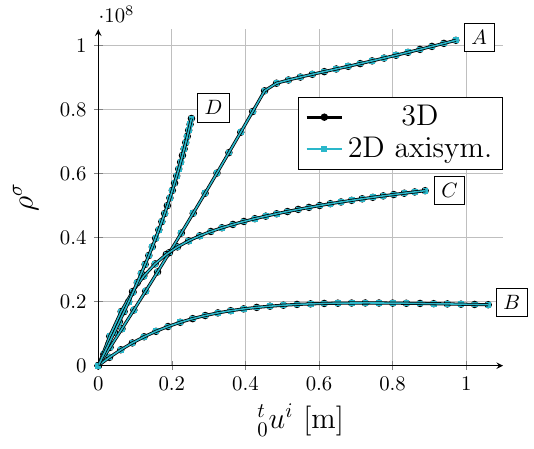}
\caption{$\rho^{\sigma}$ values along the radial displacements for points $A,B,C,$ and $D$.}
\label{subfig:rho_points}
\end{subfigure}
\\
\begin{subfigure}[t]{0.475\textwidth}
\includegraphics[width=\textwidth,  trim=0cm 0cm 0cm 0cm]{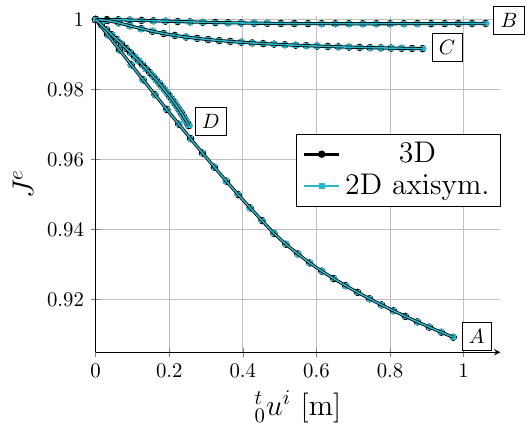}
\caption{Elastic Jacobian values along the radial displacements for points $A,B,C,$ and $D$.}
\label{subfig:J_el_points}
\end{subfigure}
\hfill
\begin{subfigure}[t]{0.475\textwidth}
\includegraphics[width=\textwidth,  trim=0cm 0cm 0cm 0cm]{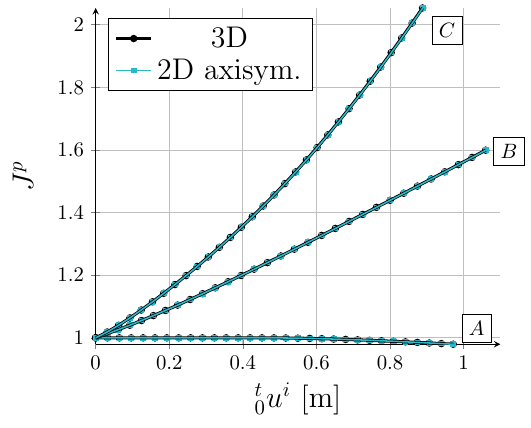}
\caption{Plastic Jacobian values along the radial displacements for points $A,B$ and $C$.}
\label{subfig:J_pl_points}
\end{subfigure}
\caption{
Invariants of Cauchy stresses (top) and Jacobians (bottom) along the radial displacements for points $A,B,C,$ and $D$.
}
\label{fig:points_comparison}
\end{figure}
To compare the results from the three-dimensional and the axisymmetric formulations, several invariant quantities from the Gauss Points closer to points $A, B, C$, and $D$ (refer to Fig.~\ref{fig:setup}) are plotted in Fig.~\ref{fig:points_comparison}. 
The considered quantities include the negative Cauchy pressure in Fig.~\ref{subfig:pressure_points}, the deviatoric Cauchy stress in Fig~\ref{subfig:rho_points}, and the elastic and plastic Jacobians in Figs.~\ref{subfig:J_el_points} and~\ref{subfig:J_pl_points}, which are plotted for various values of radial displacements they experience.
As demonstrated by these pictures, there is no appreciable difference between the values from the three-dimensional analyses and those from the axisymmetric formulation outlined in this paper\footnote{The extremely minor differences between the axisymmetric and three-dimensional analyses arise from the representation of the out-of-plane geometry. In the former case, this is an exact representation, while the latter considers a discrete out-of-plane geometry.}. 
This confirms that the axisymmetric formulation developed in this manuscript delivers three-dimensional-like results at a lower computational cost, as it involves only a two-dimensional discretisation.

\begin{figure}
\centering
\begin{subfigure}[t]{0.475\textwidth}
\includegraphics[width=1.25\textwidth, trim=30cm 8cm 5cm 8cm]{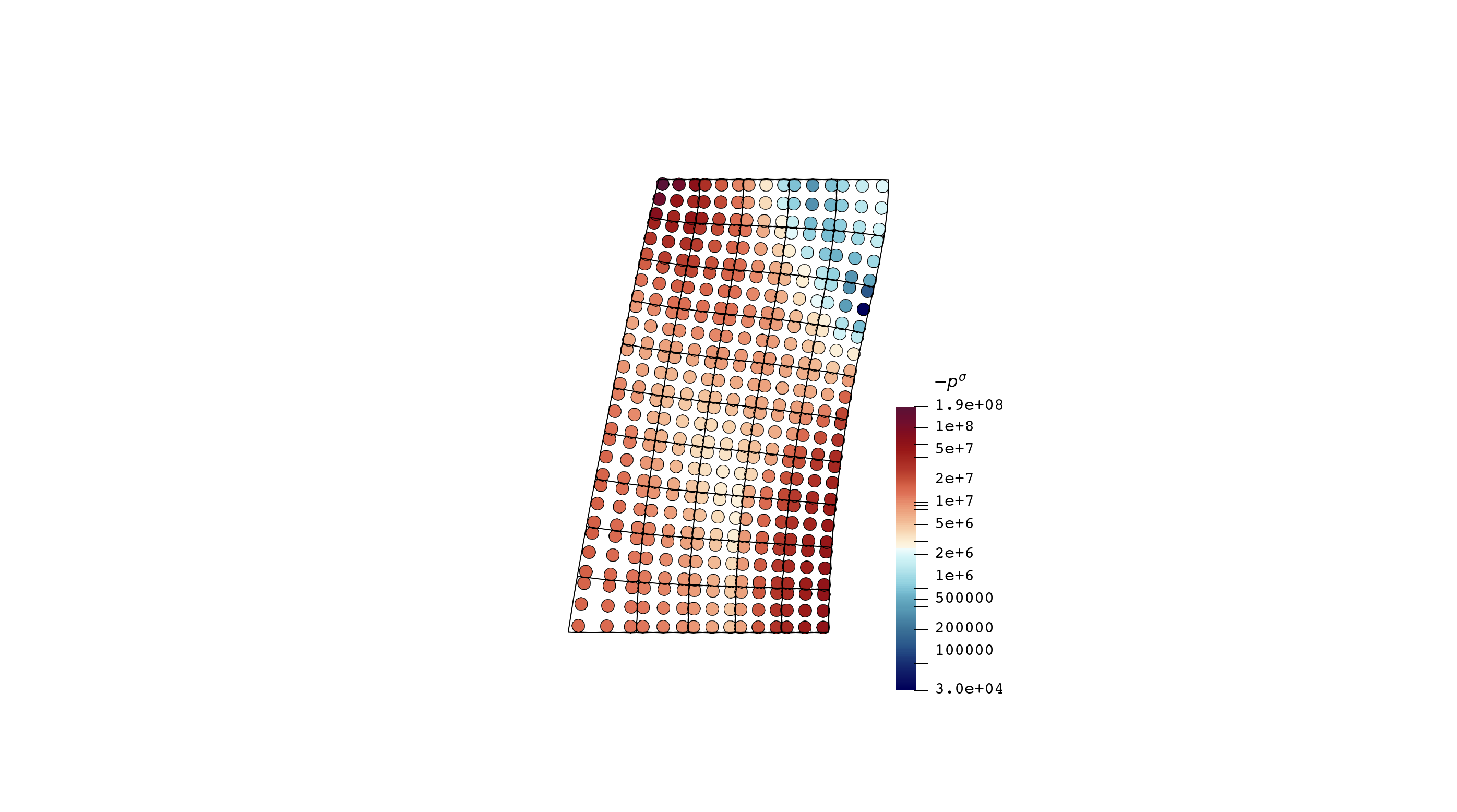}
\caption{Negative Cauchy pressure values (logarithmic scale) from the axisymmetric simulation.}
\label{subfig:pressure_2D}
\end{subfigure}
\hfill
\begin{subfigure}[t]{0.475\textwidth}
\includegraphics[width=1.25\textwidth, trim=30cm 8cm 5cm 8cm]{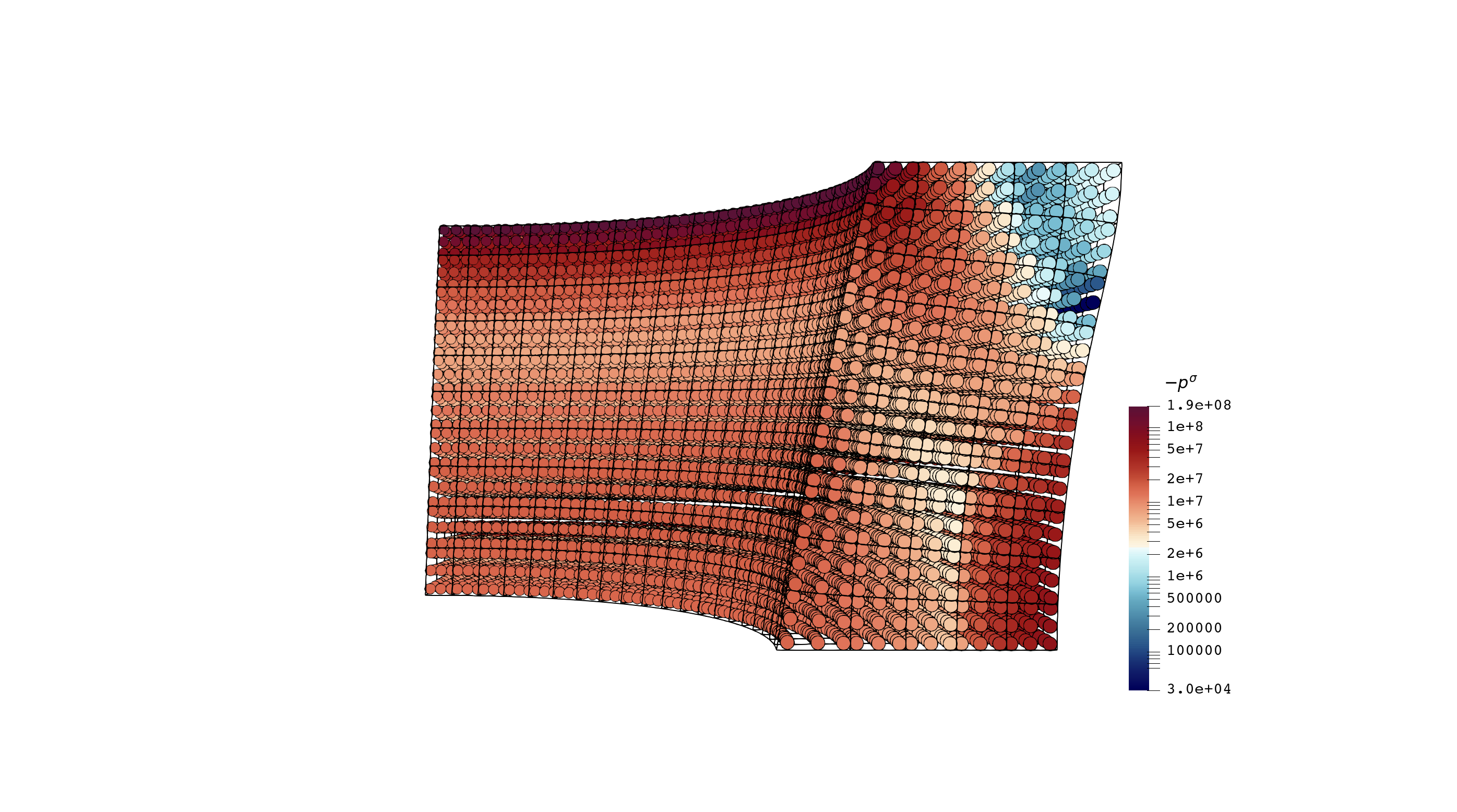}
\caption{Negative Cauchy pressure values (logarithmic scale) from the three-dimensional simulation.}
\label{subfig:pressure_3D}
\end{subfigure}
\\
\begin{subfigure}[t]{0.475\textwidth}
\includegraphics[width=1.25\textwidth, trim=30cm 8cm 5cm 8cm]{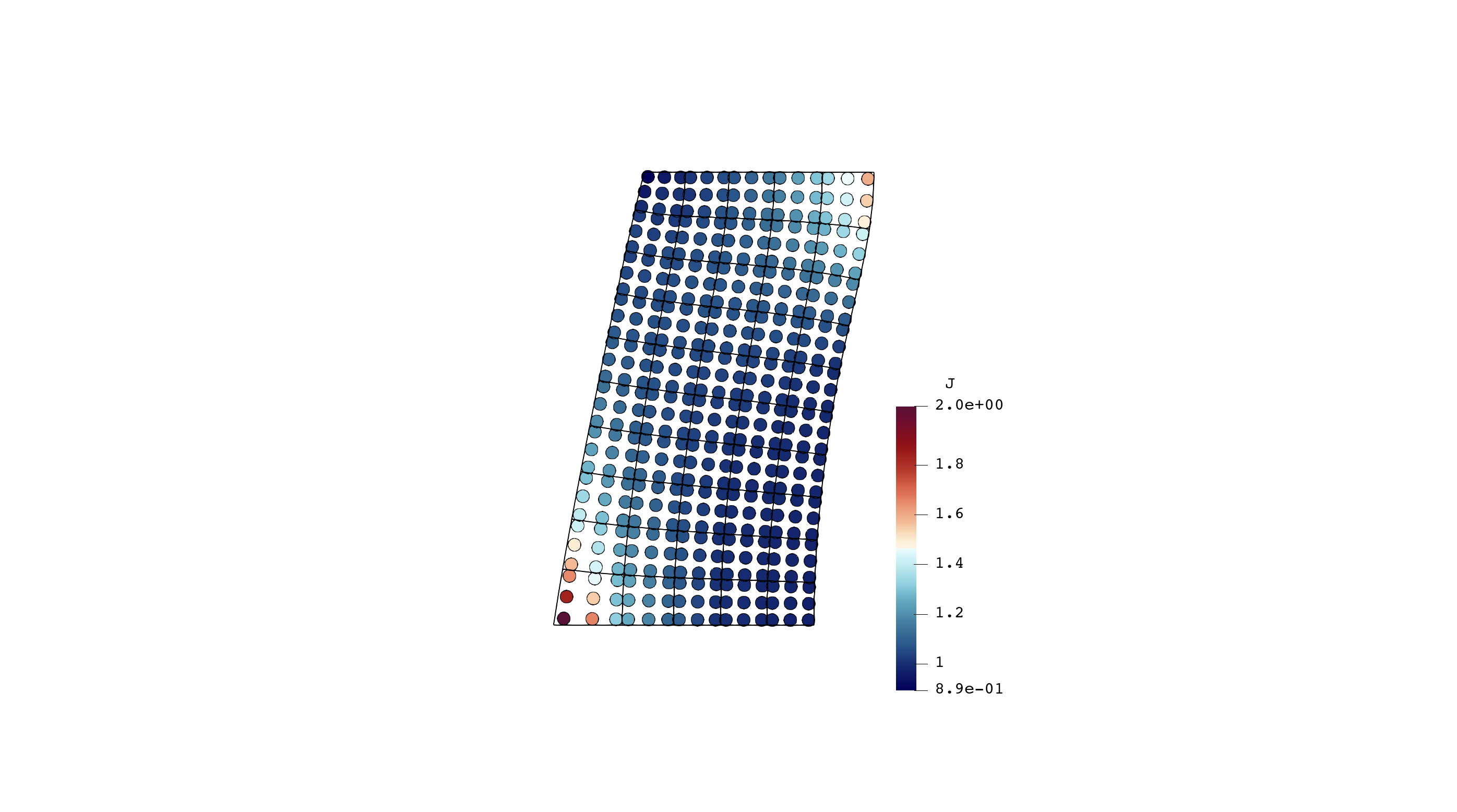}
\caption{Jacobian values from the axisymmetric simulation.}
\label{subfig:Jac_2D}
\end{subfigure}
\hfill
\begin{subfigure}[t]{0.475\textwidth}
\includegraphics[width=1.25\textwidth, trim=30cm 8cm 5cm 8cm]{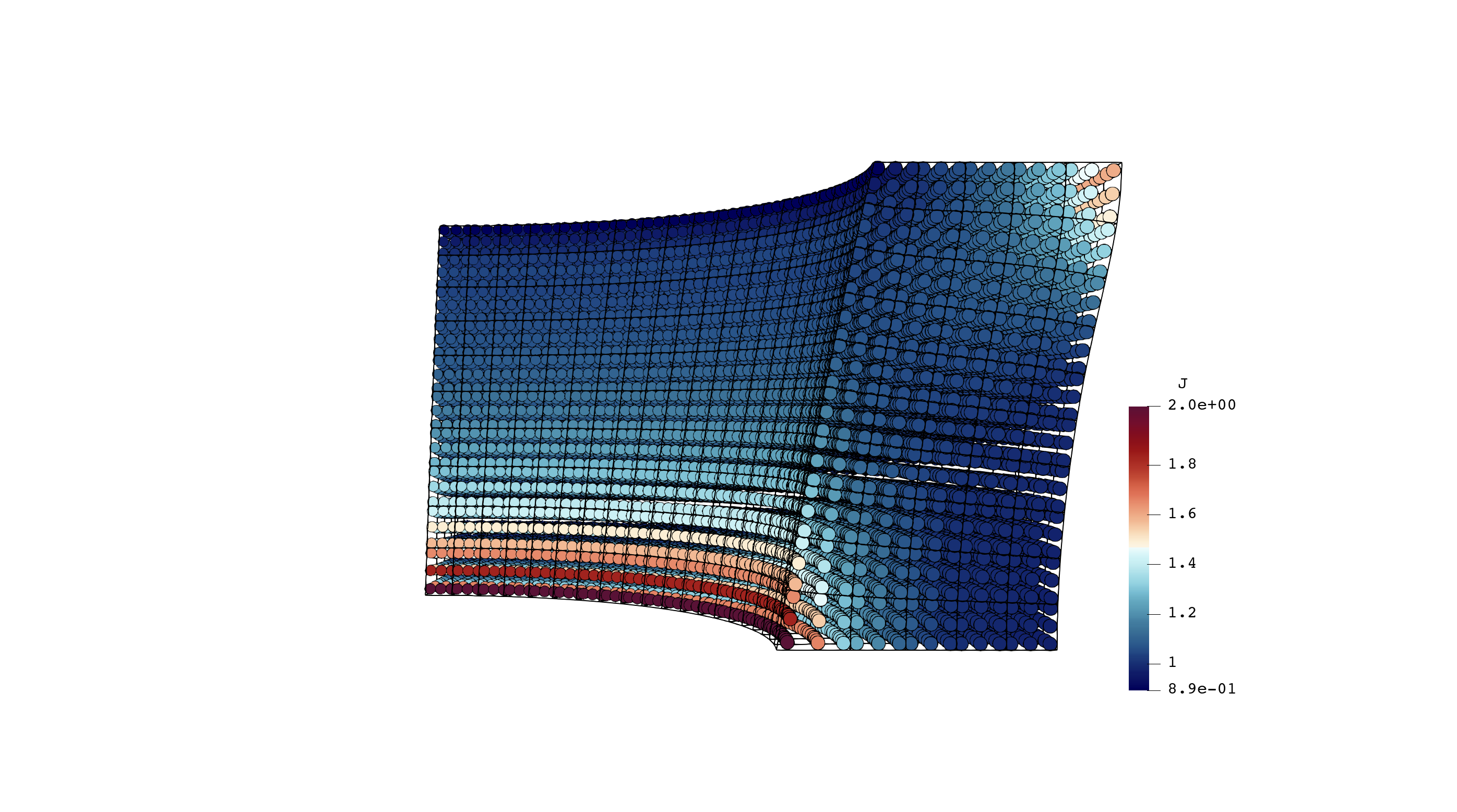}
\caption{Jacobian values from the three-dimensional simulation.}
\label{subfig:Jac_3D}
\end{subfigure}
\caption{
Qualitative comparison of the negative pressure values of Cauchy stress (top) and total Jacobian (bottom) from axisymmetric (left) and three-dimensional (right) simulations.
}
\label{fig:pressure_contours}
\end{figure}

Interestingly, the highest values of stress are registered closer to points $A$ and $D$, both in terms of compressive pressure and deviatoric stress (Fig.~\ref{subfig:pressure_points} and~\ref{subfig:rho_points}, respectively). 
In particular, the higher value of compressive pressure makes the Gauss Point closer to $A$ yield only when radial displacements reach the value of $\approx 0.45$ m, while the Gauss Point closer to $D$ does not yield (this is the reason why its corresponding Gauss Point is excluded from Fig.~\ref{subfig:J_pl_points}. 
The sharp transition to plasticity of point $A$ is clearly visible in Fig.~\ref{subfig:rho_points} and, to a lesser degree, in~\ref{subfig:J_el_points} and~\ref{subfig:J_pl_points}.
Conversely, the Gauss Points closer to points $B$ and $C$ exhibit lower compressive pressure (Fig.~\ref{subfig:pressure_points}), tensile plastic Jacobians since the beginning of the simulations (Fig.~\ref{subfig:J_pl_points}), and lower deviatoric stress (Fig.~\ref{subfig:rho_points}) than those near points $A$ and $D$. 
The qualitative confirmation of these results is also evident in Fig.~\ref{fig:pressure_contours}, which illustrates the negative Cauchy pressure (Figs.~\ref{subfig:pressure_2D} and~\ref{subfig:pressure_3D}) and the total Jacobian (Figs.~\ref{subfig:Jac_2D} and~\ref{subfig:Jac_3D}) at the Gauss Points at the end of the simulations. 

On the other hand, Fig.~\ref{fig:NR_convergence} illustrates the convergence of the NR error at four considered time-steps in the simulations, i.e., $t = 1, \, 10, \, 20$, and $30$. 
Even in this case, very mild differences can be observed between the three-dimensional simulation and the proposed axisymmetric formulation. Specifically, the first time-step requires the most iterations in the simulation (see also Table~\ref{table:NR_iterations}), with the error initially decreasing almost linearly and decreasing more significantly towards the last iterations. 
The tenth time-step requires $4$ iterations, whereas for $t = 20$ and $30$ only $3$ iterations are necessary to achieve the set tolerance. As shown in Table~\ref{table:NR_iterations}, a minimum of $3$ iterations is always required to solve the problem, highlighting the non-linear nature of the equations, and, on average, the three-dimensional and axisymmetric formulations yield very similar results.

\begin{table}
\begin{minipage}[b][7cm][c]{0.5\textwidth}
 \centering
\includegraphics[width=0.75\textwidth, trim=1.5cm 0.5cm 0.5cm 0.cm]{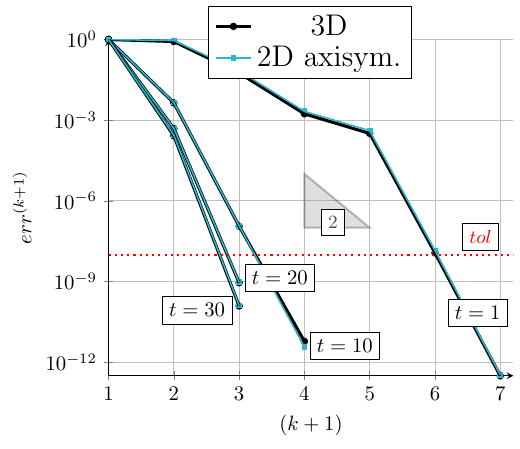}
\captionof{figure}{
Error convergence of the NR iterative loop for time-steps $t=1,10,20,30$.
}
\label{fig:NR_convergence}
\end{minipage} 
\hfill
\begin{minipage}[b][7cm][c]{0.475\linewidth}
\caption{
R\'esum\'e of the NR iterations for the simulations.}
\label{table:NR_iterations}
\centering
\begin{threeparttable} 
\renewcommand{\arraystretch}{1.25} 
{\small
\begin{tabular}{c|c|c} 
\headrow
& 3D & 2D axisym.
\\ \hiderowcolors 
\hline
$\max(k+1)$
& 
$7$
&
$7$
\\
$\min(k+1)$
& 
3
& 
3
\\
$\text{avg} (k+1)$
&
$\approx 3.47$
&
$\approx 3.5$
\\
\hline
\end{tabular}
}
\end{threeparttable}
\end{minipage}
\end{table}

\section{Conclusions}
This work has examined the definition and role of the deformation gradient in curvilinear coordinates, using a cavity‑expansion problem to expose issues that remain hidden in Cartesian orthonormal bases. Particular attention was given to the point‑to‑point mapping underlying the deformation gradient (distinguishing it from the gradient of a vector), the computation of its components from displacement fields, and the correct treatment of its inverse and transpose. The necessity of additional geometric quantities in curvilinear settings, such as the shifter and the Jacobian, and the interplay between their definitions, was also clarified. The discussion of the Jacobian was extended to isotropic elasto‑plasticity, demonstrating how the dependence of the Cauchy stress on the plastic part of the Jacobian affects the consistent linearisation.
To provide a complete kinematic framework, stress measures and the balance of linear momentum were developed in general curvilinear form and then particularised to the axisymmetric case with the practical goal of constructing a finite element implementation. A comparison between fully three‑dimensional and axisymmetric numerical tests confirmed the validity of the proposed axisymmetric formulation and the robustness of the linearisation used in the Newton–Raphson solution scheme.
It is important to emphasise that any incorrect definition of the underlying kinematic quantities compromises all subsequent relations involving their derivatives and the stress measures. Such errors propagate to the governing equations and can lead, for example, to inconsistencies between updated Lagrangian and total Lagrangian formulations. The results presented here underscore the need for careful and precise kinematic definitions when working in curvilinear coordinate systems, particularly in finite‑strain, axisymmetric analyses.

\section*{Acknowledgements}
In the preparation of this work, the first author was supported as a Postdoctoral Associate by the Department of Engineering at Durham University. 
Support from UK Research and Innovation (UKRI) under Grant EP/W000970/1 is gratefully acknowledged by the second author.
All data created during this research are openly available at \href{https://collections.durham.ac.uk/}{collections.durham.ac.uk/} (specific DOI to be confirmed if/when the paper is accepted).
For the purpose of open access, the author has applied a Creative Commons Attribution (CC BY) licence to any Author Accepted Manuscript version arising.

\bibliographystyle{elsarticle-num}
\bibliography{bibliography.bib}

\begin{thebibliography}{10}
\expandafter\ifx\csname url\endcsname\relax
  \def\url#1{\texttt{#1}}\fi
\expandafter\ifx\csname urlprefix\endcsname\relax\def\urlprefix{URL }\fi
\expandafter\ifx\csname href\endcsname\relax
  \def\href#1#2{#2} \def\path#1{#1}\fi

\bibitem{said2009axisymmetric}
I.~Said, V.~De~Gennaro, R.~Frank, Axisymmetric finite element analysis of pile
  loading tests, Computers and Geotechnics 36~(1-2) (2009) 6--19.

\bibitem{zhao2015computational}
K.~Zhao, M.~Bonini, D.~Debernardi, M.~Janutolo, G.~Barla, G.~Chen,
  Computational modelling of the mechanised excavation of deep tunnels in weak
  rock, Computers and Geotechnics 66 (2015) 158--171.

\bibitem{bird2024implicit}
R.~E. Bird, G.~Pretti, W.~M. Coombs, C.~E. Augarde, Y.~U. Sharif, M.~J. Brown,
  G.~Carter, C.~Macdonald, K.~Johnson, An implicit material point-to-rigid body
  contact approach for large deformation soil--structure interaction, Computers
  and Geotechnics 174 (2024) 106646.

\bibitem{mohapatra2025laboratory}
D.~Mohapatra, S.~Mohammadi, M.~Saresma, J.~J. Virtasalo, W.~T. So{\l}owski,
  Laboratory-scale free fall cone penetrometer test on marine clay: A numerical
  investigation using the generalized interpolation material point method,
  International Journal for Numerical and Analytical Methods in Geomechanics
  49~(4) (2025) 1299--1318.

\bibitem{williamson2008comparison}
N.~Williamson, M.~Behnia, S.~Armfield, Comparison of a 2d axisymmetric cfd
  model of a natural draft wet cooling tower and a 1d model, International
  journal of heat and mass transfer 51~(9-10) (2008) 2227--2236.

\bibitem{martinez2002simulation}
M.~Martinez, I.~Alfaro, M.~Doblare, Simulation of axisymmetric discharging in
  metallic silos. analysis of the induced pressure distribution and comparison
  with different standards, Engineering Structures 24~(12) (2002) 1561--1574.

\bibitem{yoshida2010axisymmetric}
S.~Yoshida, K.~Sekine, T.~Mitsuta, Axisymmetric finite element analysis for
  sloshing response of floating roofs in cylindrical storage tanks, Journal of
  Environment and Engineering 5~(1) (2010) 27--38.

\bibitem{lee1971analyses}
C.~H. Lee, S.~Kobayashi, Analyses of axisymmetric upsetting and plane-strain
  side-pressing of solid cylinders by the finite element method, Journal of
  Engineering for Industry 93~(2) (1971) 445--454.

\bibitem{oyinbo2015numerical}
S.~Oyinbo, O.~Ikumapayi, J.~Ajiboye, S.~Afolalu, Numerical simulation of
  axisymmetric and asymmetric extrusion process using finite element method,
  International Journal of Scientific \& Engineering Research 6~(6) (2015)
  1246--1259.

\bibitem{mcallen2007numerical}
P.~McAllen, P.~Phelan, Numerical analysis of axisymmetric wire drawing by means
  of a coupled damage model, Journal of Materials Processing Technology
  183~(2-3) (2007) 210--218.

\bibitem{kumar2025nonlinear}
A.~Kumar, A.~Yavari, Nonlinear mechanics of arterial growth, arXiv preprint
  arXiv:2511.12819 (2025).

\bibitem{de2023estimation}
A.~de~la Hoz, E.~Martinez-Enriquez, S.~Marcos, Estimation of crystalline lens
  material properties from patient accommodation data and finite element
  models, Investigative ophthalmology \& visual science 64~(11) (2023) 31--31.

\bibitem{seth1968transient}
M.~Seth, K.~Gray, Transient stresses and displacement around a wellbore due to
  fluid flow in transversely isotropic, porous media: Ii. finite reservoirs,
  Society of Petroleum Engineers Journal 8~(01) (1968) 79--86.

\bibitem{kenyon1979mathematical}
D.~E. Kenyon, A mathematical model of water flux through aortic tissue,
  Bulletin of mathematical biology 41~(1) (1979) 79--90.

\bibitem{stiles2006two}
D.~K. Stiles, D.~U. Gwost, Two-dimensional axisymmetric model for the
  sensitivity analysis of a chronic drug infusion into the brain, in: 2006
  International Conference of the IEEE Engineering in Medicine and Biology
  Society, IEEE, 2006, pp. 1521--1524.

\bibitem{davoli2015critical}
E.~Davoli, G.~A. Francfort, A critical revisiting of finite elasto-plasticity,
  SIAM Journal on Mathematical Analysis 47~(1) (2015) 526--565.

\bibitem{yavari2023direct}
A.~Yavari, F.~Sozio, On the direct and reverse multiplicative decompositions of
  deformation gradient in nonlinear anisotropic anelasticity, Journal of the
  Mechanics and Physics of Solids 170 (2023) 105101.

\bibitem{abe2024reconstruction}
Y.~Abe, S.~Yamada, K.~Hoshi, T.~Kyoya, Reconstruction of {C}am-{c}lay model
  based on multiplicative decomposition of the deformation gradient, Computers
  and Geotechnics 166 (2024) 105958.

\bibitem{dvorkin2006nonlinear}
E.~N. Dvorkin, M.~B. Goldschmit, Nonlinear continua, Springer Science \&
  Business Media, 2006.

\bibitem{ogden1997non}
R.~W. Ogden, Non-linear elastic deformations, Courier Corporation, 1997.

\bibitem{einstein2003meaning}
A.~Einstein, The meaning of relativity, Routledge, 2003.

\bibitem{hughes1983mathematical}
T.~J. Hughes, J.~Marsden, Mathematical foundations of elasticity, Citeseer,
  1983.

\bibitem{kanso2007geometric}
E.~Kanso, M.~Arroyo, Y.~Tong, A.~Yavari, J.~G. Marsden, M.~Desbrun, On the
  geometric character of stress in continuum mechanics, Zeitschrift f{\"u}r
  angewandte Mathematik und Physik 58~(5) (2007) 843--856.

\bibitem{danielson1997three}
K.~T. Danielson, A.~K. Noor, Three-dimensional finite element analysis in
  cylindrical coordinates for nonlinear solid mechanics problems, Finite
  elements in analysis and design 27~(3) (1997) 225--249.

\bibitem{celigoj1998assumed}
C.~Celigoj, An assumed enhanced displacement gradient ring-element for finite
  deformation axisymmetric and torsional problems, International journal for
  numerical methods in engineering 43~(8) (1998) 1369--1382.

\bibitem{rauchs2016direct}
G.~Rauchs, Direct-differentiation-based sensitivity analysis of an axisymmetric
  finite element formulation including torsion, Finite Elements in Analysis and
  Design 109 (2016) 65--72.

\bibitem{steigmann2017finite}
D.~J. Steigmann, Finite elasticity theory, Oxford University Press, 2017.

\bibitem{kohler2001domain}
O.~K{\"o}hler, G.~Kuhn, The domain-boundary element method ({DBEM}) for
  hyperelastic and elastoplastic finite deformation: axisymmetric and 2d/3d
  problems, Archive of Applied Mechanics 71~(6) (2001) 436--452.

\bibitem{truesdell1960classical}
C.~Truesdell, R.~Toupin, The classical field theories, in: Principles of
  classical mechanics and field theory/Prinzipien der Klassischen Mechanik und
  Feldtheorie, Springer, 1960, pp. 226--858.

\bibitem{leborgne2023objectivity}
G.~Leborgne, Objectivity in continuum mechanics, an introduction. motions,
  eulerian and lagrangian variables and functions, deformation gradient, lie
  derivatives, velocity-addition formula, coriolis, arXiv preprint
  arXiv:2301.01056 (2023).

\bibitem{steigmann2025principles}
D.~J. Steigmann, M.~Shirani, Principles of Continuum Mechanics, World
  Scientific, 2025.

\bibitem{frankel2004geometry}
T.~Frankel, The geometry of physics: an introduction, Cambridge university
  press, 2004.

\bibitem{brannon2018rotation}
R.~M. Brannon, Rotation, Reflection, and Frame Changes: Orthogonal tensors in
  computational engineering mechanics, IOP Publishing, 2018.

\bibitem{sulsky1994particle}
D.~Sulsky, Z.~Chen, H.~Schreyer, A particle method for history-dependent
  materials, Computer {M}ethods in {A}pplied {M}echanics and {E}ngineering
  118~(1-2) (1994) 179--196.

\bibitem{charlton2017igimp}
T.~Charlton, W.~Coombs, C.~Augarde, igimp: An implicit generalised
  interpolation material point method for large deformations, Computers \&
  Structures 190 (2017) 108--125.

\bibitem{bilby1957continuous}
B.~A. Bilby, L.~Lardner, A.~Stroh, Continuous distributions of dislocations and
  the theory of plasticity, in: Actes du IXe congr\'es international de
  m\'ecanique appliqu\'ee, Vol.~8, Bruxelles, 1957, pp. 35--44.

\bibitem{kroner1959allgemeine}
E.~Kr{\"o}ner, Allgemeine kontinuumstheorie der versetzungen und
  eigenspannungen, Archive for Rational Mechanics and Analysis 4~(1) (1959)
  273--334.

\bibitem{lee1969elastic}
E.~H. Lee, Elastic-plastic deformation at finite strains, Journal of Applied
  Mechanics 36~(1) (1969) 1--6.

\bibitem{goodbrake2021mathematical}
C.~Goodbrake, A.~Goriely, A.~Yavari, The mathematical foundations of
  anelasticity: existence of smooth global intermediate configurations,
  Proceedings of the Royal Society A 477~(2245) (2021) 20200462.

\bibitem{casey1980remark}
J.~Casey, P.~M. Naghdi, A remark on the use of the decomposition f= f sub e f
  sub p in plasticity., Tech. rep., California University of Berkeley (1980).

\bibitem{desouza2011computational}
E.~A. de~Souza~Neto, D.~Peric, D.~R. Owen, Computational methods for
  plasticity: theory and applications, John Wiley \& Sons, 2011.

\bibitem{simo1988frameworkI}
J.~C. Simo, A framework for finite strain elastoplasticity based on maximum
  plastic dissipation and the multiplicative decomposition: Part i. {C}ontinuum
  formulation, Computer methods in applied mechanics and engineering 66~(2)
  (1988) 199--219.

\bibitem{simo1988frameworkII}
J.~C. Simo, A framework for finite strain elastoplasticity based on maximum
  plastic dissipation and the multiplicative decomposition. part ii:
  Computational aspects, Computer methods in applied mechanics and engineering
  68~(1) (1988) 1--31.

\bibitem{sadik2025generalized}
S.~Sadik, A.~Yavari, A generalized coleman--noll procedure and the balance laws
  of hyper-anelasticity, Proceedings of the Royal Society A 481~(2320) (2025)
  20250076.

\bibitem{oden2006finite}
J.~T. Oden, Finite elements of nonlinear continua, Courier Corporation, 2006.

\bibitem{dhas2020geometric}
B.~Dhas, D.~Roy, A geometric approach to hu-washizu variational principle in
  nonlinear elasticity, arXiv preprint arXiv:2009.00275 (2020).

\bibitem{bennett2016finite}
K.~Bennett, R.~Regueiro, R.~Borja, Finite strain elastoplasticity considering
  the eshelby stress for materials undergoing plastic volume change,
  International Journal of Plasticity 77 (2016) 214--245.

\bibitem{epstein1990energy}
M.~Epstein, G.~Maugin, The energy-momentum tensor and material uniformity in
  finite elasticity, Acta Mechanica 83~(3) (1990) 127--133.

\bibitem{hencky1933elastic}
H.~Hencky, The elastic behavior of vulcanized rubber, Journal of Applied
  Mechanics (1933).

\bibitem{simo1998numerical}
J.~C. Simo, Numerical analysis and simulation of plasticity, Handbook of
  numerical analysis 6 (1998) 183--499.

\bibitem{dafalias1984plastic}
Y.~F. Dafalias, The plastic spin concept and a simple illustration of its role
  in finite plastic transformations, Mechanics of Materials 3~(3) (1984)
  223--233.

\bibitem{dafalias1985plastic}
Y.~F. Dafalias, The plastic spin, Journal of Applied Mechanics 52~(4) (1985)
  865--871.

\bibitem{bennett2019anisotropic}
K.~Bennett, R.~Regueiro, D.~J. Luscher, Anisotropic finite
  hyper-elastoplasticity of geomaterials with drucker--prager/cap type
  constitutive model formulation, International Journal of Plasticity 123
  (2019) 224--250.

\bibitem{lubliner1972thermodynamic}
J.~Lubliner, On the thermodynamic foundations of non-linear solid mechanics,
  International Journal of Non-Linear Mechanics 7~(3) (1972) 237--254.

\bibitem{hill1948variational}
R.~Hill, A variational principle of maximum plastic work in classical
  plasticity, The Quarterly Journal of Mechanics and Applied Mathematics 1~(1)
  (1948) 18--28.

\bibitem{mandel1964contribution}
J.~Mandel, Contribution th{\'e}orique {\`a} l’{\'e}tude de
  l’{\'e}crouissage et des lois de l’{\'e}coulement plastique, in: Applied
  Mechanics: Proceedings of the Eleventh International Congress of Applied
  Mechanics Munich (Germany) 1964, Springer, 1964, pp. 502--509.

\bibitem{collins1997application}
I.~Collins, G.~Houlsby, Application of thermomechanical principles to the
  modelling of geotechnical materials, Proceedings of the Royal Society of
  London. Series A: Mathematical, Physical and Engineering Sciences 453~(1964)
  (1997) 1975--2001.

\bibitem{collins2002associated}
I.~Collins, Associated and non-associated aspects of the constitutive laws for
  coupled elastic/plastic materials, The International Journal Geomechanics
  2~(2) (2002) 259--267.

\bibitem{oliynyk2020finite}
K.~Oliynyk, C.~Tamagnini, Finite deformation hyperplasticity theory for
  crushable, cemented granular materials, Open Geomechanics 2 (2020) 1--33.

\bibitem{pretti2024preserving}
G.~Pretti, W.~M. Coombs, C.~E. Augarde, M.~M. Puigvert, J.~A.~R. Guti{\'e}rrez,
  Preserving non-negative porosity values in a bi-phase elasto-plastic material
  under terzaghi’s effective stress principle, Mechanics of Materials 192
  (2024) 104958.

\bibitem{simo2005stress}
J.~Simo, J.~Marsden, Stress tensors, riemannian metrics and the alternative
  descriptions in elasticity, in: Trends and Applications of Pure Mathematics
  to Mechanics: Invited and Contributed Papers presented at a Symposium at
  Ecole Polytechnique, Palaiseau, France November 28--December 2, 1983,
  Springer, 2005, pp. 369--383.

\bibitem{menzel2007configurational}
A.~Menzel, P.~Steinmann, On configurational forces in multiplicative
  elastoplasticity, International Journal of Solids and Structures 44~(13)
  (2007) 4442--4471.

\bibitem{coombs2020lagrangian}
W.~M. Coombs, C.~E. Augarde, A.~J. Brennan, M.~J. Brown, T.~J. Charlton, J.~A.
  Knappett, Y.~G. Motlagh, L.~Wang, On lagrangian mechanics and the implicit
  material point method for large deformation elasto-plasticity, Computer
  Methods in Applied Mechanics and Engineering 358 (2020) 112622.

\bibitem{simo1998computational}
J.~C. Simo, T.~J. Hughes, Computational inelasticity, Springer, 1998.

\bibitem{deparis2004numerical}
S.~Deparis, Numerical analysis of axisymmetric flows and methods for
  fluid-structure interaction arising in blood flow simulation, Tech. rep.,
  EPFL (2004).

\bibitem{bernardi1999spectral}
C.~Bernardi, M.~Dauge, Y.~Maday, Spectral methods for axisymmetric domains, (No
  Title) (1999).

\bibitem{haughton1979bifurcation}
D.~Haughton, R.~Ogden, Bifurcation of inflated circular cylinders of elastic
  material under axial loading—ii. exact theory for thick-walled tubes,
  Journal of the Mechanics and Physics of Solids 27~(5-6) (1979) 489--512.

\bibitem{steigmann2023course}
D.~J. Steigmann, A course on plasticity theory, Vol.~7, Oxford University
  Press, 2023.

\bibitem{boffi2013mixed}
D.~Boffi, F.~Brezzi, M.~Fortin, et~al., Mixed finite element methods and
  applications, Vol.~44, Springer, 2013.

\bibitem{ortiz2004variational}
M.~Ortiz, A.~Pandolfi, A variational cam-clay theory of plasticity, Computer
  Methods in Applied Mechanics and Engineering 193~(27-29) (2004) 2645--2666.

\bibitem{mun2015compression}
W.~Mun, J.~S. McCartney, Compression mechanisms of unsaturated clay under high
  stresses, Canadian Geotechnical Journal 52~(12) (2015) 2099--2112.

\bibitem{coleman1967thermodynamics}
B.~D. Coleman, M.~E. Gurtin, Thermodynamics with internal state variables, The
  journal of chemical physics 47~(2) (1967) 597--613.

\bibitem{coleman1974thermodynamics}
B.~D. Coleman, W.~Noll, The thermodynamics of elastic materials with heat
  conduction and viscosity, in: The foundations of mechanics and
  thermodynamics: Selected papers, Springer, 1974, pp. 145--156.

\bibitem{coombs2011algorithmic}
W.~M. Coombs, R.~S. Crouch, Algorithmic issues for three-invariant hyperplastic
  critical state models, Computer methods in applied mechanics and engineering
  200~(25-28) (2011) 2297--2318.

\end{thebibliography}

\appendix
\section{Thermodynamic framework for the constitutive equations}
\label{sec:theormodynamics}
As stated in Section~\ref{subsec:constitutive_equations}, this work considers an elastically and plastically isotropic material, governed by a stored energy function (hyperelasticity) and by the principle of maximum dissipation (associated evolution laws).

Before delving into the thermodynamics, let us introduce the intermediate (stress-free) infinitesimal volume $\text{d}\bar{V}$, which is related to the initial and the current volumes by the following relationships 
\begin{equation}
\text{d}^{t}V = \ {^t_0}J \ \text{d}^{0}V = J^{e} \ \text{d}\bar{V}.
\end{equation}
Similarly, the intermediate density $\bar{\rho}$ can be expressed via the principle of mass conservation to its initial and current counterparts via 
\begin{equation}
\label{eq:densities}
{^t}\rho \ {^t_0}J = {^0}\rho = \bar{\rho} \ J^{p}.
\end{equation}

As a starting point for this appendix, let us consider the (strong) Clausius-Duhem inequality (see, e.g.,~\cite{desouza2011computational} and~\cite{bennett2016finite}) per unit initial volume for isothermal transformations
\begin{equation}
\label{eq:CD_0}
{^t_0} \underline{\underline{\boldsymbol{P}}} \boldsymbol{:} \frac{\partial {^t_0} \underline{\underline{\boldsymbol{X}}}}{\partial t} 
- 
{^0}\rho \, \frac{\partial \psi}{\partial t} \geq 0,
\end{equation}
where $\psi$ is the stored energy function per unit mass and the time derivative of the deformation gradient is given by 
\begin{align}
\nonumber
\dfrac{\partial {^t_0} \underline{\underline{\boldsymbol{X}}}}{\partial t} 
& = 
\dfrac{\partial }{\partial t} \left( {^t_0} X^{a}_{\ A}  {^t}\boldsymbol{\underline{g}}_a \right) {^0}\boldsymbol{\underline{g}}_A
\\
\nonumber
& =
\dfrac{\partial }{\partial t} 
\left( \dfrac{\partial {^t}x^a}{\partial {^0}
x^A} \right) 
{^t}\boldsymbol{\underline{g}}_a \, {^0}\boldsymbol{\underline{g}}_A
+
{^t_0} X^{a}_{\ A} \frac{\partial {^t}\boldsymbol{\underline{g}}_a}{\partial {^t}x^b} \frac{\partial {^t}x^b}{\partial t} {^0}\boldsymbol{\underline{g}}_A
\\
& =
\left( 
\dfrac{\partial {^t}V^a}{\partial {^0}
x^A} 
+ {^t_0}X^{c}_{\ A} \, ^t\gamma^{a}_{bc} \,  {^t}V^b 
\right) {^t}\boldsymbol{\underline{g}}_a \, {^0}\boldsymbol{\underline{g}}_A 
= {^t}V^a|_{A} {^t}\boldsymbol{\underline{g}}_a \, {^0}\boldsymbol{\underline{g}}_A,
\end{align}
In the above chain of equations, the material velocity has been introduced and it is defined by ${^t}V^a  \coloneqq \frac{\partial ^tx^a}{\partial t}$. Its spatial counterpart is defined by the relationship ${^t}\underline{\boldsymbol{V}} \left( \phi^{-1} \left( p \right) \right) = {^t}\underline{\boldsymbol{v}}$. 
The chain rule has also been used to derive the current curvilinear basis ${^t}\boldsymbol{\underline{g}}_a$, and the current Christoffel symbol $^t\gamma^{a}_{bc}$ has been evoked following an equation similar to Eq.~\eqref{eq:Christoffel_symbol_original}. 
Owing to their definitions, the components of the covariant derivatives of the velocities follow a relationship similar to Eq.~\eqref{eq:gradient_transformation}, i.e.,
\begin{equation}
{^t}V^a|_{A} \left( \phi^{-1} \left( p \right) \right) = {^t}v^a|_{b} \ {^t_0} X^{b}_{\ A} \left( \phi^{-1} \left( p \right) \right).    
\end{equation}

Let us defined the spatial velocity gradient as follows ${^t}\underline{\underline{\boldsymbol{l}}} \coloneqq {^t}v^a|_{b} {^t}\boldsymbol{\underline{g}}_a {^t}\boldsymbol{\underline{g}}^b$. 
This quantity is extremely advantageous in the context of finite elasto-plasticity as it can be additively decomposed into an elastic ${^t}\underline{\underline{\boldsymbol{l}}}^e$ and a plastic ${^t}\underline{\underline{\boldsymbol{l}}}^p$ parts as follows
\begin{align}
\label{eq:velocity_gradient_additive_decomposition}
\nonumber
{^t}\underline{\underline{\boldsymbol{l}}} = 
\dfrac{\partial {^t_0} \underline{\underline{\boldsymbol{X}}}}{\partial t} 
\cdot 
{^t_0}\underline{\underline{\boldsymbol{X}}}^{-1}
& =
\dfrac{\partial  }{\partial t} 
\left( {^t_0}\underline{\underline{\boldsymbol{X}}}^e \cdot {^0_0}\underline{\underline{\boldsymbol{X}}}^p \right)
\cdot
\left( 
({^t_0}\underline{\underline{\boldsymbol{X}}}^p)^{-1} \cdot ({^t_0}\underline{\underline{\boldsymbol{X}}}^e)^{-1}
\right)
\\
\nonumber
& =
\underbrace{
\dfrac{\partial {^t_0}\underline{\underline{\boldsymbol{X}}}^e}{\partial t} \cdot ({^t_0}\underline{\underline{\boldsymbol{X}}}^e)^{-1}
}_{
\coloneqq {^t}\underline{\underline{\boldsymbol{l}}}^e }
+
{^t_0}\underline{\underline{\boldsymbol{X}}}^e \cdot \underbrace{ \dfrac{\partial  {^0_0}\underline{\underline{\boldsymbol{X}}}^p}{\partial t} 
\cdot
({^t_0}\underline{\underline{\boldsymbol{X}}}^p)^{-1} }_{\coloneqq {^0}\underline{\underline{\boldsymbol{L}}}^p 
} \cdot ({^t_0}\underline{\underline{\boldsymbol{X}}}^e)^{-1}
\\
&
=
{^t}\underline{\underline{\boldsymbol{l}}}^e 
+
\underbrace{ {^t_0}\underline{\underline{\boldsymbol{X}}}^e \cdot
{^0}\underline{\underline{\boldsymbol{L}}}^p 
\cdot ({^t_0}\underline{\underline{\boldsymbol{X}}}^e)^{-1}
}_{\coloneqq {^t}\underline{\underline{\boldsymbol{l}}}^p },
\end{align}
where the multiplicative decomposition~\eqref{eq:elasto-plastic_def_gradient} has been recalled and the plastic velocity gradient ${^0}\underline{\underline{\boldsymbol{L}}}^p$ in the stress-free configuration is introduced. 
Denoting by ${^0}\underline{\underline{\boldsymbol{D}}}^p$ and ${^0}\underline{\underline{\boldsymbol{W}}}^p$ the symmetric and a skew-symmetric parts of 
${^0}\underline{\underline{\boldsymbol{L}}}^p$, i.e.,
\begin{equation}
{^0}\underline{\underline{\boldsymbol{D}}}^p \coloneqq \dfrac{{^0}\underline{\underline{\boldsymbol{L}}}^p + ({^0}\underline{\underline{\boldsymbol{L}}}^p)^T}{2}, 
\qquad
{^0}\underline{\underline{\boldsymbol{W}}}^p \coloneqq \dfrac{{^0}\underline{\underline{\boldsymbol{L}}}^p -({^0}\underline{\underline{\boldsymbol{L}}}^p)^T}{2},
\end{equation}
the \emph{plastic spin tensor} ${^0}\underline{\underline{\boldsymbol{W}}}^p$ (discussed in Section~\ref{subsec:constitutive_equations}) is defined, and setting ${^0}\underline{\underline{\boldsymbol{W}}}^p = {^0}\underline{\underline{\boldsymbol{0}}}$ implies the assumption of isotropic plasticity discussed in the same section.
Similar symmetric and skew-symmetric definitions of ${^t}\underline{\underline{\boldsymbol{l}}}^p$ are also considered below, and these are denoted by ${^t}\underline{\underline{\boldsymbol{d}}}^p$ and ${^t}\underline{\underline{\boldsymbol{w}}}^p$, respectively.

Employing the relationship between densities in different configurations given in Eq.~\eqref{eq:densities}, the definition of the spatial velocity gradient ${^t}\underline{\underline{\boldsymbol{l}}}$, and the relationship between stress measures Eqs.~\eqref{eq:first_PK_stress_def}-\eqref{eq:semi_kirchhoff_stress}, inequality~\eqref{eq:CD_0} can be expressed as
\begin{equation}
\label{eq:CD_1}
{^t_0} \underline{\underline{\boldsymbol{P}}} \boldsymbol{:} \frac{\partial {^t_0} \underline{\underline{\boldsymbol{X}}}}{\partial t} 
- 
{^0}\rho \, \frac{\partial \psi}{\partial t} 
= {^t_0}J \left(   {^t}\underline{\underline{\boldsymbol{\sigma}}} \boldsymbol{:} {^t} \underline{\underline{\boldsymbol{l}}} 
- {^t}\rho \, \frac{\partial \psi}{\partial t} \right)
= 
J^p \left(  \underline{\underline{\boldsymbol{\zeta}}} \boldsymbol{:} {^t} \underline{\underline{\boldsymbol{l}}} 
- \bar{\rho} \, \frac{\partial \psi}{\partial t} \right) \geq 0,
\end{equation}
since the following equations hold
\begin{equation}
\bar{\rho} \, \frac{\partial \psi}{\partial t}  =  \frac{\partial \left( \bar{\rho} \, \psi \right)}{\partial t} -  \frac{\partial \bar{\rho}}{\partial t} \psi
=
\frac{\partial \left( \bar{\rho} \, \psi \right)}{\partial t} +  \frac{\partial J^p}{\partial t}  \frac{\bar{\rho}}{J^p} \, \psi.
\end{equation}
From the properties of the determinant of matrices, it follows that 
\begin{align}
\nonumber
\frac{\partial J^p}{\partial t} 
& = J^p \,  ({^0_0}\underline{\underline{\boldsymbol{X}}}^{p} )^{-T} \boldsymbol{:} \dfrac{\partial {^0_0}\underline{\underline{\boldsymbol{X}}}^{p}}{\partial t} 
\\
\nonumber
& =
J^p \, \left( ({^0_0}X^{p} )^{-T} \right)^{A}_{\ B} \ ({^0_0}X^{p} )^{A}_{\ B}
\\
& =
J^p \, \left( ({^0_0}X^{p} )^{-1} \right)_{C}^{\ D} \ ({^0_0}X^{p} )^{C}_{\ D}
=
J^p \, \text{tr} \left( {^0}\underline{\underline{\boldsymbol{L}}}^p\right),
\end{align}
and, following their definitions in Eq.~\eqref{eq:velocity_gradient_additive_decomposition}, we get
\begin{equation}
\text{tr} \left( {^t}\underline{\underline{\boldsymbol{l}}}^p\right) = 
({^t}l^p)^{a}_{\ b} \, \delta^b_{\ a}
=
({^t_0}X^e)^a_{\ A} \ ({^0}L^p)^A_{\ B} \left( ({^t_0}X^e)^{-1} \right)^B_{\ a}
=
\text{tr} \left( {^0}\underline{\underline{\boldsymbol{L}}}^p\right).
\end{equation}
Adopting these manipulations on the time derivative of the plastic Jacobian, inequality~\eqref{eq:CD_1} becomes
\begin{equation}
\frac{\partial \left( \bar{\rho} \, \psi \right)}{\partial t} 
+  
\bar{\rho} \, \psi \, ( {^0_0}\underline{\underline{\boldsymbol{X}}}^{p} )^{-T} \boldsymbol{:} \dfrac{\partial {^0_0}\underline{\underline{\boldsymbol{X}}}^{p}}{\partial t} 
= 
\frac{\partial \left( \bar{\rho} \, \psi \right)}{\partial t} 
+
\bar{\rho} \, \psi \, \text{tr} ( {^t}\underline{\underline{\boldsymbol{l}}}^p ) \geq 0,
\end{equation}
or, with the additive decomposition Eq.~\eqref{eq:velocity_gradient_additive_decomposition}, 
\begin{equation}
\label{eq:CD_2}
J^p \biggl(  \underline{\underline{\boldsymbol{\zeta}}} \boldsymbol{:} {^t} \underline{\underline{\boldsymbol{l}}}^e + 
\bigl( \underbrace{\underline{\underline{\boldsymbol{\zeta}}} -
\bar{\rho} \, \psi \, {^t}\underline{\underline{\boldsymbol{1}}}}_{\coloneqq \boldsymbol{\underline{\underline{\xi}}}^{\zeta}} \bigr)  \boldsymbol{:}
{^t}\underline{\underline{\boldsymbol{l}}}^p 
- 
\frac{\partial \left( \bar{\rho} \, \psi \right)}{\partial t} 
\biggr) 
\geq 0.
\end{equation}

As assumed in Section~\ref{subsec:constitutive_equations}, the elastic isotropy of the solid under consideration permits expressing its stored energy function in the intermediate configuration, i.e., via $\bar{\rho} \, \psi$.
Moreover, the assumption of uncoupled material permits the following simplification
\begin{equation}
\label{eq:uncoupled_material}
\bar{\rho} \, \psi \left( {^t}\underline{\underline{\boldsymbol{b}}}^e, z \right) = \check{\Psi} ( {^t}\underline{\underline{\boldsymbol{b}}}^e ) + 
\tilde{\Psi} \left( {^t}z \right),
\end{equation}
resulting in the trivial time derivative
\begin{equation}
\frac{ \partial \bar{\rho} \, \psi \left( {^t}\underline{\underline{\boldsymbol{b}}}^e, {^t}z \right)}{\partial t} 
= \frac{ \partial \Psi ( {^t}\underline{\underline{\boldsymbol{b}}}^e )}{\partial t} + 
\frac{ \partial \tilde{\Psi} \left( {^t}z \right)}{\partial t}
=
\frac{ \partial \Psi ( {^t}\underline{\underline{\boldsymbol{b}}}^e )}{\partial {^t}\underline{\underline{\boldsymbol{b}}}^e} 
\boldsymbol{:}
\frac{\partial {^t}\underline{\underline{\boldsymbol{b}}}^e}{\partial t} 
+ 
\frac{ \partial \tilde{\Psi} \left( {^t}z \right)}{\partial \, {^t}z} \frac{\partial \, {^t}z}{\partial t}
\end{equation}
If introducing the principal invariants of ${^t}\underline{\underline{\boldsymbol{b}}}^e$, elastic isotropy implies that $ \Psi \left( {^t}\underline{\underline{\boldsymbol{b}}}^e \right)
= \Psi \left( I_1^{b^e}, I_2^{b^e}, I_3^{b^e} \right)$, where
with
\begin{equation}
I_1^{b^e} 
\coloneqq 
\text{tr} ({^t}\underline{\underline{\boldsymbol{b}}}^e);
\qquad
I_2^{b^e} \coloneqq \dfrac{1}{2} \left( (I_1^{b^e})^2 - \text{tr} ({^t}\underline{\underline{\boldsymbol{b}}}^e \cdot {^t}\underline{\underline{\boldsymbol{b}}}^e) \right);
\qquad
I_3^{b^e} = \det [ ( {^t}b^e )^{a}_{ \ b} ].
\end{equation}
In light of the above formulae and the following kinematic relationship 
\begin{equation}
\frac{\partial {^t}\underline{\underline{\boldsymbol{b}}}^e }{\partial t} 
= 
{^t}\underline{\underline{\boldsymbol{l}}}^e \cdot {^t}\underline{\underline{\boldsymbol{b}}}^e + {^t}\underline{\underline{\boldsymbol{b}}}^e \cdot ({^t}\underline{\underline{\boldsymbol{l}}}^e)^T,
\end{equation}
inequality~\eqref{eq:CD_2} becomes
\begin{equation}
\label{eq:CD_3}
J^p \biggl(  {^t}\underline{\underline{\boldsymbol{\zeta}}} \boldsymbol{:} {^t} \underline{\underline{\boldsymbol{l}}}^e + 
\boldsymbol{\underline{\underline{\xi}}}^{\zeta}  \boldsymbol{:}
{^t}\underline{\underline{\boldsymbol{l}}}^p 
- 
2 \frac{ \partial \Psi ( {^t}\underline{\underline{\boldsymbol{b}}}^e )}{\partial {^t}\underline{\underline{\boldsymbol{b}}}^e} \boldsymbol{:} 
\left( {^t}\underline{\underline{\boldsymbol{b}}}^e
\cdot 
{^t}\underline{\underline{\boldsymbol{l}}}^e \right)
- 
\frac{ \partial \tilde{\Psi} \left( {^t}z \right)}{\partial \, {^t}z} \frac{\partial \, {^t}z}{\partial t}
\biggr)
\geq 0.
\end{equation}
Using the standard arguments for the Coleman-Noll procedure (see, e.g.,~\cite{coleman1967thermodynamics,coleman1974thermodynamics}), it follows that
\begin{equation}
\label{eq:Colemann_Noll_1}
\boldsymbol{\underline{\underline{\zeta}}}  = 
2 \dfrac{\partial \left(  \bar{\rho} \ \psi\right)}{\partial  {^t}\boldsymbol{\underline{\underline{b}}}^{e}} \cdot  {^t}\boldsymbol{\underline{\underline{b}}}^{e}
=
2 \left( 
\dfrac{\partial \left(  \bar{\rho} \ \psi\right)}{\partial I_1^{b^e} }
{^t}\boldsymbol{\underline{\underline{b}}}^{e}
+
\dfrac{\partial \left(  \bar{\rho} \ \psi\right)}{\partial I_2^{b^e} }
\left( I_1^{b^e} \, {^t}\boldsymbol{\underline{\underline{b}}}^{e} - {^t}\boldsymbol{\underline{\underline{b}}}^{e} \cdot {^t}\boldsymbol{\underline{\underline{b}}}^{e} \right)
+
\dfrac{\partial \left(  \bar{\rho} \ \psi\right)}{\partial I_3^{b^e}}  I_3^{b^e} \,{^t}\boldsymbol{\underline{\underline{1}}}
\right)
\end{equation}
By denoting with $\beta^{\zeta} $ the stress-like internal variable which is the power-conjugate of ${^t}z$, i.e.,
\begin{equation}
\label{eq:Colemann_Noll_2}
\beta^{\zeta} 
\coloneqq 
\dfrac{\partial \left(  \bar{\rho} \ \psi\right)}{\partial \, {^t}z} ,
\end{equation}
inequality~\eqref{eq:CD_3} reduces to
\begin{equation}
\label{eq:CD_23}
D \coloneqq J^p \biggl(  \boldsymbol{\underline{\underline{\xi}}}^{\zeta}  \boldsymbol{:}
{^t}\underline{\underline{\boldsymbol{d}}}^p 
- 
\beta^{\zeta} \, \frac{\partial \, {^t}z}{\partial t}
\biggr)
\geq 0,
\end{equation}
since  $\boldsymbol{\underline{\underline{\xi}}}^{\zeta}  \boldsymbol{:}
{^t}\underline{\underline{\boldsymbol{l}}}^p 
= \boldsymbol{\underline{\underline{\xi}}}^{\zeta}  \boldsymbol{:}
{^t}\underline{\underline{\boldsymbol{d}}}^p$, owing to the symmetry of $\boldsymbol{\underline{\underline{\xi}}}^{\zeta}$.
Invoking the principle of maximum dissipation, the maximum of the functional $D$ is sought for its arguments $\boldsymbol{\underline{\underline{\xi}}}^{\zeta}$ and $\beta^{\zeta}$. Adopting the method of Lagrange multipliers, the stationarity of the Lagrangian functional provides the following evolution laws
\begin{alignat}{2}
\label{eq:evolution_laws_1}
\boldsymbol{\underline{\underline{d}}}^p &
= & \dot{\gamma} \dfrac{\partial \Phi}{\partial \boldsymbol{\underline{\underline{\xi}}}^{\zeta}};
\\
\label{eq:evolution_laws_2}
\dot{z} & = & - \dot{\gamma} \dfrac{\partial \Phi}{\partial \beta^{\zeta} },
\end{alignat}
where the rate (denoted by $\left( \dot{\bullet} \right)$) of the plastic multiplier $\dot{\gamma}$ plays the role of the Lagrange multiplier, and $\Phi \left(\boldsymbol{\underline{\underline{\xi}}}^{\zeta}, \beta^{\zeta}   \right)$ is the yield function. 
This set of equations is complemented by the Karush-Kuhn-Tucker (KKT) and consistency conditions, expressed by
\begin{equation}
\dot{\gamma} \geq 0; 
\qquad
\Phi \leq 0; 
\qquad
\dot{\gamma} \, \Phi = 0; 
\qquad
\dot{\gamma} \, \dot{\Phi} = 0.
\end{equation}

\section{Linearisation}
\label{sec:linearisation}
\subsection{Linearisation of the internal force vector}
\label{subsec:linearisation_internal_force_vector}
Employing the stress defined in Eq.~\eqref{eq:semi_kirchhoff_stress}, the components of the \emph{internal force vector} are defined by considering the first term of on the RHS of Eq.~\eqref{eq:updated_lagrangian_discrete_equilibrium} and eliminating the coefficients of the test functions, i.e., 
\begin{equation}
\left( \mathrm{f}^{int} \right)_{\tilde{f}\mathsf{N}} 
\coloneqq
\int_{{^0}\mathcal{B}^{ax}} 
J^{p}
\underbrace{\left( 
\delta_{\tilde{f}}^a \ \dfrac{\partial \uppsi_\mathsf{N} }{\partial {^t{x}^b}}
+ 
\delta^{e}_{\tilde{f}} \ {^t{\gamma}^a_{be}} \ 
\ \uppsi_\mathsf{N}
\right)}_{\coloneqq \left(\delta_{\tilde{f}}^a \ \uppsi_\mathsf{N} \right)|_{b}}  
\
\ {\zeta_a^{\ b}} 
\ \text{d} {^0}V^{ax}.
\end{equation}
When linearising, as prescribed by Eq.~\eqref{eq:NR_updated_lagrangian}, the variation of the above vector are the only ones providing a contribution. 
In particular, we will consider the variation with respect to the displacements ${^h}\mleft( {^t}\boldsymbol{\underline{u}}^{ax} \mright)$ 
\footnote{
As already mention in Section~\ref{subsec:discretisation}, we remark that the displacement solution coming from the NR scheme is the ${^h}\mleft( {^t}\boldsymbol{\underline{u}}^{ax} \mright)$, while and the one required for the deformation gradient is the ${^h}\mleft( {^t}\boldsymbol{\underline{U}}^{ax} \mright)$ (see Eq.~\eqref{eq:deformation_gradient_components_material_displ}). However, as highlighted by Eq.~\eqref{eq:same_curvilinear_components} the shifter plays a silent role between these two vectors, and shifting between these components of displacement can be skipped without practical consequences code-wise.
}
. 
If the components of this variations are considered, the linearisation of the internal force vector is given by
\begin{equation}
\delta \left( \mathscr{G}^{ax} \left( {^h}\mleft( {^t_0}\boldsymbol{\underline{u}}^{ax} \mright); {^h}\mleft({^t}\boldsymbol{\underline{w}}^{ax}\mright) \right)  \right)
\cdot
\delta {^h}\mleft( {^t}\boldsymbol{\underline{u}}^{ax} \mright) 
=
{^t}\mathrm{w}^{\tilde{f}\mathsf{N}} \
\dfrac{
\partial \left( \mathrm{f}^{int} \right)_{\tilde{f}\mathsf{N}} 
}{
\partial {^t_0}\mathrm{u}^{\tilde{g}\mathsf{M}}
} 
\
\delta {^t_0}\mathrm{u}^{\tilde{g}\mathsf{M}}, 
\end{equation}
where, in particular,
\begin{align}
\label{eq:internal_force_vector_linearised_components}
\nonumber
\dfrac{\partial \left( \mathrm{f}^{int} \right)_{\tilde{f}\mathsf{N}} }{\partial {^t_0}\mathrm{u}^{\tilde{g}\mathsf{M}}} 
& = 
\dfrac{\partial}{\partial {^t_0}\mathrm{u}^{\tilde{g}\mathsf{M}}} 
\int_{{^0}\mathcal{B}^{ax}} J^{p}
\left(\delta_{\tilde{f}}^a \ \uppsi_\mathsf{N} \right)|_{b} 
\
\ {\zeta_a^{\ b}} 
\ \text{d} {^0}V^{ax}
\\
& =
\int_{^0\Omega} 
\dfrac{\partial}{\partial \left( {^t_0{X}}^{ax} \right)^{c}_{\ A}} \left( 
J^{p}
\left(\delta_{\tilde{f}}^a \ \uppsi_\mathsf{N} \right)|_{b} 
\
\ {\zeta_a^{\ b}} \right) \ 
\dfrac{
\partial \left( {^t_0{X}^{ax}} \right)^{c}_{A}
}{\partial {^t_0}\mathrm{u}^{\tilde{g}\mathsf{M}}} \ 
\text{d} {^0}V^{ax}.
\end{align}
The components of the second derivative appearing in the above equation rely on the incremental deformation gradient Eqs.~\eqref{eq:deformation_gradient_incrementally_computed}-\eqref{eq:increment_deformation_gradient_components_explicit} to be computed, i.e.,
\begin{align}
\nonumber
\dfrac{\partial \left( {{^t_0}X}^{ax} \right)^{c}_{A} }{\partial {^t_0}\mathrm{u}^{\tilde{g}\mathsf{M}}} 
& =
\dfrac{
\partial \left( {^t_{t_n}{X}}^{ax} \right)^{c}_{\mathcal{A}}
}{
\partial {^t_0}\mathrm{u}^{\tilde{g}\mathsf{M}}} 
\ \left( {^{t_n}_0{X}}^{ax} \right)^{\mathcal{A}}_{\ A}
\\
\nonumber
& = 
\dfrac{\partial }{\partial {^t_0}\mathrm{u}^{\tilde{g}\mathsf{M}}} 
\left( 
\left( {^t_{t_n}{S}}^{ax} \right)^c_{\ \mathcal{B}} \left( \delta^{\mathcal{B}}_{\ \mathcal{A}} + \frac{\partial { \left( {^t_0}{U}^{ax} \right)^\mathcal{B}}}{\partial {^{t_n}{x}^{\mathcal{A}}}} 
+ 
\left( {^{t_n}{\Gamma}}^{ax} \right)^{\mathcal{B}}_{\mathcal{AC}} \ {^{t}{U}^{\mathcal{C}}} \right)
\right) 
\ \left( {^{t_n}_0{X}}^{ax} \right)^{\mathcal{A}}_{\ A}
\\
& = 
\nonumber
\dfrac{\partial }{\partial {^t_0}\mathrm{u}^{\tilde{g}\mathsf{M}}} 
\left( 
{{^t_0}\mathrm{u}^{c\mathsf{O}}} \ \frac{\partial \uppsi_{\mathsf{O}}}{\partial {^{t_n}{x}^{\mathcal{A}}}}
+
\left( {^t_{t_n}{S}}^{ax} \right)^c_{\ \mathcal{B}} 
\ \left( {^{t_n}{\Gamma}}^{ax} \right)^{\mathcal{B}}_{\mathcal{AC}} 
\ \left( ^{t_n}_t{S}^{ax} \right)^{\mathcal{C}}_{\ d} 
\ {{^t_0}\mathrm{u}^{d\mathsf{O}}} \ \uppsi_{\mathsf{O}}
\right)
\ \left( {^{t_n}_0{X}}^{ax} \right)^{\mathcal{A}}_{\ A}
\\
\nonumber
& =
\left( \delta^c_{\tilde{g}} \ \delta^{\mathsf{O}}_{\mathsf{M}} \ 
\frac{\partial \uppsi_{\mathsf{O}}}{\partial {^{t_n}{x}^{\mathcal{A}}}}
+
\left( {^t_{t_n}{S}}^{ax} \right)^c_{\ \mathcal{B}} 
\ \left( {^{t_n}{\Gamma}^{ax}} \right)^{\mathcal{B}}_{\mathcal{AC}} 
\ \left( ^{t_n}_t{S}^{ax} \right)^{\mathcal{C}}_{\ d} 
\ \delta^d_{\tilde{g}} \ \delta^{\mathsf{O}}_{\mathsf{M}} \  \uppsi_{\mathsf{O}}
\right)
\ \left( {^{t_n}_0{X}}^{ax} \right)^{\mathcal{A}}_{\ A}
\\
& = 
\left( \delta^c_{\tilde{g}} \ 
\frac{\partial \uppsi_{\mathsf{M}}}{\partial {^{t}{x}^a}}
+
\left( {^t}{\gamma}^{ax} \right)^{c}_{a d} \
\delta^d_{\tilde{g}} \ \uppsi_{\mathsf{M}}
\right)
\left( {^t_{t_n}{X}}^{ax} \right)^{a}_{\mathcal{A}} \ 
\ \left( {^{t_n}_0{X}}^{ax} \right)^{\mathcal{A}}_{\ A}
\\
& = 
( \delta^c_{\tilde{g}} \ 
\uppsi_{\mathsf{M}} )|_{a} \ 
\left( {^t_{0}{X}}^{ax} \right)^{a}_{\ A}.
\end{align}
Owing to the above chain of formulae, Eq.~\eqref{eq:internal_force_vector_linearised_components} becomes
\begin{align}
\label{eq:internal_force_vector_linearised_components_1}
\nonumber
\dfrac{\partial \left( \mathrm{f}^{int} \right)_{\tilde{f}\mathsf{N}} }{\partial {^t_0}\mathrm{u}^{\tilde{g}\mathsf{M}}} 
& =
\int_{{^0}\mathcal{B}^{ax}} 
\dfrac{\partial}{\partial \left( {^t_0{X}}^{ax} \right)^{c}_{\ A} }  
\biggl( J^{p} \left(\delta_{\tilde{f}}^a \ \uppsi_\mathsf{N} \right)|_{b} \
\ {\zeta_a^{\ b}} \biggr) \ 
( \delta^c_{\tilde{g}} \ 
\uppsi_{\mathsf{M}} )|_{d} \ 
\left( {^t_{0}}{X}^{ax} \right)^{d}_{\ A} \
\text{d} {^0}V^{ax} 
\\
\nonumber
& =
\int_{{^0}\mathcal{B}^{ax}} 
\Biggl(
J^{p}
\dfrac{\partial}{\partial \left( {^t_0{X}}^{ax} \right)^{c}_{\ A} } \Bigl( \left(\delta_{\tilde{f}}^a \ \uppsi_\mathsf{N} \right)|_{b} \Bigr)
\ {\zeta_a^{\ b}} 
\\
\nonumber
& \qquad 
+
\left( \delta_{\tilde{f}}^a \ \uppsi_\mathsf{N} \right)|_{b} 
\biggl( J^{p}
\dfrac{\partial {\zeta_a^{\ b}} }{\partial \left( {^t_0{X}}^{ax} \right)^{c}_{\ A}}
+ 
{\zeta_a^{\ b}}
\dfrac{\partial J^{p} }{\partial \left( {^t_0{X}}^{ax} \right)^{c}_{\ A}}
\biggr)
\Biggr)
( \delta^c_{\tilde{g}} \ 
\uppsi_{\mathsf{M}} )|_{d} \ 
\left( {^t_{0}}{X}^{ax} \right)^{d}_{\ A} \
\text{d} {^0}V^{ax}
\\
\nonumber
& = 
\int_{{^0}\mathcal{B}^{ax}} 
\Biggl(
J^{p}
\biggl( \delta^a_{\tilde{f}} \ 
\frac{\partial \uppsi_{\mathsf{N}}}{\partial {^{0}{x}^B}}
+
\dfrac{\partial}{\partial {^{0}{x}^B}} \
\frac{\partial {^t{x}^a}}{\partial {^t{z}^j}} \
\frac{\partial {^{t}{z}^j}}{\partial {^{t}{x}^e}} \
\delta^e_{\tilde{f}} \ \uppsi_{\mathsf{N}} \biggr) 
\dfrac{\partial ( \left( {^t_0}{X}^{ax} \right)^{-1})^{B}_{\ b}
}{
\partial \left( {^t_0}{X}^{ax} \right)^{c}_{A}}
\ \zeta_a^{\ b} 
\\ 
\nonumber
& \qquad 
+
\left(\delta_{\tilde{f}}^a \ \uppsi_\mathsf{N} \right)|_{b} 
\biggl( J^{p}
\dfrac{\partial {\zeta_a^{\ b}} }{\partial \left( {^t_0{X}}^{ax} \right)^{c}_{\ A}}
+ 
{\zeta_a^{\ b}}
\dfrac{\partial J^{p} }{\partial \left( {^t_0{X}}^{ax} \right)^{c}_{\ A}}
\biggr)
\Biggr)
( \delta^c_{\tilde{g}} \ 
\uppsi_{\mathsf{M}} )|_{d} \ 
\left( {^t_{0}}{X}^{ax} \right)^{d}_{\ A} \
\text{d} {^0}V^{ax}
\\
\nonumber
& = 
\int_{{^0}\mathcal{B}^{ax}} 
\Biggl(
- J^{p} \biggl( \delta^a_{\tilde{f}} \ 
\frac{\partial \uppsi_{\mathsf{N}}}{\partial {^{0}{x}^B}}
+
\dfrac{\partial}{\partial {^{0}{x}^B}} \
\frac{\partial {^t{x}^a}}{\partial {^t{z}^j}} \
\frac{\partial {^{t}{z}^j}}{\partial {^{t}{x}^e}} \
\delta^e_{\tilde{f}} \ \uppsi_{\mathsf{N}} \biggr) \ 
( \left( {^t_0}{X}^{ax} \right)^{-1})^{B}_{\ c} \
( \left( {^t_0}{X}^{ax} \right)^{-1})^{A}_{\ b} \ 
\zeta_a^{\ b}
\\ 
\nonumber
& \qquad 
+
\left(\delta_{\tilde{f}}^a \ \uppsi_\mathsf{N} \right)|_{b} 
\biggl( J^{p}
\dfrac{\partial {\zeta_a^{\ b}} }{\partial \left( {^t_0{X}}^{ax} \right)^{c}_{\ A}}
+ 
{\zeta_a^{\ b}}
\dfrac{\partial J^{p} }{\partial \left( {^t_0{X}}^{ax} \right)^{c}_{\ A}}
\biggr)
\Biggr) 
( \delta^c_{\tilde{g}} \ 
\uppsi_{\mathsf{M}} )|_{d} \ 
\left( {^t_{0}}{X}^{ax} \right)^{d}_{\ A} \
\text{d} {^0}V^{ax}
\\
\nonumber
& =
\int_{{^0}\mathcal{B}^{ax}} 
\Biggl(
- J^{p} \ ( \delta^a_{\tilde{f}} \ 
\uppsi_{\mathsf{N}} )|_{c} \
\delta^d_{b} \ 
\zeta_a^{\ b} 
\\
\nonumber
& \qquad \qquad
+
\left(\delta_{\tilde{f}}^a \ \uppsi_\mathsf{N} \right)|_{b} 
\biggl( J^{p}
\dfrac{\partial {\zeta_a^{\ b}} }{\partial \left( {^t_0{X}}^{ax} \right)^{c}_{\ A}}
+ 
{\zeta_a^{\ b}}
\dfrac{\partial J^{p} }{\partial \left( {^t_0{X}}^{ax} \right)^{c}_{\ A}}
\biggr)
\left( {^t_{0}}{X}^{ax} \right)^{d}_{\ A}
\Biggr) 
( \delta^c_{\tilde{g}} \ 
\uppsi_{\mathsf{M}} )|_{d} \ 
\
\text{d} {^0}V^{ax}
\\ 
& =
\footnotesize{
\int_{{^0}\mathcal{B}^{ax}} 
\left(\delta_{\tilde{f}}^a \ \uppsi_\mathsf{N} \right)|_{b} 
\Biggl(
\underbrace{ - J^{p} \
\zeta_a^{\ d} \ \delta_c^{b}
+
\biggl( J^{p}
\dfrac{\partial {\zeta_a^{\ b}} }{\partial \left( {^t_0{X}}^{ax} \right)^{c}_{\ A}}
+ 
{\zeta_a^{\ b}}
\dfrac{\partial J^{p} }{\partial \left( {^t_0{X}}^{ax} \right)^{c}_{\ A}}
\biggr)
\left( {^t_{0}}{X}^{ax} \right)^{d}_{\ A} }_{\coloneqq {^t}a_{a \ c}^{\ b \ d}}
\Biggr)
( \delta^c_{\tilde{g}} \ 
\uppsi_{\mathsf{M}} )|_{d} \
\text{d} {^0}V^{ax}
},
\end{align}
where the spatial fourth-order tensor ${^t}\underline{\underline{\underline{\underline{\boldsymbol{a}}}}}$ per unit original volume has also been introduced. 
The two derivative appearing in this tensor can be computed making use of the chain rule. For the first term, this is given by
\begin{align}
\nonumber
\dfrac{\partial {\zeta}_a^{\ b} }{\partial \left( {^t_0}{X}^{ax} \right)^{c}_{A}} \ \left( {^t_{0}}{X}^{ax} \right)^{d}_{\ A}
& = 
\dfrac{\partial {\zeta_a^{\ b}} }{\partial {^t{\underline{\underline{\boldsymbol{\epsilon}}}^{e, tr}}}}
\boldsymbol{:} 
\dfrac{\partial {^t{\underline{\underline{\boldsymbol{\epsilon}}}}^{e, tr}}}{\partial {^t\underline{\underline{\boldsymbol{b}}}^{e, tr}}}
\boldsymbol{:}
\dfrac{\partial {^t\underline{\underline{\boldsymbol{b}}}^{e, tr}}}{\partial \left( {^t_0}{X}^{ax} \right)^{c}_{\ A}}
\ \left( {^t_{0}}{X}^{ax} \right)^{d}_{\ A}
\\
& = 
\frac{1}{2} \left( {^t}\mathcal{D}^{alg} \right)^{\ b \ f}_{a \ e} \ {^t}\mathcal{L}^{e \ \ h}_{\ f g} \ {^t}\mathcal{B}^{g \ \ d}_{\ h c},
\end{align}
where the components of the fourth-order tensors appearing in the last row are defined as follows
\begin{gather}
\left( {^t}\mathcal{D}^{alg} \right)^{\ b \ f}_{a \ e} 
\coloneqq 
\dfrac{\partial \zeta_a^{\ b} }{\partial {( {^t_0}{\epsilon}^{e, tr}})^e_{\ f}};
\\
{^t} \mathcal{L}^{e \ \ h}_{\ f g} 
\coloneqq 
\dfrac{\partial \hbox{ln} ( {^t}b^{e,tr} )^e_{\ f} }{\partial ( {^t_0}b^{e,tr} )^g_{\ h}};
\\
\label{eq:fourth_order_B}
{^t}\mathcal{B}^{g \ \ d}_{\ h c} 
\coloneqq
\dfrac{\partial ( {^t}b^{e,tr} )^g_{\ h}}{\partial \left( {^t_0}{X}^{ax} \right)^{c}_{\ A} } \ \left( {^t_0}{X}^{ax} \right)^{d}_{\ A}
= 
\delta^g_c \ ( {^t}b^{e,tr} )^d_{\ h} + ( {^t}b^{e,tr} )^{gd} \ {^t{g}_{hc}}.
\end{gather}
In particular, the derivation of ${^t}\underline{\underline{\underline{\underline{\mathbfcal{D}}}}}^{alg}$ is provided in the following section of this appendix and it requires the constitutive relationship.

As already mentioned in Section~\ref{subsec:elasto-plasticity}, the derivatives of the plastic part of the Jacobian are not straightforward to computed and, to the authors' knowledge, this topic has never been addressed. 
This manuscript suggests of exploiting the relationship ${^t_0}J =  J^e \,
J^p$, so that the derivative of the plastic part of the Jacobian can be divided into two parts as follows
\begin{equation}
\dfrac{\partial J^{p} }{\partial \left( {^t_0{X}}^{ax} \right)^{c}_{\ A}}
=
\frac{1}{ J^e}\dfrac{\partial \, {^t_0}J}{\partial \left( {^t_0{X}}^{ax} \right)^{c}_{\ A}}
+
{^t_0}J \dfrac{\partial \left(  J^e \right)^{-1}}{\partial \left( {^t_0{X}}^{ax} \right)^{c}_{\ A}}.
\end{equation}
The first derivative appearing in the above equations can be derived from standard relationships between kinematic quantities, while the second one requires the constitutive relationship and can be developed as follows
\begin{align}
\label{eq:derivative_plastic_Jacobian}
\dfrac{\partial J^{p} }{\partial \left( {^t_0{X}}^{ax} \right)^{c}_{\ A}}
& =
\nonumber
\frac{ {^t_0} J }{  J^e } \left( ( {^t_0{X}}^{ax} )^{-1} \right)^{A}_{\ \ c}
+
 {^t_0}J
\dfrac{\partial \exp \left( -{^t}\epsilon^{e}_v \right)}{\partial {^t}\epsilon^{e}_v }
\dfrac{\partial {^t}\epsilon^{e}_v }{\partial {^t}\boldsymbol{\underline{\underline{\epsilon}}}^e}
\boldsymbol{:}
\dfrac{\partial {^t}\boldsymbol{\underline{\underline{\epsilon}}}^e}{\partial {^t}\boldsymbol{\underline{\underline{\epsilon}}}^{e, tr}}
\boldsymbol{:}
\dfrac{ \partial {^t}\boldsymbol{\underline{\underline{\epsilon}}}^{e, tr} }{ \partial {^t\underline{\underline{\boldsymbol{b}}}^{e, tr}} }
\boldsymbol{:}
\dfrac{ \partial {^t\underline{\underline{\boldsymbol{b}}}^{e, tr}} }{ \partial \left( {^t_0}{X}^{ax} \right)^{c}_{\ A} }
\\
& =
J^{p}
\left( 
\left( ( {^t_0{X}}^{ax} )^{-1} \right)^{A}_{\ \ c}
-
\dfrac{1}{2}
\delta_{m}^{n} \
\dfrac{\partial \left( {^t}\epsilon^e \right)^{m}_{\ n}}{\partial \left( {^t}\epsilon^{e, tr} \right)^{e}_{ \ f} }
{^t} \mathcal{L}^{e \ \ h}_{\ f g} 
\dfrac{\partial ( {^t}b^{e,tr} )^g_{\ h}}{\partial \left( {^t_0}{X}^{ax} \right)^{c}_{\ A} }
\right).
\end{align}
The calculation of the fourth-order $\dfrac{\partial {^t}\boldsymbol{\underline{\underline{\epsilon}}}^e}{\partial {^t}\boldsymbol{\underline{\underline{\epsilon}}}^{e, tr}}$ is explained in the next section of the appendix as it relies on the the constitutive subroutine.
In light of the above equations, the tensor ${^t}\underline{\underline{\underline{\underline{\boldsymbol{a}}}}}$ can be expressed as
\begin{align}
\nonumber
{^t_0}a_{a \ c}^{\ b \ d}
& =
J^{p}
\biggl( 
- \
\zeta_a^{\ d} \ \delta_c^{b}
+
\zeta_a^{\ b} \ \delta_{c}^{d}
+
\dfrac{1}{2}
\left(
( {^t}\mathcal{D}^{alg} )^{\ b \ f}_{a \ e} 
-
\zeta_a^{\ b} \
\delta_{m}^{n} \
\dfrac{\partial \left( {^t}\epsilon^e \right)^{m}_{\ n}}{\partial \left( {^t}\epsilon^{e, tr} \right)^{e}_{ \ f} }
\right)
{^t}\mathcal{L}^{e \ \ h}_{\ f g} \ {^t}\mathcal{B}^{g \ \ d}_{\ h c}
\biggr)
\\
& =
{^t_0}J
\biggl( 
- \
\sigma_a^{\ d} \ \delta_c^{b}
+
\sigma_a^{\ b} \ \delta_{c}^{d}
+
\dfrac{1}{2 J^e}
\left(
( {^t}\mathcal{D}^{alg} )^{\ b \ f}_{a \ e} 
-
\sigma_a^{\ b} \
\delta_{m}^{n} \
\dfrac{\partial \left( {^t}\epsilon^e \right)^{m}_{\ n}}{\partial \left( {^t}\epsilon^{e, tr} \right)^{e}_{ \ f} }
\right)
{^t}\mathcal{L}^{e \ \ h}_{\ f g} \ {^t}\mathcal{B}^{g \ \ d}_{\ h c}
\biggr).
\end{align}
It can be appreciated from the above spatial tensor that, in the case of zero plastic volumetric deformations (e.g., J-2 plasticity), the second and fourth terms in the above equation are simplified, and the spatial tensor is equal to that obtained for a standard Hencky elasto-plastic material (see, for instance,~\cite{desouza2011computational}).

\subsection{Linearisation of the constitutive relationship}
\label{subsec:linearisation_constitutive_equations}

As a standard starting point for the constitutive subroutine, the following equations are considered: the time-discrete versions of the additive decomposition of the logarithmic strain; the time-discrete version of Eq.~\eqref{eq:evolution_laws_2}; and the fulfilment with the equal sign of the yield function, i.e., respectively
\begin{equation}
\begin{dcases}
 {^t}\boldsymbol{\underline{\underline{\epsilon}}}^{e} -  {^t}\boldsymbol{\underline{\underline{\epsilon}}}^{e, tr} 
+ \Delta \gamma \dfrac{\partial \Phi}{\partial \boldsymbol{\underline{\underline{\xi}}}^{\zeta}} 
= \boldsymbol{0};
\\
{^t}z - {^{t_n}}z + \Delta {\gamma} \dfrac{\partial \Phi}{\partial \beta^{\zeta} } = 0;
\\
\Phi \left( \boldsymbol{\underline{\underline{\xi}}}^{\zeta}, \beta^{\zeta}\right) = 0,
\end{dcases}
\end{equation}
where ${^t}\boldsymbol{\underline{\underline{\epsilon}}}^{e, tr}  = {^t}\boldsymbol{\underline{\underline{\epsilon}}}$ refers to the predictor trial logarithmic strain. 
In particular, the first equation in the above system stems from Eq.~\eqref{eq:evolution_laws_2} and the exponential map integrator applied to the time-discretisation of $ {^0}\underline{\underline{\boldsymbol{L}}}^p = \frac{\partial {^0_0}\underline{\underline{\boldsymbol{X}}}^p}{\partial t} 
\cdot
({^t_0}\underline{\underline{\boldsymbol{X}}}^p)^{-1}$, which is then substituted in the multiplicative decomposition Eq.~\eqref{eq:elasto-plastic_def_gradient} (see, for instance,~\cite{desouza2011computational} for more implementation and derivation details).

Linearising the above set of equations with respect to the unknowns ${^t}\boldsymbol{\underline{\underline{\epsilon}}}^{e}, \Delta \gamma, \beta^{\zeta} $ and the trail state ${^t}\boldsymbol{\underline{\underline{\epsilon}}}^{e,tr}$ gives
\begin{equation}
\label{eq:linear_system_0}
\begin{dcases}
\text{d} {^t}\boldsymbol{\underline{\underline{\epsilon}}}^{e} 
- 
\text{d}{^t}\boldsymbol{\underline{\underline{\epsilon}}}^{e,tr} 
+ 
\text{d}\Delta \gamma \dfrac{\partial \Phi}{\partial \boldsymbol{\underline{\underline{\xi}}}^{\zeta}} 
+
\Delta \gamma 
\left( 
\dfrac{\partial^2 \Phi}{\partial \boldsymbol{\underline{\underline{\xi}}}^{\zeta} \ \partial \boldsymbol{\underline{\underline{\xi}}}^{\zeta}} 
\boldsymbol{:} \text{d} \boldsymbol{\underline{\underline{\xi}}}^{\zeta}
+
\dfrac{\partial^2 \Phi}{\partial \boldsymbol{\underline{\underline{\xi}}}^{\zeta} \ \partial \beta^{\zeta}} \text{d} \beta^{\zeta} 
\right)
= \boldsymbol{0};
\\
\text{d}  ^{t}z + \text{d}\Delta {\gamma} \dfrac{\partial \Phi}{\partial \beta^{\zeta} }
+
\Delta {\gamma} 
\left( 
\dfrac{\partial^2 \Phi}{\partial \beta^{\zeta} \ \partial \boldsymbol{\underline{\underline{\xi}}}^{\zeta} } \boldsymbol{:} \text{d}\boldsymbol{\underline{\underline{\xi}}}^{\zeta}
+
\dfrac{\partial^2 \Phi}{\left(\partial  \beta^{\zeta}\right)^2} \text{d}\beta^{\zeta}
\right)
= 0;
\\
\frac{\partial \Phi}{\partial \boldsymbol{\underline{\underline{\xi}}}^{\zeta} } 
\boldsymbol{:}
\text{d}\boldsymbol{\underline{\underline{\xi}}}^{\zeta}
+
\frac{\partial \Phi}{\partial \beta^{\zeta} } \text{d}\beta^{\zeta} = 0.
\end{dcases}
\end{equation}
To relate the increments in elastic strain with those of stresses and, similarly, the increments of strain-like internal variables with those of stress-like internal variables, the constitutive relationship Eq.~\eqref{eq:Colemann_Noll_1} and the definition Eq.~\eqref{eq:Colemann_Noll_2} are linearised too, yielding 
\begin{equation}
\label{eq:linear_system_1}
\begin{dcases}
\text{d}\boldsymbol{\underline{\underline{\zeta}}} 
=
\dfrac{\partial^2 \left(  \bar{\rho} \ \psi\right)}{\partial  {^t}\boldsymbol{\underline{\underline{\epsilon}}}^{e} \ \partial  {^t}\boldsymbol{\underline{\underline{\epsilon}}}^{e} } \boldsymbol{:}
\text{d} {^t}\boldsymbol{\underline{\underline{\epsilon}}}^{e} 
+
\dfrac{\partial^2 \left(  \bar{\rho} \ \psi\right)}{\partial  {^t}\boldsymbol{\underline{\underline{\epsilon}}}^{e} \ \partial z }
\text{d}z;
\\
\text{d}\beta^{\zeta} =
\dfrac{\partial^2 \left(  \bar{\rho} \ \psi\right)}{ \left( \partial z \right)^2} \text{d} z 
+
\dfrac{\partial^2 \left(  \bar{\rho} \ \psi\right)}{\partial z \ \partial  {^t}\boldsymbol{\underline{\underline{\epsilon}}}^{e} } 
\boldsymbol{:} \text{d}  {^t}\boldsymbol{\underline{\underline{\epsilon}}}^{e};
\end{dcases}
\end{equation}
The equations are kept general in the above case, but, in the case of uncoupled material (see Eq.~\eqref{eq:uncoupled_material}), it can be straightforwardly seen that the mixed derivatives are zero. 
Keeping these equations general and assuming enough continuity of the stored energy function to guarantee well-behaved second derivatives, the above linear system can be inverted, which gives
\begin{equation}
\label{eq:linear_system_2}
\begin{cases}
\text{d}  {^t}\boldsymbol{\underline{\underline{\epsilon}}}^{e} = 
\boldsymbol{\underline{\underline{\underline{\underline{C}}}}} 
\boldsymbol{:} 
\text{d}\boldsymbol{\underline{\underline{\zeta}}} 
+
\boldsymbol{\underline{\underline{B}}} \ \text{d}\beta^{\zeta};
\\
\text{d} z = 
\boldsymbol{\underline{\underline{A}}}  \boldsymbol{:} 
\text{d}\boldsymbol{\underline{\underline{\zeta}}} 
+
J \ \text{d}\beta^{\zeta}.
\end{cases}
\end{equation}
Employing the definition of the Eshelby-zeta stress provided in~\eqref{eq:CD_2}, its linearised form is given by
\begin{align}
\label{eq:increment_eshelby}
\nonumber
\text{d} \boldsymbol{\underline{\underline{\xi}}}^{\zeta} 
& =
\text{d} \boldsymbol{\underline{\underline{\zeta}}} - \left( \dfrac{\partial\left(  \bar{\rho} \ \psi \right)}{\partial  {^t}\boldsymbol{\underline{\underline{\epsilon}}}^{e}} 
\boldsymbol{:}
\text{d} {^t}\boldsymbol{\underline{\underline{\epsilon}}}^{e}
+
\dfrac{\partial\left(  \bar{\rho} \ \psi \right)}{\partial z} 
\text{d}z
\right){^t}\boldsymbol{\underline{\underline{1}}}
\\
\nonumber
& =
\text{d} \boldsymbol{\underline{\underline{\zeta}}} 
- 
\left( \boldsymbol{\underline{\underline{\zeta}}} 
\boldsymbol{:}
\text{d} {^t}\boldsymbol{\underline{\underline{\epsilon}}}^{e}
+
\beta^{\zeta} \
\text{d}z
\right){^t}\boldsymbol{\underline{\underline{1}}}
\\
\nonumber
& = 
\text{d} \boldsymbol{\underline{\underline{\zeta}}} 
- 
\left( 
\boldsymbol{\underline{\underline{\zeta}}} 
\boldsymbol{:}
\boldsymbol{\underline{\underline{\underline{\underline{C}}}}} 
\boldsymbol{:} 
\text{d}\boldsymbol{\underline{\underline{\zeta}}} 
+
\boldsymbol{\underline{\underline{\zeta}}}
\boldsymbol{:}
\boldsymbol{\underline{\underline{B}}} \ \text{d}\beta^{\zeta} 
+
\beta^{\zeta} \
\boldsymbol{\underline{\underline{A}}}  \boldsymbol{:} 
\text{d}\boldsymbol{\underline{\underline{\zeta}}} 
+ 
\beta^{\zeta} \
J \ \text{d}\beta^{\zeta}
\right)
{^t}\boldsymbol{\underline{\underline{1}}}
\\
& =
\underbrace{ \left( {^t}\boldsymbol{{\underline{\underline{\underline{\underline{1}}}}}}^{4,sym}
-
{^t}\boldsymbol{\underline{\underline{1}}}
\
\boldsymbol{\underline{\underline{\zeta}}} 
\boldsymbol{:}
\boldsymbol{\underline{\underline{\underline{\underline{C}}}}} 
-
\beta^{\zeta} 
\
{^t}\boldsymbol{\underline{\underline{1}}}
\
\boldsymbol{\underline{\underline{A}}} 
\right) 
}_{\coloneqq
\underline{\underline{\underline{\underline{\mbox{\textbf{\varstigma}} }}}}
}
\boldsymbol{:} 
\text{d}\boldsymbol{\underline{\underline{\zeta}}} 
\underbrace{
-
{^t}\boldsymbol{\underline{\underline{1}}}
\left( 
\boldsymbol{\underline{\underline{\zeta}}}
\boldsymbol{:}
\boldsymbol{\underline{\underline{B}}}  
+
\beta^{\zeta} \
J
\right)
}_{\coloneqq
\boldsymbol{\underline{\underline{\mbox{\textbf{\sampi}}}}}
}
\text{d}\beta^{\zeta},
\end{align}
where also Eq.\eqref{eq:linear_system_2} has been used.

In light of Eqs.~\eqref{eq:linear_system_2}-\eqref{eq:increment_eshelby}, the linear system~\eqref{eq:linear_system_0} can be coincisely expressed as follows
\begin{equation}
\label{eq:linear_system_stress_solve}
\underbrace{
\begin{bmatrix}
\boldsymbol{\underline{\underline{\underline{\underline{E}}}}}_{11} 
& 
\boldsymbol{\underline{\underline{E}}}_{12}
& 
E_{13}
\\
\boldsymbol{\underline{\underline{E}}}_{21}
&
E_{22}
&
E_{23}
\\
\boldsymbol{\underline{\underline{E}}}_{31}
&
E_{32}
&
0
\end{bmatrix}
}_{\coloneqq [A]}
\underbrace{
\begin{bmatrix}
\text{d}\boldsymbol{\underline{\underline{\zeta}}} 
\\
\text{d}\beta^{\zeta}
\\
d\Delta \gamma
\end{bmatrix}
}_{\coloneqq [x]}
=
\underbrace{
\begin{bmatrix}
\text{d} {^t}\boldsymbol{\underline{\underline{\epsilon}}}^{e, tr} 
\\
0
\\
0
\end{bmatrix}
}_{\coloneqq [b]}
,
\end{equation}
where
\begin{align}
\boldsymbol{\underline{\underline{\underline{\underline{E}}}}}_{11} 
& = 
\boldsymbol{\underline{\underline{\underline{\underline{C}}}}} 
+
\Delta \gamma 
\left( 
\dfrac{\partial^2 \Phi}{\partial \boldsymbol{\underline{\underline{\xi}}}^{\zeta} \ \partial \boldsymbol{\underline{\underline{\xi}}}^{\zeta}} 
\boldsymbol{:}
\underline{\underline{\underline{\underline{\mbox{\textbf{\varstigma}} }}}}
\right);
\\
\boldsymbol{\underline{\underline{E}}}_{12}
& =
\boldsymbol{\underline{\underline{B}}}
+
\Delta \gamma
\left( 
\dfrac{\partial^2 \Phi}{\partial \boldsymbol{\underline{\underline{\xi}}}^{\zeta} \ \partial \boldsymbol{\underline{\underline{\xi}}}^{\zeta}} 
\boldsymbol{:}
\boldsymbol{\underline{\underline{\mbox{\textbf{\sampi}}}}}
+
\dfrac{\partial^2 \Phi}{\partial \boldsymbol{\underline{\underline{\xi}}}^{\zeta} \ \partial \beta^{\zeta}} 
\right);
\\
\boldsymbol{\underline{\underline{E}}}_{13} 
& = 
\dfrac{\partial \Phi}{\partial \boldsymbol{\underline{\underline{\xi}}}^{\zeta}};
\\
\boldsymbol{\underline{\underline{E}}}_{21}
& =
\boldsymbol{\underline{\underline{A}}} 
+
\Delta \gamma
\left( 
\dfrac{\partial^2 \Phi}{\partial \beta^{\zeta} \ \partial \boldsymbol{\underline{\underline{\xi}}}^{\zeta} } \boldsymbol{:}
\boldsymbol{\underline{\underline{\underline{\underline{\mbox{\textbf{\varstigma}}}}}}}
\right);
\\
E_{22} 
& = 
J 
+
\Delta \gamma 
\left( 
\dfrac{\partial^2 \Phi}{\partial \beta^{\zeta} \ \partial \boldsymbol{\underline{\underline{\xi}}}^{\zeta} }
\boldsymbol{:}
\boldsymbol{\underline{\underline{\mbox{\textbf{\sampi}}}}}
+
\dfrac{\partial^2 \Phi}{\partial \left( \beta^{\zeta}\right)^2} 
\right);
\\
E_{23} 
& = 
\dfrac{\partial \Phi}{\partial \beta^{\zeta}};
\\
\boldsymbol{\underline{\underline{E}}}_{31}
& =
\frac{\partial \Phi}{\partial \boldsymbol{\underline{\underline{\xi}}}^{\zeta} } 
\boldsymbol{:}
\boldsymbol{\underline{\underline{\underline{\underline{\mbox{\textbf{\varstigma}}}}}}}
;
\\
E_{32} 
& = 
\frac{\partial \Phi}{\partial \boldsymbol{\underline{\underline{\xi}}}^{\zeta} } 
\boldsymbol{:}
\boldsymbol{\underline{\underline{\mbox{\textbf{\sampi}}}}}
+
\frac{\partial \Phi}{\partial \beta^{\zeta} }.
\end{align}
In particular, the fourth-order tensor corresponding to the top-left entries of the inverse of the above matrix $[A^{-1}]$, denoted by $[A^{^-1}]_{11}$, provides the relationship ${^t}\underline{\underline{\underline{\underline{\mathbfcal{D}}}}}^{alg} \coloneqq \dfrac{\partial \boldsymbol{\underline{\underline{\zeta}}} }{\partial  {^t}\boldsymbol{\underline{\underline{\epsilon}}}^{e, tr} } = [A^{^-1}]_{11}$ necessary for the linearisation of the constitutive relationship.

A linear system analogous to that in Eq.~\eqref{eq:linear_system_stress_solve} is the following one
\begin{equation}
\label{eq:linear_system_strain_solve}
\underbrace{
\begin{bmatrix}
\check{\boldsymbol{\underline{\underline{\underline{\underline{E}}}}}}_{11} 
& 
\check{\boldsymbol{\underline{\underline{E}}}}_{12}
& 
\check{\boldsymbol{\underline{\underline{E}}}}_{13}
\\
\check{\boldsymbol{\underline{\underline{E}}}}_{21}
&
\check{E}_{22}
&
\check{E}_{23}
\\
\check{\boldsymbol{\underline{\underline{E}}}}_{31}
&
\check{E}_{32}
&
0
\end{bmatrix}
}_{\coloneqq [\check{A}]}
\underbrace{
\begin{bmatrix}
\text{d}  {^t}\boldsymbol{\underline{\underline{\epsilon}}}^{e}
\\
\text{d} z
\\
d\Delta \gamma
\end{bmatrix}
}_{\coloneqq [\check{x}]}
=
\begin{bmatrix}
\text{d} {^t}\boldsymbol{\underline{\underline{\epsilon}}}^{e, tr} 
\\
0
\\
0
\end{bmatrix},\end{equation}
where, in this case
\begin{align}
\check{\boldsymbol{\underline{\underline{\underline{\underline{E}}}}}}_{11} 
& =
{^t}\boldsymbol{{\underline{\underline{\underline{\underline{1}}}}}}^{4,sym}
+
\Delta \gamma 
\left( 
\dfrac{\partial^2 \Phi}{\partial \boldsymbol{\underline{\underline{\xi}}}^{\zeta} \ \partial \boldsymbol{\underline{\underline{\xi}}}^{\zeta}} 
\boldsymbol{:}
\underline{\underline{\underline{\underline{\mbox{\textbf{\stigma}} }}}}
+
\dfrac{\partial^2 \Phi}{\partial \boldsymbol{\underline{\underline{\xi}}}^{\zeta} \ \partial \beta^{\zeta}} 
\dfrac{\partial^2 \left(  \bar{\rho} \ \psi\right)}{\partial z \ \partial  {^t}\boldsymbol{\underline{\underline{\epsilon}}}^{e} } 
\right)
;
\\
\check{\boldsymbol{\underline{\underline{E}}}}_{12}
& =
\Delta \gamma 
\left( 
\dfrac{\partial^2 \Phi}{\partial \boldsymbol{\underline{\underline{\xi}}}^{\zeta} \ \partial \boldsymbol{\underline{\underline{\xi}}}^{\zeta}} 
\boldsymbol{:}
\boldsymbol{\underline{\underline{\Xi}} } 
+
\dfrac{\partial^2 \left(  \bar{\rho} \ \psi\right)}{ \left( \partial z \right)^2}
\dfrac{\partial^2 \Phi}{\partial \boldsymbol{\underline{\underline{\xi}}}^{\zeta} \ \partial \beta^{\zeta}} 
\right)
;
\\
\check{\boldsymbol{\underline{\underline{E}}}}_{13}
& =
\dfrac{\partial \Phi}{\partial \boldsymbol{\underline{\underline{\xi}}}^{\zeta}};
\\
\check{\boldsymbol{\underline{\underline{E}}}}_{21}
& =
\Delta \gamma 
\left( 
\dfrac{\partial^2 \Phi}{\partial \beta^{\zeta} \ \partial \boldsymbol{\underline{\underline{\xi}}}^{\zeta} } \boldsymbol{:} 
\underline{\underline{\underline{\underline{\mbox{\textbf{\stigma}} }}}}
+
\dfrac{\partial^2 \Phi}{\left(\partial  \beta^{\zeta}\right)^2} 
\dfrac{\partial^2 \left(  \bar{\rho} \ \psi\right)}{\partial z \ \partial  {^t}\boldsymbol{\underline{\underline{\epsilon}}}^{e} }
\right)
;
\\
\check{E}_{22}
& =
1 + \Delta \gamma 
\left( 
\dfrac{\partial^2 \Phi}{\partial \beta^{\zeta} \ \partial \boldsymbol{\underline{\underline{\xi}}}^{\zeta} } \boldsymbol{:} 
\boldsymbol{\underline{\underline{\Xi}} }
+
\dfrac{\partial^2 \Phi}{\left(\partial  \beta^{\zeta}\right)^2}
\dfrac{\partial^2 \left(  \bar{\rho} \ \psi\right)}{ \left( \partial z \right)^2}
\right)
;
\\
\check{E}_{23}
& =
\dfrac{\partial \Phi}{\partial \beta^{\zeta} }
\\
\check{\boldsymbol{\underline{\underline{E}}}}_{31}
& =
\frac{\partial \Phi}{\partial \boldsymbol{\underline{\underline{\xi}}}^{\zeta} } 
\boldsymbol{:} \underline{\underline{\underline{\underline{\mbox{\textbf{\stigma}} }}}}
+
\frac{\partial \Phi}{\partial \beta^{\zeta} }
\dfrac{\partial^2 \left(  \bar{\rho} \ \psi\right)}{\partial z \ \partial  {^t}\boldsymbol{\underline{\underline{\epsilon}}}^{e} };
\\
\check{E}_{32}
& =
\frac{\partial \Phi}{\partial \boldsymbol{\underline{\underline{\xi}}}^{\zeta} } 
\boldsymbol{:}
\boldsymbol{\underline{\underline{\Xi}} }
+
\frac{\partial \Phi}{\partial \beta^{\zeta} } 
\dfrac{\partial^2 \left(  \bar{\rho} \ \psi\right)}{ \left( \partial z \right)^2};
\end{align}
The quantities $\underline{\underline{\underline{\underline{\mbox{\textbf{\stigma}} }}}}$ and $\boldsymbol{\underline{\underline{\Xi}} }$ appearing in the above list of equations are defined as follows
\begin{align}
\underline{\underline{\underline{\underline{\mbox{\textbf{\stigma}} }}}}
& \coloneqq
\dfrac{\partial^2 \left(  \bar{\rho} \ \psi \right)}{\partial  {^t}\boldsymbol{\underline{\underline{\epsilon}}}^{e} \ \partial  {^t}\boldsymbol{\underline{\underline{\epsilon}}}^{e} } 
-
{^t}\boldsymbol{\underline{\underline{1}}} \ {^t}\boldsymbol{\underline{\underline{\zeta}}}
;
\\
\boldsymbol{\underline{\underline{\Xi}} }
&  
\coloneqq
\dfrac{\partial^2 \left(  \bar{\rho} \ \psi \right)}{\partial  {^t}\boldsymbol{\underline{\underline{\epsilon}}}^{e} \ \partial z }
-
\beta^{\zeta} \ {^t}\boldsymbol{\underline{\underline{1}}} 
,
\end{align}
and they permit expressing the increment of Eshelby-zeta stress concisely as 
\begin{equation}
\text{d} \boldsymbol{\underline{\underline{\xi}}}^{\zeta} 
=
\underline{\underline{\underline{\underline{\mbox{\textbf{\stigma}} }}}}
\boldsymbol{:}
\text{d} {^t}\boldsymbol{\underline{\underline{\epsilon}}}^{e}
+
\boldsymbol{\underline{\underline{\Xi}} } \ \text{d}z.
\end{equation}
The main advantage of the linear system~\eqref{eq:linear_system_strain_solve} over the linear system~\eqref{eq:linear_system_stress_solve} consists in the better conditioning of the matrix $[\check{A}]$ since the vector $[\check{x}]$ contains only strain or strain-like quantities, while this is not the case of $[x]$ in Eq.~\eqref{eq:linear_system_stress_solve} (see, in this regards, the discussion in~\cite{coombs2011algorithmic}).
On top of this, solving the linear system~\eqref{eq:linear_system_stress_solve} in place of~\eqref{eq:linear_system_strain_solve} requires that the Legendre transform of $\bar{\rho} \ \psi$ is available in a closed form, which might not always be the case.

Most importantly, the inversion of matrix $[\check{A}]$ is necessary as the top-left entries of its inverse provide the relationship $\dfrac{\partial  {^t}\boldsymbol{\underline{\underline{\epsilon}}}^{e} }{\partial {^t}\boldsymbol{\underline{\underline{\epsilon}}}^{e, tr} } = [\check{A}^{^-1}]_{11}$ necessary for linearising $J^p$ as part of the linearisation process of the internal force vector, Eq.~\eqref{eq:derivative_plastic_Jacobian}.

\section{Particularisation to cylindrical coordinates of significant continuum quantities for axisymmetric problems}
\label{sec:cylindrical_coordinate_values}
To provide unambiguous matrix forms of quantities appearing throughout this manuscript and to facilitate the understanding of Algorithms~\ref{algorithm:workflow_UL} and~\ref{algorithm:workflow_TL}, the following matrices are explicitly presented:
\begin{itemize}
    \item Metric coefficients:
\begin{equation}
\label{eq:metric_coefficient_components}
\left[ {^{0}g_{AB}^{ax}} \right]
= \begin{bmatrix}
1 & 0 & 0
\\
0 & \left( {^{0}x^{I}} \right)^2 & 0
\\
0 & 0 & 1
\end{bmatrix},
\qquad
\left[ {^{t}g_{ab}^{ax}} \right]
= \begin{bmatrix}
1 & 0 & 0
\\
0 & \left( {^{t}x^{i}} \right)^2 & 0
\\
0 & 0 & 1
\end{bmatrix}
.
\end{equation}
\item 
Shifter:
\begin{equation}
\label{eq:axisymmetric_shifter}
\left[ \left( {^t_0}S^{ax} \right)^a_{\ A} \right]
=
\begin{bmatrix}
    1 & 0 & 0
    \\[8pt]
    0 & \ \dfrac{{^{0}x^{I}}}{{^{t}x^{i}}} & 0 
    \\[8pt]
    0 & 0 & 1
\end{bmatrix}.
\end{equation}
\item 
Non-zero components of the Christoffel symbols of the second kind:
\begin{equation}
\label{eq:Christoffel_symbol_original_non-zero_components}
\begin{dcases}
\left( {^0{\Gamma}}^{ax} \right)^{II}_{I \ II} 
=  
\left( {^0{\Gamma}}^{ax} \right)^{II}_{II \ I}
= \frac{1}{{^0{x}}^I};
\\
\left( {^0{\Gamma}}^{ax} \right)^{I}_{II \ II}  = - {^0{x}}^{I};
\end{dcases}
\qquad
\begin{dcases}
\left( {^t{\gamma}}^{ax} \right)^{ii}_{i \ ii}  
=  
\left( {^t{\gamma}}^{ax} \right)^{ii}_{ii \ i}  
= 
\frac{1}{{^t{x}}^i};
\\
\left( {^t{\gamma}}^{ax} \right)^{i}_{ii \ ii}
= - {^t{x}}^{i}.
\end{dcases}
\end{equation}
\item Total deformation gradient:
\begin{equation}
\left[ \left( {^t_0}X^{ax} \right)^{a}_{\ A} \right]
= 
\begin{bmatrix}
1 + \dfrac{\partial {^t{u}^i}}{\partial {^0{x}^I}} & 0 & \dfrac{\partial {^t{u}^i}}{\partial {^0{x}^{III}}}
\\[8pt]
0 & 1 & 0
\\[8pt]
\dfrac{\partial {^t{u}^{iii}}}{\partial {^0{x}^I}} & 0 & 1 + \dfrac{\partial {^t{u}^{iii}}}{\partial {^0{x}^{III}}}
\end{bmatrix}
\end{equation}
\item 
Transpose of the total deformation gradient:
\begin{equation}
\left[ \left( \left( {^t_0}X^{ax} \right)^T \right)^{A}_{\ \ a}  \right] = 
\begin{bmatrix}
1 + \dfrac{\partial {^t{u}^i}}{\partial {^0{x}^I}} 
& 0 
& \dfrac{\partial {^t{u}^{iii}}}{\partial {^0{x}^I}}
\\[8pt]
0 & \left(\dfrac{{^{t}x^{i}}}{{^{0}x^{I}}}\right)^2 & 0
\\[8pt]
\dfrac{\partial {^t{u}^i}}{\partial {^0{x}^{III}}}
& 0 
& 1 + \dfrac{\partial {^t{u}^{iii}}}{\partial {^0{x}^{III}}}
\end{bmatrix}
\end{equation}
\item Jacobian: 
\begin{equation}
{^t_0}J^{ax} = \det [ \left( {^t_0}X^{a}_{\ A} \right)^{ax}] \dfrac{{^t}x^i}{{^0}x^I}
\end{equation}
\item 
Gradient of test functions:
\begin{alignat}{2}
\underline{\textbf{grad}} \left( {^t}\underline{\boldsymbol{w}}^{ax} \right)
& =
{^t}{w}^a|_{b}
\ {^t}\underline{\boldsymbol{g}}_a \ {^t}\underline{\boldsymbol{g}}^b
& =
\begin{bmatrix}
\dfrac{\partial {\left( {^t}{w}^{ax} \right)^i}}{\partial {{^t}x^i}} & 0 &  \dfrac{\partial {\left( {^t}{w}^{ax} \right)^{i}}}{\partial {{^t}x^{iii}}} 
\\[8pt]
0 & \dfrac{{\left( {^t}{w}^{ax} \right)^i}}{{{^t}x^i}}& 0
\\[8pt]
\dfrac{\partial {\left( {^t}{w}^{ax} \right)^{iii}}}{\partial {{^t}x^i}} & 0 &  \dfrac{\partial {\left( {^t}{w}^{ax} \right)^{iii}}}{\partial {{^t}x^{iii}}}
\end{bmatrix}
{^t}\underline{\boldsymbol{g}}_a \ {^t}\underline{\boldsymbol{g}}^b;
\\
\underline{\textbf{Grad}} \left( {^t}\underline{\boldsymbol{W}}^{ax} \right)
& 
= {^t}{W}^a|_{B}
\ {^t}\underline{\boldsymbol{g}}_a \ {^0}\underline{\boldsymbol{g}}^B
&
=
\begin{bmatrix}
\dfrac{\partial {\left( {^t}{W}^{ax} \right)^i}}{\partial {{^0}x^I}} & 0 &  \dfrac{\partial {\left( {^t}{W}^{ax} \right)^{i}}}{\partial {{^0}x^{III}}} 
\\[8pt]
0 & \dfrac{{\left( {^t}{W}^{ax} \right)^i}}{{{^t}x^i}}& 0
\\[8pt]
\dfrac{\partial {\left( {^t}{W}^{ax} \right)^{iii}}}{\partial {{^0}x^I}} & 0 &  \dfrac{\partial {\left( {^t}{W}^{ax} \right)^{iii}}}{\partial {{^0}x^{III}}}
\end{bmatrix}
{^t}\underline{\boldsymbol{g}}_a \ {^0}\underline{\boldsymbol{g}}^B.
\end{alignat}
\end{itemize}

\pagebreak
\end{document}